\documentclass[10pt]{amsart}
\usepackage{amsmath}
\usepackage{amsfonts}
\usepackage{amsthm}
\usepackage{amssymb}
\usepackage{url}
\usepackage{verbatim}
\usepackage[usenames,dvips]{color}
\usepackage{graphicx}
\usepackage{xcolor}
\usepackage{mathrsfs}
\usepackage{lmodern}
\usepackage{bm}
\usepackage{mathtools}
\usepackage[english]{babel}

\numberwithin{equation}{section}
\numberwithin{subsection}{section}
\newtheorem{theorem}{Theorem}
\newtheorem{lemma}{Lemma}[section]
\newtheorem{corollary}[lemma]{Corollary}
\newtheorem{definition}[lemma]{Definition}
\newtheorem{remark}[lemma]{Remark}

\newtheorem{proposition}[lemma]{Proposition}

\newcommand{\di}{\mathrm{d}}

\newcommand{\bmzeta}{\bm{\zeta}}

\newcommand{\bmD}{\bm{D}}
\newcommand{\bme}{\bm{e}}

\newcommand{\bmF}{\bm{F}}
\newcommand{\bmg}{\bm{g}}

\newcommand{\bmh}{\bm{h}}

\newcommand{\bmk}{\bm{k}}

\newcommand{\bmN}{\bm{N}}

\newcommand{\bmu}{\bm{u}}

\newcommand{\bmx}{\bm{x'}}

\newcommand{\bmzero}{\bm{0}}

\DeclareMathOperator{\Div}{\mathrm{div}}
\DeclareMathOperator{\grad}{\mathrm{\nabla}}

\DeclareMathOperator{\pd}{\partial}

\usepackage[colorlinks=true,breaklinks]{hyperref}
\usepackage{breakurl}

\title[]{ Two-dimensional water waves with constant vorticity and general bottom topography }

\author{S. Pasquali$^{(\ast)}$}
\address[$\ast$]{International School for Advanced Studied (SISSA), via Bonomea 265, 34136, Trieste, Italy}
\email[$\ast$]{stefano.pasquali@sissa.it}

\begin{document}

\begin{abstract}
In this paper we consider two-dimensional water waves with constant vorticity, under the action of gravity and surface tension, in a fluid domain with finite depth and general bottom topography. We present a formulation which generalizes the one by Zakharov--Craig--Sulem for irrotational water waves, and the one by Constantin--Ivanov--Prodanov for water waves with constant vorticity and flat bottom topography. We study in detail an operator which appears in such formulation, extending well-known results for the classical Dirichlet--Neumann operator, such as an analyticity result, the Taylor expansion in homogeneous powers of the wave profile, and a paralinearization formula. As an application, we prove a local well-posedness result. 
\end{abstract}

\maketitle

\tableofcontents

\section{Introduction} \label{sec:intro}

This paper is concerned with the existence of two-dimensional gravity-capillary waves on the free surface of a perfect, incompressible and inviscid fluid with constant density and with \emph{constant vorticity} $\gamma$. We assume that the fluid fills a domain with \emph{finite depth and general bottom topography}, with space-periodic boundary conditions. \\

The problem of existence of two-dimensional water waves has a long history. The first existence results for two-dimensional travelling water waves have been proved in the irrotational setting by Stokes \cite{stokes1880theory}, and later for periodic waves by Levi-Civita \cite{levi1925determination}, and in the case of nonzero vorticity by Gerstner \cite{gerstner1809theorie}. Since then the problem has been extensively studied in the irrotational setting, see \cite{iooss2005standing,wu2009almost,alazard2015sobolev,ionescu2015global,baldi2018time,ai2022two,berti2022analyticity,feola2024quasi} and references therein for pure gravity water waves, and \cite{alazard2015gravity,berti2018amost} for gravity-capillary water waves. We also mention the monography \cite{lannes2013water} about irrotational water waves.

Also the study of two-dimensional water waves with nonzero vorticity has received a lot of attention in more recent years. For the case of water waves with constant vorticity we mention \cite{wahlen2007hamiltonian} regarding the Hamiltonian formulation of such a model, \cite{berti2021traveling,berti2024pure,barbieri2026bifurcation} regarding the existence of travelling waves, \cite{ifrim2018two} for a local well-posedness result for pure gravity waves and a long time existence result for small amplitude pure gravity waves, and \cite{berti2024hamiltonian} for an almost global existence result of small amplitude gravity-capillary waves. We also mention \cite{constantin2004exact,wahlen2006steady,ehrnstrom2022smooth,wahlen2023global} for the case of an arbitrary vorticity distribution. We point out that in the aforementioned works the authors assume that the fluid domain has either a flat bottom topography or infinite depth. 

The local well-posedness of the free-surface incompressible Euler equation without the irrotationality assumption has received more interest in recent years: regarding the gravity-capillary case, we mention \cite{coutand2007well,shatah2011local}; regarding the pure gravity case, we mention \cite{castro2015well} and the low-regularity results in \cite{wang2021local,ifrim2025sharp}.\\

Allowing for general (and possibly time-dependent) bottom topography in fluid models, as described by Lannes \cite{lannes2013water} in the irrotational case, could be useful in order to describe some natural phenomena in a more precise way; we mention for example \cite{iguchi2011mathematical}, where the author considers a fluid domain with time-dependent bottom topography in order to describe tsunamis generated by the deformation of the seabed caused by submarine earthquakes.

In the irrotational setting, we mention \cite{alazard2011water}, where the authors prove a local well-posedness result for gravity-capillary water waves (see also \cite{alazard2014cauchy} for pure gravity waves).  One remarkable aspect of their result is that no regularity assumption on the time-independent bottom of the fluid domain is assumed (if the bottom is time-dependent, their result is valid under a Lipschitz regularity assumption on the bottom).\\

In this paper we study two-dimensional gravity-capillary waves for an incompressible, inviscid and homogeneous fluid and with constant vorticity $\gamma \in \mathbb{R}$, assuming that the fluid domain has a time-independent bottom with finite depth and general topography, with space-periodic boundary conditions. Since the vorticity is constant, the velocity field decomposes as the sum of a Couette flow and an irrotational velocity field $\varphi$ solving the elliptic boundary value problem \eqref{eq:EllBVP}. 

Moreover, denoting by $\eta$ the function describing the free surface of the wave and by $\beta$ the function describing the displacement of the bottom of the fluid domain with respect to the flat bottom case, we introduce an operator $G(\eta,\beta,\gamma)$ which reduces to the classical Dirichlet--Neumann operator in the irrotational case $\gamma=0$, and which appears naturally in the system \eqref{eq:WWsys} describing the time evolution of the fluid; such formulation is based on an approach developed by Iguchi in \cite{iguchi2011mathematical} for dealing with elliptic boundary value problems with nonhomogeneous Neumann boundary conditions at the bottom. We study in detail the operator $G(\eta,\beta,\gamma)$, proving an analyticity result in Sec. \ref{subsec:analOpG} and a paralinearization result in Sec. \ref{subsec:ParaOpG}; our result generalize the ones obtained for the classical Dirichlet--Neumann operator in \cite{craig1993numerical,craig2005hamiltonian,alazard2011water,lannes2013water}.

Next, by following the approach of \cite{alazard2011water}, we use the above results on the operator $G(\eta,\beta,\gamma)$ in order to deduce a local well-posedness result. Compared to \cite{coutand2007well,shatah2011local}, our result allows for less regular initial data.

\subsection{The model} \label{subsec:model}

We consider two-dimensional gravity-capillary waves on the free surface of a perfect, incompressible and inviscid fluid with constant density equal to $1$, and with \emph{constant vorticity} $\gamma$. We assume that the fluid fills a domain with finite depth and general bottom topography, with space-periodic boundary conditions; writing $\mathbb{T} := \mathbb{R}/(2\pi\mathbb{Z})$, we denote the fluid domain in the following way,
\begin{align} \label{eq:FluidDom}
D_{\eta,\beta}(t) &:= \{ (x,y) \in \mathbb{T} \times \mathbb{R} : -h + \beta(x) < y < \eta(t,x) \},
\end{align}
where $\beta: \mathbb{T} \to \mathbb{R}$ denotes the variation with respect to the flat bottom of finite depth $h>0$, and where the free surface is given by the graph of an unknown function $\eta:\mathbb{R} \times \mathbb{T} \to \mathbb{R}$. We also denote by $D_0$ the fixed domain $\mathbb{T} \times (-h,0)$. \\

The equations describing the flow are given by
\begin{align}
\Div \bmu &= 0, \; \; \text{in} \; \; D_{\eta,\beta}(t), \label{eq:Incompr} \\
\bmu_t + (\bmu \cdot \nabla)  \bmu &= - \nabla \mathscr{P} - g \bme_2, \; \; \text{in} \; \; D_{\eta,\beta}(t), \label{eq:Euler} \\
\partial_x u_{2} - \partial_y u_1 &= \gamma, \; \; \text{in} \; \; D_{\eta,\beta}(t), \label{eq:ConstVort} \\
\bmu \cdot \bmN_{b} &= 0 , \; \; \text{at} \; \; y=-h+\beta, \label{eq:Imperm} \\
\eta_t - \bmu \cdot \bmN &= 0, \; \; \text{at} \; \; y=\eta, \label{eq:KinFree} \\
\mathscr{P} - \mathscr{P}_0 &=  - \kappa \, \left(  \frac{\eta_x}{ \sqrt{1+\eta_x^2}} \right)_x , \; \; \text{at} \; \; y=\eta, \label{eq:DynBC}
\end{align}
where $\bmu = 
\begin{pmatrix}
u_1 \\ u_2
\end{pmatrix}
: \overline{D_{\eta,\beta}(t)} \to \mathbb{R}^2$ is the velocity field, $\bme_2=(0,1)^T$, $\bmN=(-\eta_x,1)^T$ denotes the outward normal vector at the free surface, $\bmN_{b}=(\beta_x,-1)^T$ denotes the outward normal vector at the bottom, $g>0$ is the gravity constant, $\gamma \in \mathbb{R}$ is the value of the vorticity (assumed to be constant), $\mathscr{P}_0$ is the constant atmospheric pressure and $\kappa>0$ is the surface tension coefficient.  The divergence operator used above is defined as $\Div \bmu := \nabla \cdot \bmu$. 

Whenever it is not important, we omit the time-dependence of $D_{\eta,\beta}(t)$ and we will denote it by $D_{\eta,\beta}$. 

Due to the constant value of the vorticity, the velocity field $\bmu$ is given by the sum of the Couette flow 
$ \begin{pmatrix}
-\gamma y \\ 0
\end{pmatrix} $
and an irrotational velocity field, that can be expressed as the gradient of a function $\varphi$, called velocity potential. Denoting by
\begin{align*}
\psi(t,x) &:= \varphi(t,x,\eta(t,x))
\end{align*}
the evaluation of the velocity potential at the free surface, we can recover $\varphi$ by solving the following elliptic boundary value problem,
\begin{equation} \label{eq:EllBVP}
\begin{cases}
\Delta \varphi = 0 \; \; &\text{in} \; \; D_{\eta,\beta}, \\
\varphi = \psi \; \; &\text{at} \; \; y=\eta, \\
\gamma \, y \, \beta_x - \varphi_x \, \beta_x + \varphi_y = 0  \; \; &\text{at} \; \; y=-h+\beta .
\end{cases}
\end{equation}

For $s \in \mathbb{R}$, we define the Sobolev space 
\begin{align*}
H^s(\mathbb{T}) &\coloneqq \left\{ f(x) = \sum_{j \in \mathbb{Z}} f_j e^{\mathrm{i} jx} : \|f\|_{H^s(\mathbb{T})} \coloneqq \left( \sum_{j \in \mathbb{Z}} |f_j|^2 \langle j \rangle^{2s} \right)^{1/2} < +\infty \right\}, \\
\langle j \rangle &\coloneqq (1+j^2)^{1/2} ,
\end{align*}
and the corresponding Sobolev space of functions with zero average
\begin{align*}
H^s_0(\mathbb{T}) &:= \left\{ f \in H^s(\mathbb{T}) : \int_{\mathbb{T}} f(x) \di x = 0 \right\},
\end{align*}
with the associated norm
\begin{align*}
\|f\|_{H^s_0(\mathbb{T})} &\coloneqq \left( \sum_{j \in \mathbb{Z}\setminus \{0\} } |f_j|^2 \langle j \rangle^{2s} \right)^{1/2} , \;\; \forall f \in H^s_0(\mathbb{T}),
\end{align*}
while we denote by $\dot{H}^s(\mathbb{T})$ the homogeneous Sobolev space of order $s$, namely the quotient space $H^s(\mathbb{T})/\mathbb{R}$; for simplicity we denote the equivalent classes $[f] = \{ f+ c, c \in \mathbb{R} \}$ just by $f$.

In order to study the boundary value problem \eqref{eq:EllBVP}, we follow the approach developed by Iguchi (see Sec. 3 in \cite{iguchi2011mathematical}; see also Appendix A.5 in \cite{lannes2013water}), introducing a suitable family of operators. Given $s \in \mathbb{R}$ and $\eta, \beta , \psi \in H^s(\mathbb{T})$, we define the \emph{Dirichlet--Neumann operator} $G^{DN}(\eta,\beta)$, the \emph{Neumann--Neumann operator} $G^{NN}(\eta,\beta)$, the \emph{Dirichlet--Dirichlet operator} $G^{DD}(\eta,\beta)$ and the \emph{Neumann--Dirichlet operator} $G^{ND}(\eta,\beta)$ such that
\begin{align}
G^{DN}(\eta,\beta)\psi + \gamma \, G^{NN}(\eta,\beta) \left( (-h+\beta)\beta_x \right) &\coloneqq \nabla\varphi \cdot \bmN |_{y=\eta(x)}, \label{eq:DNO} \\
G^{DD}(\eta,\beta)\psi + \gamma \, G^{ND}(\eta,\beta) \left( (-h+\beta)\beta_x \right) &\coloneqq \varphi |_{y=-h+\beta(x)}, \label{eq:DDO}
\end{align}
where $\varphi$ is the solution of the boundary value problem \eqref{eq:EllBVP}.
In the irrotational case $\gamma=0$ the Dirichlet--Neumann operator $G^{DN}(\eta,\beta)$ has been studied in detail, see Sec.2 and Appendix A in \cite{craig2005hamiltonian}. We defer the study of  $G^{NN}(\eta,\beta)$, $G^{DD}(\eta,\beta)$ and $G^{ND}(\eta,\beta)$ to Appendix \ref{sec:BVP}. To simplify the notation, we also introduce the operator $G(\eta,\beta,\gamma)$ by
\begin{align} \label{eq:OpG}
G(\eta,\beta,\gamma)(\psi) &\coloneqq G^{DN}(\eta,\beta)\psi + \gamma \, G^{NN}(\eta,\beta) \left( (-h+\beta)\beta_x \right) .
\end{align}
We defer the study of the operator $G(\eta,\beta,\gamma)$ to Sec. \ref{sec:OpG}.

\subsection{Formulation for the flat bottom case} \label{subsec:HamFlat}

We now present a formulation introduced by Constantin et al. in \cite{constantin2007nearly} for the system \eqref{eq:Incompr}-\eqref{eq:DynBC} in the case $\beta \equiv 0$. In this case the time evolution of the fluid can be described by 
\begin{equation} \label{eq:WWsysF}
\begin{cases}
\eta_t &= G^{DN}(\eta)\psi + \gamma \, \eta \, \eta_x \\
\psi_t &= - g \, \eta - \frac{1}{2} \, \psi_x^2 + \frac{ ( G^{DN}(\eta)\psi + \eta_x \psi_x )^2 }{ 2(1+\eta_x^2)} + \kappa \, \left(  \frac{\eta_x}{ \sqrt{1+\eta_x^2}} \right)_x  \\
&\qquad + \gamma \, \left( \eta \, \psi_x +  \partial_x^{-1}G^{DN}(\eta)\psi \right)  ,
\end{cases}
\end{equation}
where $G^{DN}(\eta)$ is the classical \emph{Dirichlet--Neumann operator}
\begin{align} \label{eq:DNOF}
G^{DN}(\eta)\psi &:= \nabla\varphi \cdot \bmN |_{y=\eta(x)}.
\end{align}

\begin{remark} \label{rem:ZCS}

For $\gamma=0$ the system \eqref{eq:WWsysF} reduces to the classical Zakharov--Craig--Sulem formulation for irrotational fluids, see \cite{craig1993numerical}.

\end{remark}

Since the bottom of the fluid domain is flat, the system \eqref{eq:WWsysF} is invariant by space translation. Notice that since $G^{DN}(\eta)\psi$ has zero average, then the quantity $\int_{\mathbb{T}} \eta(x) \mathrm{d}x$ is a constant of motion of \eqref{eq:WWsysF}; moreover, since $G^{DN}(\eta)1 =0$, then the vector field on the right-hand side of \eqref{eq:DNOF} depends only on $\eta$ and on $\psi - \int_{\mathbb{T}} \psi \frac{ \di x }{2\pi}$. 

We assume that the variables $(\eta,\psi)$ in \eqref{eq:WWsysF} belong to the space $H^1_0(\mathbb{T}) \times \dot{H}^1(\mathbb{T})$.

The system \eqref{eq:WWsysF} admits also the following Hamiltonian structure: let us consider the space $H^1_0(\mathbb{T}) \times \dot{H}^1(\mathbb{T})$, endowed with the non canonical Poisson tensor
\begin{align} \label{eq:NonCPoisson}
J_{\gamma} &:=
\begin{pmatrix}
0 & \mathrm{Id} \\ - \mathrm{Id} & \gamma \, \partial_x^{-1}
\end{pmatrix}
.
\end{align}
The non canonical Poisson tensor $J_{\gamma}$ is well-defined as an operator from (a subspace of)  $( L^2_0(\mathbb{T}) \times \dot{L}^2(\mathbb{T}) )^{\ast}$ to $L^2_0(\mathbb{T}) \times \dot{L}^2(\mathbb{T})$, since $\partial_x^{-1}$ maps a dense subspace of  $L^2_0(\mathbb{T})$ to $\dot{L}^2(\mathbb{T})$; for simplicity, throughout the paper we omit this detail, and we identify with $\dot{L}^2(\mathbb{T}) \times L^2_0(\mathbb{T})$ the dual space $( L^2_0(\mathbb{T}) \times \dot{L}^2(\mathbb{T}) )^{\ast}$ with respect to the scalar product in $L^2(\mathbb{T})$. Let us consider the Hamiltonian
\begin{align} 
& H_{0,\gamma}(\eta,\psi) \nonumber \\
&= \frac{1}{2} \int_{\mathbb{T}} \psi \, G^{DN}(\eta)\psi + g \, \eta^2  + \gamma \, \left( -\psi_x \, \eta^2 + \frac{\gamma}{3} \, \eta^3 \right) \di x + \kappa \, \int_{\mathbb{T}} \sqrt{1+\eta_x^2} \, \di x , \label{eq:HamF}
\end{align}
which is well defined on $H^1_0(\mathbb{T}) \times \dot{H}^1(\mathbb{T})$; then \eqref{eq:WWsysF} are the Hamilton's equations associated to \eqref{eq:HamF} with respect to the Poisson tensor \eqref{eq:NonCPoisson}, namely
\begin{align} \label{eq:HamEqNC}
\pd_t \, 
\begin{pmatrix}
\eta \\ \psi
\end{pmatrix}
&=  J_{\gamma}
\begin{pmatrix}
\nabla_{\eta} H_{0,\gamma} \\ \nabla_{\psi} H_{0,\gamma}
\end{pmatrix}
,
\end{align}
where $(\nabla_{\eta} H_{0,\gamma} , \nabla_{\psi} H_{0,\gamma} ) \in \dot{L}^2(\mathbb{T}) \times L^2_0(\mathbb{T})$ denote the $L^2$-gradients.

The system \eqref{eq:WWsysF} is also time reversible: if we define on the space $H^1_0(\mathbb{T}) \times \dot{H}^1(\mathbb{T})$ the involution
\begin{align} \label{eq:revers}
\mathscr{S}
\begin{pmatrix}
\eta \\ \psi
\end{pmatrix}
(x) &\coloneqq
\begin{pmatrix}
\eta(-x) \\ -\psi(-x)
\end{pmatrix}
,
\end{align}
then $H_{0,\gamma}$ is invariant under $\mathcal{S}$, namely $H_{0,\gamma} \circ \mathscr{S} = H_{0,\gamma}$ (see Sec. 2 of \cite{berti2021traveling}). 

We mention that Wahl\'en introduced another Hamiltonian formulation, in which the Poisson tensor takes the canonical form (see \cite{wahlen2007hamiltonian} and Sec. 2 of \cite{berti2021traveling}).

\subsection{Formulation for general bottom topography} \label{subsec:VarDepth}

We now study the system \eqref{eq:Incompr}-\eqref{eq:DynBC} for general bottom topography. 

First, we observe that \eqref{eq:KinFree} can be rewritten as
\begin{align} \label{eq:KinFree2}
\eta_t + \eta_x (\varphi_x - \gamma \, y) - \varphi_y &= 0, \;\; \text{at}\;\; y = \eta .
\end{align}
We also notice that the incompressibility equation \eqref{eq:Incompr} implies the existence of a stream function $\Psi: D_{\eta,\beta}(t) \to \mathbb{R}$ such that
\begin{align*}
\Psi_x = - u_2 , &\;\; \Psi_y = u_1, \;\; \text{in} \;\; D_{\eta,\beta}(t) ,
\end{align*}
so that the kinematic boundary condition takes the form
\begin{align*} 
\eta_t + ( \Psi|_{y=\eta} )_x  &= 0, \;\; \text{at}\;\; y = \eta .
\end{align*}
Moreover, \eqref{eq:Euler} takes the form
\begin{align*}
\nabla \, \left[ \varphi_t + \frac{1}{2} |\nabla\Psi|^2 + \gamma \, \Psi + g \, y + \mathscr{P} \right] &= \bmzero,  \;\; \text{in} \;\; D_{\eta,\beta}(t) ,
\end{align*}
so that by evaluating the above equation at the surface and using that the velocity potential $\varphi$ is defined up to an arbitrary function of time, we obtain that
\begin{align} \label{eq:DynBC2}
 \varphi_t + \frac{1}{2} |\nabla\Psi|^2 + \gamma \, \Psi + g \, y - \kappa \, \left(  \frac{\eta_x}{ \sqrt{1+\eta_x^2}} \right)_x  &= 0,  \;\; \text{at} \;\; y = \eta ,
\end{align}
which is equivalent to the dynamic boundary condition \eqref{eq:DynBC}. Therefore, we obtain that $(\eta,\varphi,\Psi)$ satisfy
\begin{equation*} 
\begin{cases}
\Delta \varphi = 0 \; \; &\text{in} \; \; D_{\eta,\beta}, \\
\varphi = \psi \; \; &\text{at} \; \; y=\eta, \\
\gamma \, y \, \beta_x - \varphi_x \, \beta_x + \varphi_y = 0  \; \; &\text{at} \; \; y=-h+\beta , \\
\eta_t + \eta_x (\varphi_x - \gamma \, y) - \varphi_y = 0, \;\; &\text{at}\;\; y = \eta \\
\varphi_t + \frac{1}{2} |\nabla\Psi|^2 + \gamma \, \Psi + g \, y  - \kappa \, \left(  \frac{\eta_x}{ \sqrt{1+\eta_x^2}} \right)_x  = 0,  \;\; &\text{at} \;\; y = \eta .
\end{cases}
\end{equation*}

Using the boundary value problem \eqref{eq:EllBVP} and \eqref{eq:DNO}, together with the equality
\begin{align*}
\varphi_y |_{y=\eta} &= \frac{ G(\eta,\beta,\gamma)(\psi)  + \psi_x \eta_x  }{1+\eta_x^2} ,
\end{align*}
we have that
\begin{align*}
\frac{1}{2} |\nabla\Psi|^2 + \gamma \, \Psi |_{y=\eta} 
&= \frac{1}{2} |\nabla\varphi|^2 -\gamma \eta \varphi_x + \gamma \, \left( \Psi + \frac{1}{2} \gamma \eta^2 \right) |_{y=\eta}
\end{align*}
where 
\begin{align*}
\Psi + \frac{1}{2} \gamma \eta^2 |_{y=\eta} &= - \partial_x^{-1} G(\eta,\beta,\gamma)(\psi) ,
\end{align*}
and that
\begin{align*}
\psi_t &= \varphi_t + \eta_t \varphi_y |_{y=\eta} ,\\
\psi_x &= \varphi_x + \eta_x \varphi_y |_{y=\eta} ,\\
-\gamma \, \eta \, \varphi_x |_{y=\eta} &= -\gamma \, \eta \, \psi_x  + \gamma \, \eta  \eta_x \, \varphi_y |_{y=\eta} , \\
|\nabla\varphi|^2 |_{y=\eta} &= \psi_x^2 + (1+\eta_x^2) \varphi_y^2 |_{y=\eta} - 2 \eta_x \psi_x \, \varphi_y |_{y=\eta} ,
\end{align*}
so that
\begin{align*}
\varphi_t +\frac{1}{2} |\nabla\varphi|^2 - \gamma \eta \varphi_x |_{y=\eta} 
&= \psi_t - \left( G(\eta,\beta,\gamma)(\psi)  \right)\, \varphi_y |_{y=\eta} - \gamma \, \eta \eta_x \, \varphi_y |_{y=\eta} \\
&\quad +\frac{1}{2} \psi_x^2 + \frac{1}{2(1+\eta_x^2)} \left[ G(\eta,\beta,\gamma)(\psi)  + \psi_x \eta_x \right]^2 \\
&\quad - \eta_x \psi_x  \, \varphi_y |_{y=\eta} -\gamma \eta \psi_x + \gamma \, \eta \eta_x \varphi_y |_{y=\eta} .
\end{align*}

Therefore the time evolution of the fluid can be described by the following system,
\begin{equation} \label{eq:WWsys}
\begin{cases}
\eta_t &= G(\eta,\beta,\gamma)(\psi)  + \gamma \, \eta \, \eta_x \\
\psi_t &= - \frac{1}{2} \, \psi_x^2 + \frac{ 1 }{ 2(1+\eta_x^2)}   \left[ G(\eta,\beta,\gamma)(\psi)  + \eta_x \psi_x \right]^2 \\
&\quad + \gamma \, \left[ \eta \, \psi_x +  \partial_x^{-1} G(\eta,\beta,\gamma)(\psi)  \right] - g \, \eta + \kappa \, \left(  \frac{\eta_x}{ \sqrt{1+\eta_x^2}} \right)_x .
\end{cases}
\end{equation}

We briefly discuss the phase space of \eqref{eq:WWsys}. Recalling the boundary value problem \eqref{eq:EllBVP} (more precisely, the interior equation and the condition at the bottom) and using the divergence theorem, we have that also $G(\eta,\beta,\gamma)(\psi)$ has zero average, which  implies that $\int_{\mathbb{T}} \eta(x) \mathrm{d}x$ is a constant of motion for \eqref{eq:WWsys}. Hence, we assume that the variables $(\eta,\psi)$ in \eqref{eq:WWsys} belong to the space $H^1_0(\mathbb{T}) \times H^1(\mathbb{T})$.

Finally, if we denote by $D_{\eta,\beta}$ the domain
\begin{align*} 
D_{\eta,\beta} &= \{ (x,y) \in \mathbb{T} \times \mathbb{R}: -h+\beta(x) < y < \eta(x) \} ,
\end{align*}
we say that $D_{\eta,\beta}$ is \emph{strictly connected} if there exists $h_0>0$ such that
\begin{align} \label{eq:StrConnected}
h - \beta(x) + \eta(x) &\geq h_0, \; \; \forall \; x \in \mathbb{T} .
\end{align}

\section{Main results} \label{sec:main}

\subsection{Properties of the operator $G(\eta,\beta,\gamma)$} \label{subsec:OpGmain}

We first state an analyticity result for $G(\eta,\beta,\gamma)$, which extends the corresponding analyticity result for the classical Dirichlet--Neumann operator in the irrotational case (see Theorem A.11 in \cite{lannes2013water}).

\begin{theorem} \label{thm:AnalOpGMain}

Let $s > 3/2$, and fix $\psi \in H^{s}(\mathbb{T})$, $\gamma \in \mathbb{R}$. If we consider the operator $G(\eta,\beta,\gamma)$ defined in \eqref{eq:OpG}, then the map
\begin{align*}
G(\cdot,\cdot,\gamma)(\psi) &: \{ (\eta,\beta) \in H^{s}(\mathbb{T}) \times H^{s}(\mathbb{T}) : \eqref{eq:StrConnected} \; \text{holds true} \} \to H^{s-1}(\mathbb{T})  
\end{align*}
is analytic.
\end{theorem}

We defer the proof of the above result to Sec. \ref{subsec:analOpG}.

We mention that in Sec.\ref{subsec:analOpG} we also prove a more explicit analyticity result (see Corollary \ref{cor:AnalOpGsmall}) for $(\eta,\beta)$ belonging to a neighbourhood of the origin in a suitable Sobolev space; such result allows us to study the homogeneous expansion of the operator $G(\eta,\beta,\gamma)$ around $\eta=0$ 
\begin{align*}
G(\eta,\beta,\gamma) &= \sum_{j=0}^{\infty} G_j[\eta](\beta,\gamma) ,
\end{align*}
where $G_j[\eta](\beta,\gamma)$ is homogeneous of degree $j$ in $\eta$. In particular, we prove that
\begin{align}
G_{0}(\beta,\gamma)(\psi) &= ( D \, \tanh(h D) + D \, L(\beta) ) \psi + \gamma \, \nu(\beta) , \label{eq:G0newMain} \\
D \, L(\beta)  \psi &\coloneqq   \mathscr{F}^{-1} \, \left\{ \left[  - \xi \, \frac{ D\beta }{ \cosh(h \, \xi) \, ( \cosh(h \, \xi) - (D\beta) \, \sinh(h \, \xi) ) } \right] \, (\mathscr{F} \psi) \, \right\} , \nonumber \\
\nu(\beta) &\coloneqq - \mathscr{F}^{-1} \, \bigg[   \frac{ 1 }{  \cosh(h \, \xi) - \, (D\beta) \, \sinh(h \, \xi) } \, \mathscr{F} \left( (-h+\beta) \beta_x \right) \bigg] . \nonumber
\end{align}

For more details see Sec. \ref{subsec:HomExp}. Our expansion reduces to the one computed in \cite{craig2005hamiltonian} for the irrotational case $\gamma=0$. \\

Finally, in Sec. \ref{subsec:ParaOpG} we prove a paralinearization formula for the operator $G(\eta,\beta,\gamma)$. The paradifferential calculus was introduced by Bony in order to deal with the quantization of symbols $a(x,\xi)$ of degree $m$ with respect to $\xi$ and limited regularity in the $x$-variable, to which are associated operators of order $m$ denoted by $T_a$ (we defer to Appendix \ref{sec:Paradiff} for recalling some standard rules of paradifferential calculus).

For any $u \in \mathcal{S}'(\mathbb{R})$ we denote its Fourier transform by $\mathscr{F}u$. If $p: \mathcal{S}(\mathbb{R}) \to \mathcal{S}'(\mathbb{R})$, then 
\begin{align*}
p(D)u &:= \mathscr{F}^{-1}(p(\xi) \; \mathscr{F}u).
\end{align*}
In the following we write $\mathbb{R}_{\neq} = \mathbb{R}\setminus\{0\}$. The above definitions apply also to periodic distributions.

\begin{definition} \label{def:symbolsPara}
Let $\varrho \geq 0$ and $m \in \mathbb{R}$, we denote by $\Gamma^m_\varrho(\mathbb{T})$ the space of functions $a(x,\xi)$ on $\mathbb{T} \times \mathbb{R}_{\neq}$ which are of class $C^\infty$ with respect to $\xi \neq 0$ and such that for all $\alpha \in \mathbb{N}$ and all $\xi \neq 0$ the map $x \mapsto \pd_{\xi}^\alpha a(x,\xi)$ belongs to $W^{\varrho,\infty}(\mathbb{T})$ and there exists $C_\alpha >0$ such that
\begin{align}
\|\pd_{\xi}^\alpha a(\cdot,\xi)\|_{ W^{\varrho,\infty}(\mathbb{T}) } &\leq C_\alpha (1 + |\xi|)^{m-|\alpha|}. \label{eq:symbolsEst}
\end{align}
\end{definition}

We now introduce the space of pluri-homogeneous symbols.

\begin{definition} \label{def:symbolsHomPara}
Let $\varrho \geq 1$ and $m \in \mathbb{R}$, we denote by $\Sigma^m_\varrho(\mathbb{T})$ the space of symbols of the form
\begin{align*}
a(x,\xi) &= \sum_{0 \leq j < \varrho} a_{m-j}(x,\xi),
\end{align*}
where $a_{m-j} \in \Gamma^{m-j}_{\varrho-j}(\mathbb{T})$ is homogeneous of degree $m-j$ in the variable $\xi$ and is of class $C^\infty$ in the variable $\xi \in \mathbb{R}_{\neq}$, and with regularity $W^{\varrho -j,\infty}$ in the variable $x$. We say that $a_m$ is the principal symbol of $a$.
\end{definition}

Let $m \in \mathbb{R}$, $\varrho \in [0,1]$, $a \in \Gamma^m_{\varrho}(\mathbb{T})$, we also define the semi-norm
\begin{align} \label{eq:NormPara}
M^m_{\varrho}(a) &:= \sup_{|\alpha| \leq 6+\varrho} \sup_{\xi \in \mathbb{R}_{\neq}} \| (1+|\xi|)^{\alpha-|m|} \, a(\cdot, \xi) \|_{W^{\varrho,\infty}(\mathbb{T})} .
\end{align}

In order to define rigorously paradifferential operators we first introduce a fixed cutoff functions $\chi$. Fix sufficiently small $\varepsilon_1$, $\varepsilon_2$ such that $0 <  \varepsilon_1 < \varepsilon_2 $ and a function $\chi \in C^\infty$ homogeneous function of degree $0$ such that:
\begin{itemize}
\item[i.] the following conditions are satisfied,
\begin{align*}
\chi(\xi_1 , \xi_2 ) &= 1, \; \; \text{if} \; \; |\xi_1 | \leq \varepsilon_1 |\xi_2 | , \\
\chi(\xi_1 , \xi_2 ) &= 0, \; \; \text{if} \; \; |\xi_1 | \geq \varepsilon_2 |\xi_2 | ;
\end{align*}
\item[ii] the following symmetry property holds true,
\begin{align} \label{eq:chiSymm}
\chi(\xi_1,\xi_2) &= \chi(-\xi_1,-\xi_2) = \chi(-\xi_1,\xi_2) .
\end{align}
\end{itemize}

Given a symbol $a \in \Sigma^m_\varrho(\mathbb{T})$, we define the paradifferential operator $T_a$ by
\begin{align}
\mathscr{F}(T_au)(\xi) &:= (2\pi)^{-1} \sum_{\eta \in \mathbb{Z}} \chi(\xi-\eta,\eta) \; \; \mathscr{F}a(\xi-\eta,\eta) \; \; \mathscr{F}u(\eta) , \label{eq:ParadiffOp}
\end{align}
where $\mathscr{F}a$ is the Fourier transform of $a$ with respect to the first variable.  \\

Now we state the paralinearization result.

\begin{theorem} \label{thm:OpGparaMain}

Let $s > 5/2$. If $\eta, \beta \in H^{s+1/2}(\mathbb{T})$ are such that \eqref{eq:StrConnected} holds true, if $\gamma \in \mathbb{R}$ and if $\psi \in H^s(\mathbb{T}) $, then  
\begin{align} \label{eq:OpGparaMain}
G(\eta,\beta,\gamma)(\psi) &= T_{\lambda} \omega  - T_{ V(\eta,\beta,\gamma)(\psi) } \eta_x + \mathcal{R}(\eta,\beta,\gamma,\psi) ,
\end{align}
where 
\begin{align*}
\omega(\eta,\beta,\gamma)(\psi) &= \psi - T_{B(\eta,\beta,\gamma)(\psi)} \eta \in H^s(\mathbb{T}) , \\
V(\eta,\beta,\gamma)(\psi) &= \psi_x - B(\eta,\beta,\gamma)(\psi) \, \eta_x \in H^{s-1}(\mathbb{T}), \\
B(\eta,\beta,\gamma)(\psi) &= \frac{ G(\eta,\beta,\gamma)(\psi) + \eta_x \psi_x }{1+\eta_x^2} \in H^{s-1}(\mathbb{T}) , 
\end{align*}
and where the symbol $\lambda = \lambda^{(1)} + \lambda^{(0)}$ does not depend on $\gamma$ and on $\beta$, with $\lambda^{(1)} =| \xi |$ and $\lambda^{(0)}$ being a symbol of order $0$.

Moreover,the remainder $\mathcal{R}(\eta,\beta,\gamma,\psi) \in H^{s+1/2}(\mathbb{T})$ satisfies
\begin{align*}
& \| \mathcal{R}(\eta,\beta,\gamma,\psi) \|_{ H^{s+1/2}(\mathbb{T}) } \\
&\leq C_1 \left( \| \eta \|_{ H^{s+1/2}(\mathbb{T}) } ,\|\beta\|_{ H^{s+1/2}(\mathbb{T}) }  \right) \, \| \psi \|_{ H^{s}(\mathbb{T}) } + |\gamma| \, C_2 \left( \| \eta \|_{ H^{s+1/2}(\mathbb{T}) } , \|\beta\|_{ H^{s+1/2}(\mathbb{T}) }  \right)  ,
\end{align*}
for some non-decreasing functions $C_1, C_2$ depending on $h$ and $h_0$.
\end{theorem}

As mentioned in the statement of the above result, the principal symbol of the symbol $\lambda$ is given by $\lambda^{(1)} = |\xi|$, which is the principal symbol of the classical Dirichlet--Neumann operator. In the case of a fluid domain with finite depth, our formula \eqref{eq:OpGparaMain} generalizes the corresponding one proved in Proposition 3.14 of \cite{alazard2011water} for the classical Dirichlet--Neumann operator.

\subsection{Local well-posedness} \label{subsec:LWP}

We consider \eqref{eq:WWsys}, namely
\begin{equation*} 
\begin{cases}
\eta_t &= G(\eta,\beta,\gamma)(\psi)  + \gamma \, \eta \, \eta_x \\
\psi_t &= -g \, \eta - \frac{1}{2} \, \psi_x^2 + \frac{ 1 }{ 2(1+\eta_x^2)}   \left[ G(\eta,\beta,\gamma)(\psi)  + \eta_x \psi_x \right]^2 \\
&\quad + \gamma \, \left[ \eta \, \psi_x +  \partial_x^{-1} G(\eta,\beta,\gamma)(\psi) \right]  + \kappa \, \left(  \frac{\eta_x}{ \sqrt{1+\eta_x^2}} \right)_x .
\end{cases}
\end{equation*}
with initial data
\begin{align} \label{eq:InData}
\eta|_{t=0} = \eta_{0} , &\;\; \psi|_{t=0} = \psi_{0} .
\end{align}
We recall that the average of $\eta$ is conserved in time, and we set it equal to $0$ without loss of generality. We recall that we denote by $H^s_0(\mathbb{T})$ the Sobolev spaces of functions with zero average.

\begin{theorem} \label{thm:LWPMain}

Let $\gamma \in \mathbb{R}$, let $h,\kappa>0$, and let $s_0 \coloneqq \frac{5}{2}$. Then for any $s > s_0$, for any $(\eta_{0},\psi_{0}) \in H^{s+1/2}_0(\mathbb{T}) \times H^{s}(\mathbb{T})$ and for any $\beta \in H^{s+1/2}(\mathbb{T})$ such that 
\begin{itemize}
\item[(A1)] the initial fluid domain is strictly connected, namely there exists $h_0 >0$ such that
\begin{align*} 
h - \beta(x) + \eta_{0}(x) &\geq h_0, \;\; \forall x \in \mathbb{T} ;
\end{align*}
\end{itemize}
then there exists $T>0$ such that the Cauchy problem \eqref{eq:WWsys}-\eqref{eq:InData} has a unique solution
\begin{align*}
( \eta , \psi ) \in C^0([0,T]; H^{s+1/2}_0(\mathbb{T}) \times H^s(\mathbb{T}) ) .
\end{align*}

\end{theorem}

Some comments regarding Theorem \ref{thm:LWPMain} are in order.

\begin{remark} \label{rem:RegLWP}

The proof of Theorem \ref{thm:LWPMain} follows the strategy of \cite{alazard2011water}; compared to the local well-posedness results in \cite{coutand2007well,shatah2011local}, our result allows for less regular initial data. In the non-periodic setting, we expect that the exponent $s_0$ could be lowered (for such a result in the irrotational case, see \cite{alazard2014cauchy}).

\end{remark}

\begin{remark} \label{rem:comparLWP}

Compared to usual results in local well-posedness for water waves (e.g. Theorem 9.6 in \cite{lannes2013water} and Theorem 1.1 in \cite{berti2024hamiltonian}), we assumed $\psi \in H^s(\mathbb{T})$ instead of $\psi \in \dot{H}^s(\mathbb{T})$. If we assume that $\gamma \beta_x \equiv 0$ and we write $\langle \psi \rangle \coloneqq (2\pi)^{-1} \int_{\mathbb{T}} \psi(x) \mathrm{d}x$, we see from \eqref{eq:WWsys} that the evolution of $(\eta,\psi- \langle \psi \rangle )$ is independent from $\langle \psi \rangle$, and is described by the local well-posedness results in $H^s_0(\mathbb{T}) \times \dot H^s(\mathbb{T})$; on the other hand, the evolution of $\langle \psi \rangle$ is determined by taking the average of the second equation in \eqref{eq:WWsys} and by noticing that the right-hand side does not depend on $\langle \psi \rangle$. In the case $\gamma=0$ we refer also to \cite{alazard2016cauchy} and \cite{alazard2018control} for well-posedness and controllability results, respectively, in this and more generic functional settings.
\end{remark}

\begin{remark} \label{rem:AssLWP}

We point out that there exists $\epsilon_0>0$ such that for any $\epsilon \leq \epsilon_0$, if 
\begin{align*}
\| \eta_{0} \|_{ H^{s+1/2}_0(\mathbb{T}) } + \| \psi_{0} \|_{ H^{s}(\mathbb{T}) } &\leq \epsilon , 
\end{align*}
and if $\beta \in H^{s+1/2}(\mathbb{T})$ satisfies
\begin{align*}
| \beta(x) | &\leq h - h_0 - \epsilon_0 , \;\; \forall x \in \mathbb{T} ,
\end{align*}
then the assumption $\mathrm{(A1)}$ of Theorem \ref{thm:LWPMain} is satisfied. Regarding assumption $\mathrm{(A1)}$ and its generalizations in the irrotational case $\gamma=0$, we mention the local well-posedness result given by Theorem 1.1 in \cite{alazard2011water}  (see also Sec.4.4.2 of \cite{lannes2013water}). 

\end{remark}

\begin{remark}  \label{rem:AssTaylor}

As pointed out in Remark 9.7 of \cite{lannes2013water} for the irrotational case $\gamma=0$, the existence time $T$ in Theorem \ref{thm:LWPMain} is such that $T^{-1}$ depends on the Bond number $\mathrm{Bo} \sim g/\kappa$, which may be too small for some applications. In order to extend the time $T$ in such a way that $T^{-1}$ depends on $\mathrm{Bo}^{-1}$, one should impose the additional assumption 
\begin{align} \label{eq:Rayleigh}
\mathfrak{a}(\eta_0,\psi_0,\gamma,\beta) &\geq \mathfrak{a}_0 > 0 ,
\end{align}
where $\mathfrak{a}$ is a suitable function (for the case $\gamma=0$ see Theorem 9.6 in \cite{lannes2013water}). The assumption \eqref{eq:Rayleigh}, called Rayleigh-Taylor condition, is necessary for the local well-posedness in the pure gravity case $\kappa=0$ (see Theorem 4.16 of \cite{lannes2013water} in the irrotational case; see \cite{ifrim2018two} and \cite{castro2015well} for water waves with vorticity).



\end{remark}

\section{Properties of the operator $G(\eta,\beta,\gamma)$} \label{sec:OpG}

In this section we study the operator $G(\eta,\beta,\gamma)$ defined in \eqref{eq:OpG}, namely
\begin{align*}
G(\eta,\beta,\gamma)(\psi) &= \nabla\varphi \cdot \bmN |_{y=\eta(x)},
\end{align*}
where $\eta$ and $\beta$ satisfy \eqref{eq:StrConnected}, and where $\varphi$ is the solution of the boundary value problem \eqref{eq:EllBVP}, namely
\begin{equation*} 
\begin{cases}
\Delta \varphi = 0 \; \; &\text{in} \; \; D_{\eta,\beta}, \\
\varphi = \psi \; \; &\text{at} \; \; y=\eta, \\
\gamma \, y \, \beta_x - \varphi_x \, \beta_x + \varphi_y = 0  \; \; &\text{at} \; \; y=-h+\beta .
\end{cases}
\end{equation*}

We just mention that the operators $G(\eta,0,0)$ and $G(\eta,\beta,0)$ have been studied in \cite{craig1993numerical} and \cite{craig2005hamiltonian}, respectively (see also Chap. 3 of \cite{lannes2013water}). In the following we use the notation $D \coloneqq -\mathrm{i}\partial_{x}$.

We first introduce the notion of straightening diffeomorpshisms in Sec. \ref{subsec:straight}. Next, in Sec. \ref{subsec:analOpG} we prove Theorem \ref{thm:AnalOpGMain}, while in Sec. \ref{subsec:HomExp} we consider the homogeneous expansion of the operator $G(\eta,\beta,\gamma)$ around $\eta=0$; finally, in Sec. \ref{subsec:ParaOpG} we state the paralinearization formula for the operator $G(\eta,\beta,\gamma)$, which we prove in Sec. \ref{subsec:proofOpGpara}.

\subsection{Straightening Diffeomorphisms} \label{subsec:straight}

Observe that systems \eqref{eq:Incompr}-\eqref{eq:DynBC} and \eqref{eq:EllBVP} involve functions defined on the time-dependent fluid domain $D_{\eta,\beta}(t)$. In order to deal with this difficulty we introduce the following straightening diffeomorphism (also called admissible diffeomorphism in the literature, see Definition 2.13 in \cite{lannes2013water}) ,
\begin{align}
\Sigma(t) &: D_0 \to D_{\eta,\beta}(t) , \nonumber \\
\Sigma(t)(x,w) &:= ( x, w+ \sigma(t,x,w) ), \nonumber \\
\sigma(t,x,-h) = \beta(x)\in H^{s}(\mathbb{T}) , &\; \; \sigma(t,x,0)=\eta(t,x)\in H^{s}(\mathbb{T}), \; \; s > 3/2,  \label{eq:straight} 
\end{align}
where $\eta$ and $\beta$ satisfy condition \eqref{eq:StrConnected}, the coefficients of the Jacobian matrix $J_{\Sigma} = \nabla_{x,w}\Sigma$ belong to $L^\infty(D_0)$, and there exists a constant 
\begin{align*}
M_0 &= M_0( h_0^{-1}, \|\eta\|_{ H^{s}(\mathbb{T})} , \|\beta\|_{ H^{s}(\mathbb{T})},  \|\eta\|_{ C^{1}(\mathbb{T})} , \|\beta\|_{ C^{1}(\mathbb{T})} )>0
\end{align*}
such that
\begin{align*}
\| J_{\Sigma,i,j} \|_{L^\infty(D_0)} &\leq M_0 , \; \; \forall i,j=1,2,
\end{align*}
and such that the determinant of the Jacobian matrix $J_{\Sigma}$ is uniformly bounded from below on $D_0$ by $c_0>0$ satisfying $M_0 \geq 1/c_0$. 

We now make some remarks about straightening diffeomorphisms, and we discuss two common examples of them.

\begin{remark} \label{rem:Straight}

From the above definition of the straightening diffeomorphism, and in particular from the fact that the determinant of the Jacobian matrix $J_{\Sigma}$ is bounded from below, we deduce that the function $\sigma$ satisfies
\begin{align*}
1 + \partial_w \sigma(t) &\neq 0, \; \; \text{in} \; \; D_0.
\end{align*}
\end{remark}

\begin{remark} \label{rem:TrivialDiffeo}

A so-called \emph{trivial diffeomorphism} between $D_0$ and $D_{\eta,\beta}(t)$ is given by choosing
\begin{align} \label{eq:TrivialDiffeo}
\sigma(t,x,w) &:= \left( 1 + \frac{w}{h} \right) \eta(t,x) - \frac{w}{h} \beta(x) 
\end{align}
in the formula \eqref{eq:straight} (see also Remark 2.4 in  \cite{lannes2005well}).

\end{remark}

We now introduce another example of straightening diffeomorphism, the so-called \emph{regularizing diffeomorphism} (see Proposition 2.16 and Proposition 2.18 of \cite{lannes2013water}; see also Proposition 2.13 of \cite{lannes2005well}). We recall that $D=-\mathrm{i}\partial_x$, and we write $\langle D \rangle \coloneqq (1-\partial_{xx})^{1/2}$; we also define, for any $s \in \mathbb{R}$ and $k\in \mathbb{N}$
\begin{align*}
H^{s,k}(D_0) &\coloneqq \bigcap_{j=0}^k H^{j}( (-h,0) ;H^{s-j}(\mathbb{T}) ), \;\; \| f \|_{ H^{s,k}(D_0) } \coloneqq \sum_{j=0}^k \| \langle D \rangle \, \partial_w^j f \|_{ L^2(\mathbb{T}) } .
\end{align*}

\begin{proposition} \label{prop:RegDiffeo}

Let $s>3/2$ and assume that $\eta \in H^{s}(\mathbb{T})$ and $\beta \in H^{s}(\mathbb{T})$ satisfy condition \eqref{eq:StrConnected}. Let $\chi \in C^{\infty}_c(\mathbb{R})$ be a positive function equal to $1$ in a neighborhood of the origin, and consider the diffeomorphism \eqref{eq:straight}, where
\begin{align}
\sigma(t,x,w) &\coloneqq \left( 1 + \frac{w}{h} \right) \eta^{(\delta)}(t,x,w) - \frac{w}{h} \beta_{(\delta)}(x) ,  \label{eq:straightReg} \\
\eta^{(\delta)}(t,x,w) &\coloneqq \chi(\delta \, w \, |D|)\eta(t,x)  , \;\; \beta_{(\delta)}(x,w) \coloneqq \chi(\delta \, (w+h) \, |D|)\beta(x) . \nonumber
\end{align}
Then, if $\delta >0$ in \eqref{eq:straightReg} satisfies
\begin{align*}
\delta &< \frac{ h_0 }{C(\chi) \left[ \|\eta\|_{ H^{s}(\mathbb{T})} + \|\beta\|_{ H^{s}(\mathbb{T})} \right] } , \nonumber \\
C(\chi) &:= \|\chi'\|_{L^\infty(\mathbb{R})} \, \left[ \int_{\mathbb{R}^2} (1+|\xi|^2)^{-(s-1)} \, \mathrm{d}\xi \right]^{1/2} < +\infty ,
\end{align*}
then the diffeomorphism defined in \eqref{eq:straightReg} is such that $J_{\Sigma} = \nabla_{x,w}\Sigma$ belong to $( L^\infty(D_0) )^0$, and there exists $M_0>0$ such that
\begin{align*}
\| J_{\Sigma,i,j} \|_{L^\infty(D_0)} &\leq M_0 , \; \; \forall i,j=1,2.
\end{align*}
Moreover, the determinant of the Jacobian matrix $J_{\Sigma}$ is uniformly bounded from below on $D_0$ such that $M_0 \geq 1/c_0$, where $c_0 = h_0 - \delta \, C(\chi) \, \left[ \|\eta\|_{ H^{s}(\mathbb{T}) } + \|\beta\|_{ H^{s}(\mathbb{T})} \right]$.

\end{proposition}

We also mention the following smoothing property (see Lemma 2.20 in \cite{lannes2013water}).

\begin{lemma} \label{lem:smoothing}

Let $s \in \mathbb{R}$ and $\lambda_1>0$, then for all $\chi \in C_c(\mathbb{R})$ there exists a constant $C=C(\chi)>0$ such that
\begin{align} \label{eq:smoothing}
\left\| \chi(\lambda_1^{1/2} w |D|) f \right\|_{H^s(D_0)}^2 &\leq C \, \left\| (1+\lambda_1^{1/2} |D|)^{-1/2} f \right\|_{H^s(\mathbb{T})}^2 .
\end{align}
\end{lemma}

Observe that by Lemma \ref{lem:smoothing} (with $\lambda_1=\delta^2$) we have that 
\begin{align*}
\| \eta^{(\delta)} \|_{H^{s+1/2,2}(D_0)} &\leq C(1/\delta) \|\eta\|_{H^{s}(\mathbb{T})} . 
\end{align*}
As a comparison, recall that for the trivial diffeomorphism of Remark \ref{rem:TrivialDiffeo} one need to control $H^s(\mathbb{T})$-norm of the surface parameterization $\eta$ in order to control the $H^s(D_0)$-norm of $\sigma$ in \eqref{eq:TrivialDiffeo}. \\

For simplicity, in the following discussion we assume that $\eta$, $\beta$, the function $\sigma$ in \eqref{eq:straight}, the function $f: D_{\eta,\beta} \to \mathbb{R}$ and the vector field $\bmF: D_{\eta,\beta} \to \mathbb{R}^2$ are sufficiently smooth, so that the quantities we introduce below are well-defined.

We define $\tilde{\bmF}(x,w) := \bmF \circ \Sigma(x,w)$ and $\tilde{f}(x,w)=f \circ \Sigma(x,w)$; we also write
\begin{align*}
\partial_i^\Sigma \tilde{\bmF} &:= ( \partial_i \bmF ) \circ \Sigma, \;\; i=t,x,w ,
\end{align*}
and similarly for $\partial_i^\Sigma \tilde{f}$. From \eqref{eq:straight} we obtain that
\begin{align*}
\partial_j^\Sigma = \partial_j - \frac{ \partial_j \sigma }{ 1+\partial_w \sigma} \partial_w , \; \; j =x,t, &\; \; \partial_w^\Sigma = \frac{ \partial_w }{ 1+\partial_w \sigma } .
\end{align*}

\begin{remark} \label{rem:flatOp}

Similarly, we can define
\begin{align*}
\Div^{\Sigma} \tilde{\bmF} = \nabla^\Sigma \cdot \tilde{\bmF} := (\Div \bmF) \circ \Sigma , &\;\; \Delta^\Sigma  \tilde{f} := (\Delta f) \circ \Sigma .
\end{align*}

One can check that the flattened version of the operators $\Div$ and $\Delta$ applied to $\tilde{\bmF}$ and to $\tilde{f}$ are given by
\begin{align*}
\Div^\Sigma \tilde{\bmF} = \partial_x^\Sigma \tilde{F}_1 + \partial_w^\Sigma \tilde{F}_2 &= \Div \tilde{\bmF} - \frac{1}{1+\partial_w \sigma} \, \partial_x \sigma \; \partial_w \tilde{F}_1  - \frac{ \partial_w \sigma }{ 1+\partial_w \sigma}  \; \partial_w \tilde{F}_2 , \\
\Delta^\Sigma \tilde{f} = (\partial_x^\Sigma)^2 \tilde{f}  + (\partial_w^\Sigma)^2 \tilde{f} &= \left[ \mathtt{a} \, \partial_w^2 + \partial_x^2 + \mathtt{b} \, \partial_x \partial_w - \mathtt{c} \, \partial_w \right] \tilde{f},
\end{align*}
\begin{align} 
\mathtt{a} := \frac{ 1+ (\partial_{x} \sigma)^2 }{ (1+\partial_w \sigma)^2 } , \; \; \mathtt{b} &:= -2 \frac{ \partial_{x} \sigma }{ 1+\partial_w \sigma } , \;\; \mathtt{c} := \frac{1}{ 1+\partial_w \sigma } \left[  \partial_x^2 \sigma  + \mathtt{b} \,  \partial_{x} \, \partial_w \sigma  +  \mathtt{a} \,  \partial_w^2 \sigma \right] . \label{eq:CoeffFlatLap}
\end{align}

Moreover, if we denote by $J_\Sigma$ the Jacobian matrix of $\Sigma$, we denote by $P(\Sigma)$ the symmetric matrix given by
\begin{align}
P(\Sigma) &:= |\det(J_\Sigma)| \, J_{\Sigma}^{-1}  J_\Sigma^{-T}  = (1+\partial_w \sigma) \,
\begin{pmatrix}
1 & \frac{- \partial_{x} \sigma }{ 1+ \partial_w \sigma} \\
& \\
\frac{- ( \partial_{x} \sigma )^T }{ 1+\partial_w \sigma} & \frac{ 1+ (\partial_{x} \sigma)^2 }{ (1+\partial_w \sigma)^2 } 
\end{pmatrix}
, \label{eq:PSigma} 
\end{align}
so that 
\begin{align*}
\Delta^{\Sigma} \tilde{f} &= (1 + \partial_w\sigma)^{-1} \, \nabla_{x,w} \cdot P(\Sigma)\nabla_{x,w} \tilde{f} ,
\end{align*}
see Lemma 2.5 in \cite{lannes2005well}.

\end{remark}

\subsection{Analyticity} \label{subsec:analOpG}

First we prove an analyticity result for the operator $G(\eta,\beta,\gamma)$ where $(\eta,\beta)$ are in a neighbourhood of the origin in suitable Sobolev spaces, by following the approach of Sec. 4 in \cite{groves2020variational}; this allows us to derive a homogeneous expansion around $\eta=0$ in Sec. \ref{subsec:HomExp}. 

In the following we use the straightening diffeomorphisms \eqref{eq:straight} introduced in Sec. \ref{subsec:straight}. We consider the flattened version of \eqref{eq:EllBVP}. Let $s > 3/2$, and fix $\beta \in H^s(\mathbb{T})$, $\gamma \in \mathbb{R}$. We obtain that $\tilde{\varphi} =\varphi \circ \Sigma$ satisfies

\begin{align}
(D^2 - \partial_w^2) \tilde{\varphi} &= r_{1,0}[\eta,\beta,\tilde{\varphi}] , \; \; \text{in} \; \; D_0, \label{eq:systphieq1} \\
\tilde{\varphi}(\cdot,0) &=  \psi , \label{eq:systphieq2} 
\end{align}
\begin{align}
\gamma (-h+\beta) \beta_x - \beta_x \, \partial_x \tilde{\varphi}(\cdot,-h) + \partial_w \tilde{\varphi}(\cdot,-h) &= r_{2,0}[\eta,\beta,\tilde{\varphi}], \label{eq:systphieq3} 
\end{align}
where
\begin{align*}
r_{1,0}[\eta,\beta,\tilde{\varphi}] &:=   \Delta^{\Sigma} \tilde{\varphi} - \Delta_{x,w}\tilde{\varphi} , \\
r_{2,0}[\eta,\beta,\tilde{\varphi}] &:= - ( \partial^{\Sigma}_w \tilde{\varphi} ) (\cdot,-h) + \partial_w \, \tilde{\varphi}(\cdot,-h) \\
&\quad  -\beta_x \, \left[ - ( \partial^{\Sigma}_x \tilde{\varphi} ) (\cdot,-h) + \partial_x \, \tilde{\varphi}(\cdot,-h) \right] .
\end{align*}

In the following we denote by $\mathscr{F}$ the Fourier transform with respect to $x$, and we denote by $\tilde{f}'$ the partial derivative $\partial_w \, \tilde{f}$.

\begin{lemma} \label{lem:SolFlatBVP}

Let $s > 3/2$, and fix $\beta \in H^s(\mathbb{T})$, $\gamma \in \mathbb{R}$. Then the boundary value problem
\begin{equation} \label{eq:BVPtildevarphi} 
\begin{cases}
(D^2 - \partial_w^2) \tilde{\varphi} &= r_{1} , \; \; \text{in} \; \; D_0, \\
\tilde{\varphi}(\cdot,0) &=  \psi , \\
\gamma (-h+\beta) \beta_x - \beta_x \, \partial_x \tilde{\varphi}(\cdot,-h)  + \partial_w \, \tilde{\varphi}(\cdot,-h) &= r_{2} ,
\end{cases}
\end{equation}
admits a unique solution $\tilde{\varphi} \in  H^{s+1/2}( D_0 )$ for any $\psi \in H^{s}(\mathbb{T})$, $r_1 \in H^{s-3/2}(D_0)$ and $r_2 \in H^{s-1}(\mathbb{T})$. 
\end{lemma}

\begin{proof}

We rewrite the first boundary value problem for $\tilde{g} := \mathscr{F}\tilde{\varphi}$ in the following way,
\begin{align*}
\mathtt{M} \tilde{g} &= \mathscr{F} r_1 \; \; \text{in} \; \; D_0, \\
\mathtt{L}_{0}\tilde{g} &=  \mathscr{F} \psi ,\; \; \text{at} \; \; w = 0, \\
\mathtt{L}_{-h}\tilde{g} &= \mathscr{F} \left[ -\gamma(-h+\beta) \, \beta_x + r_2 \right], \; \; \text{at} \; \; w = -h,  
\end{align*}
where
\begin{align*}
\mathtt{M} = (\xi^2-\partial_w^2) , \;\; \mathtt{L}_{0} = 1 , \; \; \mathtt{L}_{-h} &= ( D\beta ) \, \xi + \partial_w .
\end{align*}
Let us define the Green's matrix for the boundary value problem for $\tilde{g}$,
\begin{align*}
\mathtt{G}_{\tilde{\varphi}}(w,\zeta) &=
\begin{cases}
\mathtt{u}(w+h) \, \overline{\mathtt{c}^{-1}} \, \overline{\mathtt{w}(\zeta)}, \; \; \text{for} \; \; -h \leq w \leq \zeta \leq 0, \\
\mathtt{w}(w) \, \mathtt{c}^{-1} \, \overline{\mathtt{u}(\zeta+h)}, \; \; \text{for} \; \; -h \leq \zeta \leq w \leq 0, 
\end{cases}
\end{align*}
where $\mathtt{u}(\cdot+h)$ (resp. $\mathtt{w}$) is a solution of 
\begin{align*}
(\xi^2-\partial_w^2)\tilde{g}=0, \; \; \text{in}\; \; D_0, &\; \; ( D\beta ) \, \xi \, \tilde{g}(-h) + \tilde{g}'(-h)=0, \\
(\text{resp}. \; \; (\xi^2-\partial_w^2)\tilde{g}=0, \; \; \text{in}\; \; D_0, &\; \; \tilde{g}(0)=0  )
\end{align*}
and
\begin{align*}
\mathtt{c} &= -\overline{\mathtt{u}}(w+h) \, \mathtt{w}'(w) + \overline{\mathtt{u}'}(w+h) \, \mathtt{w}(w) .
\end{align*}
More explicitly, we choose
\begin{align*}
\mathtt{u}(w) = \cosh( \xi w ) - (D\beta) \, \sinh( \xi w ) , &\; \; \mathtt{w}(w) = \sinh( \xi w ) , 
\end{align*}
so that
\begin{align*}
\mathtt{c} &= -\xi \, \cosh(h \, \xi) + \xi \, (D\beta) \, \sinh(h \, \xi) ,
\end{align*}
and
\begin{align*}
\mathtt{G}_{\tilde{\varphi}}(w,\zeta) &=
\begin{cases}
\frac{ [ \cosh(\xi(w+h) ) - (D\beta) \, \sinh( \xi(w+h) ) ] \, \sinh(\xi \zeta) }{ - \xi \, \cosh(h \, \xi) + \xi \, (D\beta) \, \sinh(h \, \xi) } , \;\;  \text{for} \; \; -h \leq w \leq \zeta \leq 0, \\
\\
\frac{ [ \cosh(\xi(\zeta+h) ) - (D\beta) \, \sinh( \xi(\zeta+h) ) ] \, \sinh(\xi w) }{ - \xi \, \cosh(h \, \xi) + \xi \, (D\beta) \, \sinh(h \, \xi) } , \;\;  \text{for} \; \; -h \leq \zeta \leq w \leq 0.
\end{cases}
\end{align*}

Hence the solution $\tilde{\varphi}$ of the boundary value problem \eqref{eq:BVPtildevarphi}  is given by
\begin{align}
\tilde{\varphi}(\cdot,w) &= \mathscr{F}^{-1} \, \bigg[ 
\int_{-h}^0 \mathtt{G}_{\tilde{\varphi}}(w,\zeta) \, \mathscr{F}r_1(\cdot,\zeta) \, \mathrm{d}\zeta \nonumber \\
&\qquad \qquad - \mathtt{G}_{\tilde{\varphi}}(w,-h) \, \mathscr{F} [ -\gamma (-h+\beta) \, \beta_x + r_2 ]  - \mathtt{G}_{\tilde{\varphi},\zeta}(w,0) \, (\mathscr{F} \psi) \, \bigg] . \nonumber \\
&\phantom{} \label{eq:BVPsolphi}
\end{align}
Notice that by using the explicit formula for $\mathtt{G}_{\tilde{\varphi}}$ we have
\begin{equation*}
\mathtt{G}_{\tilde{\varphi}}(w,0) = 0, \;\; \mathtt{G}_{\tilde{\varphi},\zeta}(w,0) = - \frac{ \cosh(\xi (w+h)) - (D\beta) \, \sinh( \xi (w+h) ) }{ \cosh(h \, \xi) - (D\beta) \, \sinh(h \, \xi) } \, (\mathscr{F} \psi) ,
\end{equation*}
so that \eqref{eq:BVPsolphi} takes the form
\begin{align}
\tilde{\varphi}(\cdot,w) &= \mathscr{F}^{-1} \, \bigg[ 
\int_{-h}^0 \mathtt{G}_{\tilde{\varphi}}(w,\zeta) \, \mathscr{F}r_1(\cdot,\zeta) \, \mathrm{d}\zeta \nonumber \\
&\qquad \qquad +  \frac{ \sinh(\xi \, w) }{ - \xi \, \cosh(h \, \xi) + \xi \, (D\beta) \, \sinh(h \, \xi) } \, \mathscr{F} \left[ \gamma (-h+\beta) \beta_x - r_2 \right] \nonumber \\
&\qquad \qquad +  \frac{ \cosh(\xi (w+h)) - (D\beta) \, \sinh( \xi (w+h) ) }{ \cosh(h \, \xi) - (D\beta) \, \sinh(h \, \xi) } \, (\mathscr{F} \psi)  \, \bigg] . \nonumber \\
&\phantom{} \label{eq:BVPsolphi2}
\end{align}

Moreover, notice that
\begin{align*}
& \frac{ \cosh(\xi (w+h)) - (D\beta) \, \sinh( \xi (w+h) ) }{ \cosh(h \, \xi) - (D\beta) \, \sinh(h \, \xi) }  \\
&= \frac{ \cosh(\xi (w+h)) }{ \cosh(h \, \xi) } - \sinh( \xi \, w ) \frac{ D\beta }{ \cosh(h \, \xi) \, \left[ \cosh(h \, \xi) - (D\beta) \, \sinh(h \, \xi) \right] }  .
\end{align*}

\end{proof}

Following an argument analogous to the proof of Theorem 4.10 in \cite{groves2020variational}, one can deduce from Lemma \ref{lem:SolFlatBVP} the following result. 

\begin{proposition} \label{prop:ExSolFBVP}

\begin{itemize}
\item[i.] Let $s > 3/2$, $\gamma \in \mathbb{R}$ and $\psi \in H^s(\mathbb{T})$. There exist open neighbourhoods $W_1,W_2$ of the origin in $H^{s+1}(\mathbb{T})$ such that for any $\eta \in W_1$ and $\beta \in W_2$ satisfying \eqref{eq:StrConnected}, and for any function $\sigma \in H^{s+1}(D_0)$ in the diffeomorphism \eqref{eq:straight} for which there exists $C>0$ such that
\begin{align} \label{DiffeoAss}
\| \sigma \|_{ H^{s+1}(D_0) } &\leq C \left[ \| \eta \|_{ H^{s+1}(\mathbb{T}) } + \| \beta \|_{ H^{s+1}(\mathbb{T}) } \right] ,
\end{align}
then the boundary value problem \eqref{eq:systphieq1}-\eqref{eq:systphieq3} admits a unique solution 
\begin{align*}
\tilde{\varphi} &= \tilde{\varphi}(\eta,\beta,\gamma,\psi) \in H^{s+1/2}(D_0) .
\end{align*}
Furthermore, $\tilde{\varphi}(\eta,\beta,\gamma,\psi)$ depends analytically upon $\eta$ and $\beta$, and $\psi \mapsto \tilde{\varphi}(\eta,\beta,\gamma,\psi)$ is an affine map from $H^s(\mathbb{T})$ to $H^{s+1/2}(D_0)$.

\item[ii.] If we choose $\sigma$ as in \eqref{eq:straightReg}, then in the statement of part $\mathrm{i.}$ we can take $W_1,W_2$ as open neighbourhoods of the origin in $H^{s+1/2}(\mathbb{T})$.
\end{itemize}

\end{proposition}

\begin{proof}

First, we prove part $\mathrm{i.}$: let $s > 3/2$ and $\gamma \in \mathbb{R}$. By Sobolev embedding (see Proposition 3.3 in \cite{taylor2011partial}; this holds true also for Sobolev spaces on $D_0$, see also pag.16 of \cite{groves2020variational}), one can see that the maps

\begin{equation*}
(\eta,\beta,\tilde{\varphi}) \mapsto r_{1,0}[\eta,\beta,\tilde{\varphi}]  , \;\; (\eta,\beta,\tilde{\varphi}) \mapsto r_{2,0}[\eta,\beta,\tilde{\varphi}] 
\end{equation*}
are analytic in $W_1 \times W_2 \times  H^{s+1/2}(D_0)$, where $W_1,W_2$ are open neighbourhoods of $H^{s+1}(\mathbb{T})$. Hence the map
\begin{align*}
\mathcal{F}(\eta,\beta,\tilde{\varphi},\psi) &:= 
\begin{pmatrix}
(D^2 - \partial_w^2) \tilde{\varphi} - r_{1,0}[\eta,\beta,\tilde{\varphi}] \\
\tilde{\varphi}(\cdot , 0) -  \psi \\
\gamma (-h+\beta) \, \beta_x - \beta_x \, \partial_x \tilde{\varphi}(\cdot , -h) +  \partial_w \tilde{\varphi}(\cdot , -h) - r_{2,0}[\eta,\beta,\tilde{\varphi}] 
\end{pmatrix}
\end{align*}
define an analytic map  
\begin{equation*}
\mathcal{F} :  W_1 \times W_2 \times H^{s+1/2}(D_0) \times H^{s}(\mathbb{T})  \to  H^{s-3/2}(D_0)  \times H^{s}(\mathbb{T}) \times H^{s-1}(\mathbb{T}) , 
\end{equation*}
where $\mathcal{F}(0,0,0)  = (0,0,0)$. Moreover, Lemma \ref{lem:SolFlatBVP} implies that the map 
 $\mathrm{d}_1 \mathcal{F}[0,0,0,0] :  H^{s+1}(\mathbb{T}) \times H^{s+1}(\mathbb{T}) \times H^{s+1/2}(D_0) \times  H^{s}(\mathbb{T})  \to  H^{s-3/2}(D_0)  \times H^{s}(\mathbb{T}) \times  H^{s-1}(\mathbb{T}) $ is an isomorphism. The thesis follows by the analytic implicit function theorem (see Sec. 4.5 of \cite{buffoni2003analytic}).

Part $\mathrm{ii.}$ is proved similarly, using the smoothing property \eqref{eq:smoothing} of Lemma \ref{lem:smoothing}.

\end{proof}

We point out that the assumption \eqref{DiffeoAss} in Proposition \ref{prop:ExSolFBVP} and Corollary \ref{cor:AnalOpG} holds true if we choose the trivial diffeomorphism \eqref{eq:straight}-\eqref{eq:TrivialDiffeo}.

From Proposition \ref{prop:ExSolFBVP} we can deduce the following result (see also Corollary 4.12 of \cite{groves2020variational}).

\begin{corollary} \label{cor:AnalOpGsmall}

\begin{itemize}
\item[i.] Let $s > 3/2$, $\gamma \in \mathbb{R}$ and $\psi \in H^s(\mathbb{T})$. There exist open neighbourhoods $W_1,W_2$ of the origin in $H^{s+1}(\mathbb{T})$ such that for any $\eta \in W_1$ and $\beta \in W_2$ satisfying \eqref{eq:StrConnected}, and for any function $\sigma \in H^{s+1}(D_0)$ in the diffeomorphism \eqref{eq:straight} for which the estimate \eqref{DiffeoAss} holds true, the formula
\begin{align*}
G(\eta,\beta,\gamma)(\psi) &= \nabla^{\Sigma}\tilde{\varphi}(\eta,\beta,\gamma,\psi) \cdot \bmN |_{w=0} 
\end{align*}
defines an analytic map
\begin{align*}
G(\cdot,\cdot,\gamma)(\psi) &: \{ (\eta,\beta) \in W_1 \times W_2 : \; \eqref{eq:StrConnected} \; \text{holds true} \} \to H^{s-1}(\mathbb{T}) .
\end{align*}

\item[ii.] If we choose $\sigma$ as in \eqref{eq:straightReg}, then in the statement of part $\mathrm{i.}$ we can take $W_1,W_2$ as open neighbourhoods of the origin in $H^{s+1/2}(\mathbb{T})$.
\end{itemize}

\end{corollary}

Next, following the argument of Sec. A.1 in \cite{lannes2013water} we can prove an analyticity result for the operator $G(\eta,\beta,\gamma)$ which is not restricted to small perturbations of a flat free surface and of a flat bottom topography. In the following results we use the regularizing diffeomorphism \eqref{eq:straightReg} for simplicity, but more general choices of regularizing diffeomorphisms are possible (see Definition 2.17 in \cite{lannes2013water}). In the following we write $\langle \partial_x \rangle \coloneqq (1-\partial_{xx})^{1/2}$, and we denote by $D_{0,bott}$ the set $D_0 \cup \{ w = -h \}$.

\begin{theorem} \label{thm:AnSolBVP}

Let $s > 3/2$, $\gamma \in \mathbb{R}$ and $\psi \in H^s(\mathbb{T})$. Then, for any $\eta,\beta \in H^{s}(\mathbb{T})$ satisfying \eqref{eq:StrConnected} and for any $\sigma$ of the form  \eqref{eq:straightReg} in the diffeomorphism \eqref{eq:straight}, the boundary value problem \eqref{eq:systphieq1}-\eqref{eq:systphieq3} admits a unique solution 
\begin{align*}
\tilde{\varphi} &= \tilde{\varphi}(\eta,\beta,\gamma,\psi) \in H^{s+1/2}(D_0) .
\end{align*}
Furthermore, $\tilde{\varphi}(\eta,\beta,\psi,\gamma)$ depends analytically upon $\eta$ and $\beta$, and $\psi \mapsto \tilde{\varphi}(\eta,\beta,\gamma,\psi)$ is an affine map from $H^s(\mathbb{T})$ to $H^{s+1/2}(D_0)$.
\end{theorem}

\begin{proof}

Let $\bmzeta_0 = (\eta_0,\beta_0) \in ( H^{s}(\mathbb{T}) )^2$, and let us consider the map 
\begin{equation*}
Q:\bmzeta=(\eta,\beta) \mapsto Q(\Sigma) \coloneqq P(\Sigma) - \mathbb{I}_2 ,
\end{equation*}
where $P(\Sigma)$ is the matrix defined in \eqref{eq:PSigma}. Then in a small neighbourhood of $\bmzeta_0$ in $( H^{s}(\mathbb{T}) )^2$ we can write
\begin{align*}
Q(\bmzeta) &= \sum_{k =0}^{+\infty} Q_k(\bmzeta-\bmzeta_0)^k ,
\end{align*}
where each $Q_k: ( H^{s}(\mathbb{T}) )^{2k} \to M_2( H^{s-1/2,2}(D_0) )$ is bounded, symmetric, $k$-linear and satisfies
\begin{align*}
\| Q_k \|_{ ( H^{s}(\mathbb{T}) )^{2k} \to M_2( H^{s-1/2,2}(D_0) ) } &\leq N_0 r_0^{-k} , 
\end{align*}
for some $N_0,r_0>0$. 

We want to prove that in a possibly smaller neighbourhood of $\bmzeta_0$ we have
\begin{align*}
\tilde{\varphi}(\bmzeta,\gamma,\psi) &= \sum_{k =0}^{+\infty} \tilde{\varphi}_k(\bmzeta-\bmzeta_0,\gamma,\psi)^k ,
\end{align*}
where each $\tilde{\varphi}_k(\cdot,\gamma,\psi): ( H^{s}(\mathbb{T}) )^{2k} \to H^{s+1/2}(D_0) $ is bounded, symmetric, $k$-linear and satisfies
\begin{align} \label{eq:varphikEst}
\| \tilde{\varphi}_k(\cdot,\gamma,\psi) \|_{ ( H^{s}(\mathbb{T}) )^{2k} \to H^{s+1/2}(D_0) } &\leq N_1 r_1^{-k} , 
\end{align}
for some $N_1,r_1>0$. 

By identifying terms with the same degree of homogeneity, we have that, if they exist, then $a_k \coloneqq \tilde{\varphi}(\bmzeta-\bmzeta_0,\gamma,\psi)^k$ solve the following systems
\begin{equation} \label{eq:IndSys0}
\begin{cases}
\nabla \cdot (\mathbb{I}_2 + Q_0)\nabla\tilde{a}_0 =0, &\text{in} \;\; D_0, \\
\tilde{a}_0 =\psi, &\text{at} \;\; w=0, \\
-\bme_2 \cdot (\mathbb{I}_2 + Q_0)\nabla\tilde{a}_0 =0, &\text{at} \;\; w=-h,
\end{cases}
\end{equation}

\begin{equation} \label{eq:IndSys1}
\begin{cases}
\nabla \cdot (\mathbb{I}_2 + Q_0)\nabla\tilde{a}_1 = - \nabla \cdot Q_1 \nabla \tilde{a}_0, &\text{in} \;\; D_0, \\
\tilde{a}_1 =0, &\text{at} \;\; w=0, \\
-\bme_2 \cdot (\mathbb{I}_2 + Q_0)\nabla\tilde{a}_1 = -\bme_2 \cdot Q_1\nabla\tilde{a}_0 -\gamma h (\beta-\beta_0)_x, &\text{at} \;\; w=-h,
\end{cases}
\end{equation}

\begin{equation} \label{eq:IndSys2}
\begin{cases}
\nabla \cdot (\mathbb{I}_2 + Q_0)\nabla\tilde{a}_2 = - \nabla \cdot \sum_{j=1,2} Q_j \nabla \tilde{a}_{k-j}, &\text{in} \;\; D_0, \\
\tilde{a}_2 =0, &\text{at} \;\; w=0, \\
-\bme_2 \cdot (\mathbb{I}_2 + Q_0)\nabla\tilde{a}_2 = -\bme_2 \cdot \sum_{j=1,2} Q_j\nabla\tilde{a}_{k-j} +\gamma (\beta-\beta_0)(\beta-\beta_0)_x, &\text{at} \;\; w=-h,
\end{cases}
\end{equation}

\begin{equation} \label{eq:IndSysk}
\begin{cases}
\nabla \cdot (\mathbb{I}_2 + Q_0)\nabla\tilde{a}_k = - \nabla \cdot \sum_{j=1}^k Q_j \nabla \tilde{a}_{j-k}, &\text{in} \;\; D_0, \\
\tilde{a}_k =0, &\text{at} \;\; w=0, \\
-\bme_2 \cdot (\mathbb{I}_2 + Q_0)\nabla\tilde{a}_k = -\bme_2 \cdot \sum_{j=1}^k Q_j\nabla\tilde{a}_{j-k}, &\text{at} \;\; w=-h,
\end{cases}
\;\; k \geq 3,
\end{equation}
where we denoted by $Q_k$ the term $Q_k(\bmzeta-\bmzeta_0)^k$.

\begin{itemize}
\item Claim: There exist a unique solution $(a_k)_{k \in \mathbb{N}}$ of \eqref{eq:IndSys0}-\eqref{eq:IndSysk} such that $a_k \in H^{s+1/2}(D_0)$. Moreover, there exist 
\begin{align*}
N_1 &= N_1( \|\eta_0\|_{ H^{s}(\mathbb{T}) } , \|\beta_0\|_{ H^{s}(\mathbb{T}) }, \|\psi\|_{H^{s}(\mathbb{T})} ) > 0
\end{align*}
and $r_1 >0$ such that
\begin{align*}
\|  \langle \partial_x \rangle^m \, \nabla \tilde{a}_k \|_{L^2(D_0)} &\leq  N_1 r_1^{-k} \|\bmzeta-\bmzeta_0\|_{ ( H^{s}(\mathbb{T}) )^2 }^k ,  \forall 0 \leq m \leq s-1/2, \\
\|  \nabla \tilde{a}_k \|_{ H^{s-1/2,1}(D_0)} &\leq  N_1 r_1^{-k} \|\bmzeta-\bmzeta_0\|_{ ( H^{s}(\mathbb{T}) )^2 }^k , 
\end{align*}
for all $k \in \mathbb{N}$. \\

We prove the claim: by Proposition A.20 of \cite{lannes2013water} there exists a unique variational solution to \eqref{eq:IndSys0} such that
\begin{align*}
\|  \langle \partial_x \rangle^m \, \nabla \tilde{a}_0 \|_{L^2(D_0)} &\leq  N_1  ,  \forall 0 \leq m \leq s-1/2 .
\end{align*}

We argue by induction: assume that there exist $(\tilde{a}_j)_{0 \leq j \leq k-1}$ such that 
\begin{align*}
\|  \langle \partial_x \rangle^m \, \nabla \tilde{a}_j \|_{L^2(D_0)} &\leq  N_1 r_1^{-j} \|\bmzeta-\bmzeta_0\|_{ ( H^{s}(\mathbb{T}) )^2 }^j ,  \forall 0 \leq m \leq s-1/2, \\
\|  \nabla \tilde{a}_j \|_{ H^{s-1/2,1}(D_0)} &\leq  N_1 r_1^{-j} \|\bmzeta-\bmzeta_0\|_{ ( H^{s}(\mathbb{T}) )^2 }^j , 
\end{align*}
for all $0 \leq j \leq k-1$, we construct $\tilde{a}_k$ satisfying the thesis of the claim. We just prove the first estimate for $m \leq s-1$, since the others follow from the embedding $H^{s-1/2,1}(D_0) \subset L^{\infty}([-h,0]; H^{s-1}(\mathbb{T}) )$ and by interpolation.

Let $\bmg_k = \sum_{j=1}^k Q_j \nabla\tilde{a}_{j-k}$, then by Corollary B.5 of \cite{lannes2013water} we have
\begin{align*}
&\|  \langle \partial_x \rangle^m \, \bmg_k \|_{L^2(D_0)} \\
&\leq \sum_{j=1}^k \|Q_j\|_{M_2( L^{\infty}([-h,0] ; H^{s-1/2}(\mathbb{T}) ) )} \, \| \langle \partial_x \rangle^m \, \nabla \tilde{a}_{k-j} \|_{L^2(D_0)} , \;\; \forall m \leq s-1/2,
\end{align*}
so that by the assumption on the norm of $Q_j$ and by the induction assumption we obtain that
\begin{align*}
&\| \langle \partial_x \rangle^m \, \bmg_k \|_{L^2(D_0)} \\
&\leq \left[\sum_{j=1}^k N_0 r_0^{-j} N_1 r_1^{-(k-j)} \right] \| \bmzeta-\bmzeta_0 \|^k_{ ( H^s(\mathbb{T}) )^2} ,  \;\; \forall m \leq s-1/2 .
\end{align*}

Therefore, by Proposition A.20 of \cite{lannes2013water} there exists a unique $\tilde{a}_k$ solving \eqref{eq:IndSysk}; moreover, we have
\begin{align*}
&\|  \langle \partial_x \rangle^m \, \nabla\tilde{a}_k \|_{L^2(D_0)} \\
&\leq C\left( \|\bmzeta_0\|_{ ( H^s(\mathbb{T}) )^2} , |\gamma|  \right) \, \left[\sum_{j=1}^k N_0 r_0^{-j} N_1 r_1^{-(k-j)} \right] \| \bmzeta-\bmzeta_0 \|^k_{( H^s(\mathbb{T}) )^2} ,  \;\; \forall m \leq s-1/2 ,
\end{align*}
so that by choosing $r_1$ such that 
\begin{align*}
C\left( \|\bmzeta_0\|_{ ( H^s(\mathbb{T}) )^2} , |\gamma|  \right) \, N_0 \sum_{j=1}^{\infty} (r_1/r_0)^{j} &\leq 1 ,
\end{align*}
we can deduce the thesis for $j=k$, which proves the claim.
\end{itemize}

Notice that the dependence of $\tilde{a}_k$ on $\bmzeta-\bmzeta_0$ is homogeneous of order $k$, hence by the above claim we can find a symmetric $k$-linear mapping $\tilde{\varphi}_k$ such that $\tilde{\varphi}_k(\bmzeta-\bmzeta_0)^k = \tilde{a}_k$ and satisfying \eqref{eq:varphikEst}; in particular, if we denote by $H^1_{0,surf}(D_0)$ the completion of $C^{\infty}_c(D_{0,bott} )$ in $H^1(D_{0})$, we have that $\tilde{b}=\sum_{k \geq 1} \tilde{a}_k$ converges absolutely in $H^1_{0,surf}(D_{0})$ for $\bmzeta$ sufficiently close to $\bmzeta_0$. Next, by \eqref{eq:IndSys0}-\eqref{eq:IndSysk} we have
\begin{align*}
\int_{D_0} \nabla_{x,w}(\tilde{a}_0+\tilde{b}) \cdot P(\Sigma)\nabla_{x,w}\tilde{f} &= - \int_{D_0} \nabla\psi \cdot P(\Sigma)\nabla_{x,w}\tilde{f} - \gamma \int_{\mathbb{T}} (-h+\beta)\beta_x \tilde{f} |_{w=-h} , \\
&\;\; \forall \tilde{f} \in H^1_{0,surf}(D_0),
\end{align*}
so that $\tilde{a}_0+\tilde{b}$ is a variational solution to \eqref{eq:systphieq1}-\eqref{eq:systphieq3}. By uniqueness of the variational solution, we obtain that $\tilde{a}_0+\tilde{b} = \tilde{\varphi}$, and we can deduce the thesis.
\end{proof}

Next, we prove an estimate for the derivatives of the map $(\eta,\beta) \mapsto \tilde{\varphi}$.

\begin{lemma} \label{lem:TechDerEst}

Let $s > 3/2$, $\gamma \in \mathbb{R}$, $\psi \in H^s(\mathbb{T})$ and $\bmzeta=(\eta,\beta) \in ( H^{s}(\mathbb{T}) )^2$. Moreover, let us choose $\sigma$ of the form  \eqref{eq:straightReg} in the diffeomorphism \eqref{eq:straight} and assume that the condition \eqref{eq:StrConnected} holds true. Then for all $j \in \mathbb{N}$ and $(\bmh,\bmk)= (h_1,\ldots,h_j,k_1,\ldots,k_j) \in ( H^s(\mathbb{T}) )^{2j}$ the solution of the boundary value problem \eqref{eq:systphieq1}-\eqref{eq:systphieq3}  satisfies
\begin{align*}
& \| \nabla d^j \tilde{\varphi}(\bmzeta,\gamma,\psi)(\bmh,\bmk) \|_{ ( H^{s-1/2,1}(D_0) )^2 } \\
&\leq C_1( \|\eta\|_{H^s(\mathbb{T})} ,  \|\beta\|_{H^s(\mathbb{T})} , |\gamma| ) \left[ \prod_{m=1}^j ( \| h_m \|_{ H^s(\mathbb{T}) } + \| k_m \|_{ H^s(\mathbb{T}) } ) \right] \, \|\psi\|_{H^s(\mathbb{T})}  \\
&\quad + C_2( \|\eta\|_{H^s(\mathbb{T})} ,  \|\beta\|_{H^s(\mathbb{T})} , |\gamma| ) .
\end{align*}

\end{lemma}

\begin{proof}

For simplicity we write $d^j\tilde{\varphi}(\bmh,\bmk)$ instead of $d^j \tilde{\varphi}(\bmzeta,\gamma,\psi)(\bmh,\bmk)$.

We differentiate \eqref{eq:systphieq1}-\eqref{eq:systphieq3} with respect to $(\eta,\beta)$: for $j=0$ we have $d^0\tilde{\varphi}(\bmh,\bmk)=\psi$, and we can conclude by Theorem \ref{thm:AnSolBVP};  by arguing as in the proof of Proposition A.7 in \cite{lannes2013water} for $j \geq 1$, we can deduce the thesis.
\end{proof}

Finally, from the above results we can deduce an analyticity result for the operator $G(\eta,\beta,\gamma)$.

\begin{corollary} \label{cor:AnalOpG}

Let $s > 3/2$, $\gamma \in \mathbb{R}$ and $\psi \in H^s(\mathbb{T})$. Moreover, let us choose $\sigma$ of the form  \eqref{eq:straightReg} in the diffeomorphism \eqref{eq:straight}. Then the formula
\begin{align*}
G(\eta,\beta,\gamma)(\psi) &= \nabla^{\Sigma}\tilde{\varphi}(\eta,\beta,\gamma,\psi) \cdot \bmN |_{w=0} 
\end{align*}
defines an analytic map
\begin{align*}
G(\cdot,\cdot,\gamma)(\psi) &: \{ (\eta,\beta) \in H^{s}(\mathbb{T}) \times H^{s}(\mathbb{T}) : \eqref{eq:StrConnected} \; \text{holds true} \} \to H^{s-1}(\mathbb{T}) 
\end{align*}

\end{corollary}

\begin{proof}

Let $\bmzeta=(\eta,\beta)$, and let us write $\mathcal{G}(\bmzeta) \coloneqq G(\bmzeta,\gamma)(\psi)$. We want to prove that for all $\bmzeta_0$ in $( H^{s}(\mathbb{T}) )^2$ we can write for all $\bmzeta$ in a small neighbourhood of $\bmzeta_0$ 
\begin{align*}
\mathcal{G}(\bmzeta) &= \sum_{k =0}^{+\infty} \mathcal{G}_k(\bmzeta-\bmzeta_0)^k ,
\end{align*}
where each $\mathcal{G}_k: ( H^{s}(\mathbb{T}) )^{2k} \to H^s(\mathbb{T})$ is bounded, symmetric, $k$-linear and satisfies
\begin{align*}
\| \mathcal{G}_k \|_{ ( H^{s}(\mathbb{T}) )^{2k} \to H^s(\mathbb{T}) } &\leq N_0 r_0^{-k} , 
\end{align*}
for some $N_0,r_0>0$. 

If such a result is true, then by looking at terms with same degree of homogeneity we have that for all $0 \leq m \leq s-1$
\begin{align} \label{eq:OpGrel1}
\left( \langle \partial_x \rangle^m g_k , f \right) &= \sum_{j=0}^k \int_{D_0} \langle \partial_x \rangle^{m+1/2} ( P_j \nabla \tilde{a}_{k-j} ) \langle \partial_x \rangle^{-1/2} \nabla (\chi(w |D|) f ),
\end{align}
where 
\begin{align*}
g_j = \mathcal{G}_j(\bmzeta-\bmzeta_0)^j , \;\; Q_j &= Q_j(\bmzeta-\bmzeta_0)^j ,\;\; \tilde{a}_j = \tilde{\varphi}_j(\bmzeta-\bmzeta_0)^j ,
\end{align*}
and $P_0 = \mathbb{I}_2+Q_0$, $P_j = Q_j$ $\forall j \geq 1$. Hence, we want to prove that
\begin{align} \label{eq:OpGrel2}
\| g_k \|_{ H^m(\mathbb{T}) } &\leq N_0 r_0^{-k} \, \| \bmzeta - \bmzeta_0 \|_{( H^s(\mathbb{T}) )^2 }^k , \;\; \forall k \in \mathbb{N}.
\end{align}
From \eqref{eq:OpGrel1} we get that for all $0 \leq m \leq s-3/2$
\begin{align*}
\| g_k \|_{ H^m(\mathbb{T}) } &\leq \sum_{j=0}^k \| Q_j \nabla \tilde{a}_{j-k} \|_{ H^{m+1/2,0}(D_0) } \\
&\leq \sum_{j=0}^k \| Q_j \|_{ M_3( L^{\infty}([-h,0]; H^{s-1}(\mathbb{T}) ) )} \| \langle \partial_x \rangle^{m+1/2} \nabla \tilde{a}_{k-j} \|_{( L^2(D_0) )^2 } ,
\end{align*}
but using the properties of the regularizing diffeomorphism and Theorem \ref{thm:AnSolBVP} we have
\begin{align*}
 \| Q_j \|_{ M_3( L^{\infty}([-h,0]; H^{s-1}(\mathbb{T}) ) )} + \| \langle \partial_x \rangle^{m+1/2} \nabla \tilde{a}_{j} \|_{( L^2(D_0) )^2 } &\leq N_1 r_1^{-j} \|\bmzeta-\bmzeta_0\|_{( H^s(\mathbb{T}) )^2}^j ,
\end{align*}
hence arguing as in the proof of Theorem \ref{thm:AnSolBVP} we can deduce the thesis.

\end{proof}

From Corollary \ref{cor:AnalOpG} we can deduce Theorem \ref{thm:AnalOpGMain}.

\subsection{ Homogeneous expansion around $\eta=0$ } \label{subsec:HomExp}

Now we consider the expansion of the operator $G(\eta,\beta,\gamma)$ around $\eta=0$,
\begin{align*}
G(\eta,\beta,\gamma) &= \sum_{j=0}^{\infty} G_j[\eta](\beta,\gamma)
\end{align*}
where $G_j[\eta](\beta,\gamma)$ is homogeneous of degree $j$ in $\eta$. From \eqref{eq:BVPsolphi} we have that for $\eta=0$
\begin{align*}
\tilde{\varphi}_0(\cdot,w) &= \tilde{\varphi}(\cdot,w)|_{r_{1}=r_{2}=0} \\
&= \mathscr{F}^{-1} \, \bigg[ \mathtt{G}_{\tilde{\varphi}}(w,-h) \mathscr{F} [ \gamma (-h+\beta) \, \beta_x  ] - \mathtt{G}_{\tilde{\varphi},\zeta}(w,0) \, (\mathscr{F} \psi) \, \bigg] \\
&=  \mathscr{F}^{-1} \, \bigg[  \gamma \, \frac{ \sinh(\xi \, w) }{ - \xi \, \cosh(h \, \xi) + \xi \, (D\beta) \, \sinh(h \, \xi) } \, \mathscr{F} \left( (-h+\beta) \beta_x \right) \nonumber \\
&\qquad \qquad + \frac{ \cosh(\xi (w+h)) - (D\beta) \, \sinh( \xi (w+h) ) }{ \cosh(h \, \xi) - (D\beta) \, \sinh(h \, \xi) } \, (\mathscr{F} \psi)  \bigg] ,
\end{align*}
so that 
\begin{align}
G_{0}(\beta,\gamma)(\psi) &=  \partial_{w}\tilde{\varphi}_0|_{w=0}  \nonumber \\
&=  \gamma \, \mathscr{F}^{-1} \, \bigg[ \mathtt{G}_{\tilde{\varphi},w }(0,-h) \, \mathscr{F} \left( (-h+\beta) \beta_x \right) \, \bigg]  - \mathscr{F}^{-1} \, \bigg[ \mathtt{G}_{\tilde{\varphi},w \zeta}(0,0) \, (\mathscr{F} \psi) \, \bigg] \nonumber \\
&= -\gamma \, \mathscr{F}^{-1} \, \bigg[   \frac{ 1 }{  \cosh(h \, \xi) - \, (D\beta) \, \sinh(h \, \xi) } \, \mathscr{F} \left( (-h+\beta) \beta_x \right) \bigg] \nonumber \\
&\;\; + \mathscr{F}^{-1} \, \bigg[ \xi \,  \frac{ \sinh(h \, \xi) - (D\beta) \, \cosh( h \, \xi ) }{ \cosh(h \, \xi) - (D\beta) \, \sinh(h \, \xi) } \, (\mathscr{F} \psi) \, \bigg]  .  \label{eq:G0}
\end{align}

\begin{remark} \label{rem:IrrDNO}

In the irrotational case $\gamma=0$ we have
\begin{align}
& G_{0}(\beta,0)(\psi) \nonumber \\
&=   \mathscr{F}^{-1} \, \left\{ \left[ \xi \, \tanh(h \, \xi) - \xi \, \frac{ D\beta }{ \cosh(h \, \xi) \, ( \cosh(h \, \xi) - (D\beta) \, \sinh(h \, \xi) ) } \right] \, (\mathscr{F} \psi) \, \right\}  \nonumber \\
&=: ( D \, \tanh(h D) + D \, L(\beta) ) \psi , \label{eq:G0Irr}
\end{align}
which recovers formula \emph{(A.1)} in \cite{craig2005hamiltonian}. 

Using Remark \ref{rem:IrrDNO}, we can rewrite the operator $G_{0}(\beta,\gamma)$ defined in \eqref{eq:G0} in the following way,
\begin{align}
G_{0}(\beta,\gamma)(\psi) &= ( D \, \tanh(h D) + D \, L(\beta) ) \psi + \gamma \, \nu(\beta) , \label{eq:G0new} \\
\nu(\beta) &= - \mathscr{F}^{-1} \, \bigg[   \frac{ 1 }{  \cosh(h \, \xi) - \, (D\beta) \, \sinh(h \, \xi) } \, \mathscr{F} \left( (-h+\beta) \beta_x \right) \bigg] . \nonumber
\end{align}

\end{remark}

\begin{remark} \label{rem:HigherTerms}

Using Proposition \ref{prop:SysOpDer} we have 
\begin{align} 
\partial_{\eta} G^{DN}(\eta,\beta,\gamma)(\delta\eta)(\psi) &= - G^{DN}(\eta,\beta) \left[ \delta\eta \, B(\eta,\beta,\gamma)(\psi)  \right] - \partial_x  \left[ \delta\eta \, V(\eta,\beta,\gamma)( \psi ) \right] , \label{eq:ShapeDerOpG}
\end{align}
where
\begin{align*}
V(\eta,\beta,\gamma)(\psi) &\coloneqq \psi_x - B(\eta,\beta,\gamma)(\psi) \, \eta_x , \\
B(\eta,\beta,\gamma)(\psi ) &\coloneqq \frac{ G(\eta,\beta,\gamma)(\psi) + \eta_x \psi_x  }{1+\eta_x^2} .
\end{align*}
By plugging into \eqref{eq:ShapeDerOpG} the expansions
\begin{align*}
G(\eta,\beta,\gamma) &= \sum_{j=0}^{\infty} G_j[\eta](\beta,\gamma), \\
V(\eta,\beta,\gamma) &= \sum_{j=0}^{\infty} V_j[\eta](\beta,\gamma), \, B(\eta,\beta,\gamma) = \sum_{j=0}^{\infty} B_j[\eta](\beta,\gamma),
\end{align*}
where $V_j[\eta](\beta,\gamma)$ and $B_j[\eta](\beta,\gamma)$ are homogeneous of degree $j$ in $\eta$, we can deduce recursive formulae for the terms $G_j[\eta](\beta,\gamma)$, for all $j \geq 1$ (see Sec 2.4 of \cite{groves2024analytical} and Sec. 6 of \cite{pasquali2025analytical} for similar arguments for three-dimensional water waves).

\end{remark}

\subsection{Paralinearization of the operator $G(\eta,\beta,\gamma)$} \label{subsec:ParaOpG}

Recalling \eqref{eq:OpG}, we have that if $s > 3/2$, $\beta \in H^s(\mathbb{T})$, $\gamma \in \mathbb{R}$, $(\eta,\psi) \in H^s(\mathbb{T}) \times H^s(\mathbb{T})$ and if \eqref{eq:StrConnected} holds true, then by Corollary \ref{cor:AnalOpG} and by Proposition 2.3 in \cite{lannes2013water} we can deduce the estimate
\begin{align} 
\| G(\eta,\beta,\gamma)(\psi) \|_{ H^{s-1}(\mathbb{T}) } &\leq C_1 \left( \| \eta \|_{ H^s(\mathbb{T}) } , \| \beta \|_{ H^s(\mathbb{T}) }  \right) \, \| \psi \|_{ H^s(\mathbb{T}) } \nonumber \\
&\quad +|\gamma| \, C_2 \left( \| \eta \|_{ H^s(\mathbb{T}) } , \| \beta \|_{ H^s(\mathbb{T}) }  \right) . \label{eq:EstOpG}
\end{align}

By using the estimate \eqref{eq:EstOpG} we can deduce the following lemma.

\begin{lemma} \label{lem:EstGU}

Let $s > 3/2$, $\beta \in H^s(\mathbb{T})$, $\gamma \in \mathbb{R}$ and $(\eta,\psi) \in H^s(\mathbb{T}) \times H^s(\mathbb{T})$. Assume that \eqref{eq:StrConnected} holds true, then 
\begin{align}
B(\eta,\beta,\gamma)(\psi) &= \frac{ G(\eta,\beta,\gamma)(\psi) + \eta_x \psi_x }{1+\eta_x^2} \in H^{s-1}(\mathbb{T}) , \nonumber \\
V(\eta,\beta,\gamma)(\psi) &= \psi_x - B(\eta,\beta,\gamma)(\psi) \, \eta_x \in H^{s-1}(\mathbb{T}), \nonumber \\
\omega(\eta,\beta,\gamma)(\psi) &= \psi - T_{B(\eta,\beta,\gamma)(\psi)} \eta \in H^s(\mathbb{T}) , \label{eq:defGU}
\end{align}
and
\begin{align} 
& \| B(\eta,\beta,\gamma)(\psi) \|_{ H^{s-1}(\mathbb{T}) } + \| V(\eta,\beta,\gamma)(\psi) \|_{ H^{s-1}(\mathbb{T}) } + \| \omega(\eta,\beta,\gamma)(\psi) \|_{H^s(\mathbb{T})} \nonumber \\
&\leq C_1 \left( \| \eta \|_{ H^s(\mathbb{T}) } , \| \beta \|_{ H^s(\mathbb{T}) }  \right) \, \| \psi \|_{ H^s(\mathbb{T}) } + |\gamma| \, C_2 \left( \| \eta \|_{ H^s(\mathbb{T}) } , \| \beta \|_{ H^s(\mathbb{T}) }  \right) . \label{eq:EstGUSob}
\end{align}

\end{lemma}	

\begin{proof}

The estimate \eqref{eq:EstGUSob} for $B(\eta,\beta,\gamma)(\psi)$ and for $V(\eta,\beta,\gamma)(\psi)$ follow from classical product estimates for Sobolev spaces. Moreover, the Sobolev embedding implies that $B(\eta,\beta,\gamma)(\psi) \in L^{\infty}(\mathbb{T})$; hence, by Lemma \ref{lem:Paraprod1} we obtain
\begin{align}
\| T_{B(\eta,\beta,\gamma)(\psi)} \eta  \|_{H^s(\mathbb{T})} &\leq C \, \| B(\eta,\beta,\gamma)(\psi) \|_{L^{\infty}(\mathbb{T})} \, \|\eta\|_{H^s(\mathbb{T})} \nonumber \\
&\leq C_1 \left( \| \eta \|_{ H^s(\mathbb{T}) } , \| \beta \|_{ H^s(\mathbb{T}) }  \right)  \, \| \psi \|_{ H^s(\mathbb{T}) }  \, \| \eta \|_{ H^s(\mathbb{T}) } \nonumber \\
&\quad + |\gamma| \, C_2 \left( \| \eta \|_{ H^s(\mathbb{T}) } , \| \beta \|_{ H^s(\mathbb{T}) }  \right)  \, \| \eta \|_{ H^s(\mathbb{T}) } , \label{eq:EstTBeta}
\end{align}
which allows to deduce estimate \eqref{eq:EstGUSob} for $\omega(\eta,\beta,\gamma)(\psi)$.

\end{proof}

Moreover, by exploiting \eqref{eq:EstTBeta} we also have the following result (see lemma 2.2 in \cite{alazard2018control} for an analogous result for the flat bottom topography case).

\begin{lemma} \label{lem:InvGU}

Let $s > 3/2$, $\eta,\beta \in H^s(\mathbb{T})$. Assume that \eqref{eq:StrConnected} holds true. Then there exists $\epsilon_0>0$ such that if 
\begin{align} \label{eq:SmallAss}
\| \eta \|_{ H^s(\mathbb{T}) }  &< \epsilon_0 ,
\end{align}
then there exists an operator $\Psi(\eta,\beta)$ such that 
\begin{align*}
\Psi(\eta,\beta) \left( \omega(\eta,\beta,0)(\psi) \right) &= \psi, \;\; \forall \psi \in H^{1/2}(\mathbb{T}) .
\end{align*}
Moreover, if $\omega \in H^s(\mathbb{T})$, then $\Psi(\eta,\beta)(\omega) \in H^s(\mathbb{T})$, and
\begin{align*}
\| \Psi(\eta,\beta)(\omega) \|_{H^s(\mathbb{T})} &\leq C \left( \| \eta \|_{ H^s(\mathbb{T}) } , \| \beta \|_{ H^s(\mathbb{T}) } \right) \, \| \omega \|_{ H^s(\mathbb{T}) } .
\end{align*}

\end{lemma}

By Lemma \ref{lem:InvGU} we have that under the smallness assumption \eqref{eq:SmallAss}, then $B=B(\eta,\beta,0)(\psi)$ and $V=V(\eta,\beta,0)(\psi)$ can be rewritten in terms of $\eta$, $\beta$ and $\omega$,
\begin{align*}
B = B(\eta,\beta,0) \left( \Psi(\eta,\beta)\omega \right), &\;\; V = V(\eta,\beta,0) \left( \Psi(\eta,\beta)\omega \right) .
\end{align*}

Now we state the main result of the section.

\begin{proposition} \label{prop:OpGpara}

Let $s > 5/2$. If $\eta, \beta \in H^{s+1/2}(\mathbb{T})$ are such that \eqref{eq:StrConnected} holds true, if $\gamma \in \mathbb{R}$ and if $\psi \in H^s(\mathbb{T}) $, then  
\begin{align} \label{eq:OpGpara}
G(\eta,\beta,\gamma)(\psi) &= T_{\lambda} \omega  - T_{ V(\eta,\beta,\gamma)(\psi) } \eta_x + \mathcal{R}(\eta,\beta,\gamma,\psi) ,
\end{align}
where the function $\omega$ defined in \eqref{eq:defGU} is Alinhac's good unknown, $\lambda$ is a symbol given by 
\begin{align*}
\lambda(x,\xi) &= \lambda^{(1)} + \lambda^{(0)}(x,\xi) , \\
\lambda^{(1)} &= |\xi| ,\\
\lambda^{(0)}(x,\xi) &= \frac{ 1+\eta_x^2 }{ 2\lambda^{(1)} }  \left( ( \eta_x + \mathrm{i} \partial_{\xi}\lambda^{(1)} ) \, \partial_x + \eta_{xx} \right) \left( \frac{ \mathrm{i} \, \xi \, \eta_x + \lambda^{(1)} }{ 1+\eta_x^2 } \right) , 
\end{align*}
and $\mathcal{R}(\eta,\beta,\gamma,\psi) \in H^{s+1/2}(\mathbb{T})$ satisfies
\begin{align*}
\| \mathcal{R}(\eta,\beta,\gamma,\psi) \|_{ H^{s+1/2}(\mathbb{T}) } &\leq C_1 \left( \| \eta \|_{ H^{s+1/2}(\mathbb{T}) } , \| \beta \|_{ H^{s+1/2}(\mathbb{T}) }    \right) \, \| \psi \|_{ H^{s}(\mathbb{T}) }  \\
&\quad + |\gamma| \, C_2 \left( \| \eta \|_{ H^{s+1/2}(\mathbb{T}) } , \| \beta \|_{ H^{s+1/2}(\mathbb{T}) }    \right) ,
\end{align*}
for some non-decreasing functions $C_1, C_2$.
\end{proposition}

We defer the proof of the above result to Sec. \ref{subsec:proofOpGpara}; we point out that the principal symbol of $\lambda$ in \eqref{eq:OpGpara} is given by the principal symbol of the classical Dirichlet--Neumann operator (see the proof of Lemma \ref{lem:factorisation}). We mention Proposition 4.12 of \cite{alazard2009paralinearization} and Lemma 3.14 of \cite{alazard2011water} regarding analogous results for the case $\gamma=0$.

\subsection{Proof of Proposition \ref{prop:OpGpara} } \label{subsec:proofOpGpara}

In this section we prove Proposition \ref{prop:OpGpara}, namely we paralinearize the operator $G(\eta,\beta,\gamma)$ defined in \eqref{eq:OpG}. We recall that
\begin{align*}
G(\eta,\beta,\gamma)(\psi) &= \nabla\varphi \cdot \bmN |_{y=\eta(x)},
\end{align*}
where $\varphi$ is the solution of the boundary value problem \ref{eq:EllBVP}, namely
\begin{equation*} 
\begin{cases}
\Delta \varphi = 0 \; \; &\text{in} \; \; D_{\eta,\beta}, \\
\varphi = \psi \; \; &\text{at} \; \; y=\eta, \\
\gamma \, y \, \beta_x - \varphi_x \, \beta_x + \varphi_y = 0  \; \; &\text{at} \; \; y=-h+\beta .
\end{cases}
\end{equation*}

We follow the strategy of Sec. 3 of \cite{groves2024analytical} and Sec. 3 of \cite{alazard2011water} (see also Sec. 4 of \cite{alazard2009paralinearization}). In the following we use the notations introduced in Appendix \ref{sec:Paradiff}.

First of all, we start by paralinearizing the interior equation in \eqref{eq:EllBVP}. We want to recast the interior equation in the fixed domain $D_0$ by using a suitable localizing diffeomorphism (which differs from the straightening diffeomorphism defined in \eqref{eq:straight}): we choose $\delta >0$ such that the fluid domain $D_\eta$ contains the strip
\begin{align} \label{eq:strip}
\Omega_{\delta} &:= \{ (x,y) \in \mathbb{T} \times \mathbb{R}| \eta(x)-\delta h \leq y < \eta(x) \} .
\end{align}

\begin{remark} \label{rem:locTrans}

Because of \eqref{eq:StrConnected} we can take any $\delta \in (0, h_0/h)$ in \eqref{eq:strip}: indeed, $y \geq \eta(x)-\delta h > - h+\beta(x)$.
\end{remark}

We define $\hat\Sigma:D_0 \to \Omega_{\delta}$,
\begin{align*} 
	\hat\Sigma: (x,w) \mapsto (x,\varrho(x,w)), &\; \; \varrho(x,w):=\delta w + \eta(x),
\end{align*}
and write $\hat\varphi(x,w) := \varphi \left( x, \varrho(x,w) \right)$, where $\varphi$ is the solution of the boundary value problem \eqref{eq:EllBVP} restricted to $\Omega_{\delta}$.

We denote by $\grad_{x,y}=(\pd_x,\pd_y)^T$, and
\begin{align*}
	\grad_{\bmx,w}^\varrho := \left( J_{\hat\Sigma}^{-1} \right)^T \grad_{x,y}, &\; \; \left( J_{\hat\Sigma}^{-1} \right)^T :=
	\begin{pmatrix}
		1 & - \frac{\eta_x}{\delta} \\
		0 & \frac{1}{\delta} \\
	\end{pmatrix}, \\
\pd_x^\varrho := \pd_x - \frac{ \eta_x }{\delta} \pd_w, &\quad \pd_w^\varrho := \frac{ \pd_w }{\delta} .
\end{align*}

The flattened version of the operator $\Delta$ applied to $\hat{f}(x,w)=f(x,y)$ is given by
\begin{align*}
-\Delta^\varrho \hat{f} &= -\Delta \hat{f} + \frac{2}{\delta} \eta_x \; \pd_{xw}^2\hat{f}  - \left( \frac{1}{\delta^2} - 1\right) \pd_w^2\hat{f} + \frac{1}{\delta} \, \eta_{xx} \, \pd_w\hat{f} - \frac{1}{\delta^2} \, \eta_x^2 \, \pd_w^2\hat{f} .
\end{align*}

Therefore, if we denote by $\varphi$ the solution of \eqref{eq:EllBVP}, we have that $\hat{\varphi} = \varphi \circ \hat{\Sigma}$ satisfies

\begin{align} 
E \hat{\varphi} &= 0, \;\; \text{in} \; \; \Omega_{\delta} , \label{eq:StrIntEq} \\
E &:=  \check{\mathtt{a}} \partial_w^2  + \partial_x^2 + \check{\mathtt{b}} \, \partial_x \partial_w - \check{\mathtt{c}} \, \partial_w  , \nonumber 
\end{align}
where 
\begin{align*}
\check{\mathtt{a}} = \frac{1+\eta_x^2}{\delta^2} , \quad \check{\mathtt{b}} &= -2 \frac{\eta_x}{\delta} , \quad \check{\mathtt{c}} = \frac{\eta_{xx}}{\delta} ,
\end{align*}
while
\begin{align}
\hat{\varphi} &= \psi, \; \; \text{at} \; \; w=0. \label{eq:StrBC1} 
\end{align}

\begin{remark} \label{rem:StrIntEq}

We point out that, due to the fact that the bottom of the fluid domain is sufficiently regular (namely, $\beta \in H^{s+1/2}(\mathbb{T})$ ), we can restrict the interior equation in \eqref{eq:EllBVP} to the strip $\Omega_{\delta}$ so that \eqref{eq:StrIntEq} does not present a forcing term depending on $\beta$ in the right-hand side; see a similar approach in three-dimensional Beltrami flows in Eq. (3.1) of \cite{groves2024analytical}. By using a more careful localization argument, one could allow less regular bottom topographies, obtaining a forcing bottom-dependent term in the right-hand side of the straightened interior equation (see for example Sec. 2.2 and Sec. 3.4 of \cite{alazard2011water}, where the authors consider the irrotational case without any regularity assumption on the bottom).

\end{remark}

\begin{remark} \label{rem:Alinhac}

From \eqref{eq:OpG} we have
\begin{align*}
G(\eta,\beta,\gamma)(\psi) &= - \eta_x \partial_x \hat{\varphi} + \frac{1+\eta_x^2}{\delta} \, \partial_w \hat{\varphi} |_{w=0} ,
\end{align*}
and
\begin{align*}
\omega(\eta,\beta,\gamma)(\psi) = \psi - T_{ B(\eta,\beta,\gamma)(\psi) } \eta &= \hat\varphi - T_{\delta^{-1} \partial_w \hat{\varphi} } \varrho |_{w=0} ,
\end{align*}
hence the good unknown
\begin{align*}
\hat{\Phi} &= \hat{\varphi} - T_{ \partial_w^{\varrho} \hat{\varphi} }\eta ,
\end{align*}
satisfies $\hat{\Phi}|_{w=0} = \psi - T_{ B(\eta,\beta,\gamma)(\psi) } \eta$.

\end{remark}

\begin{proposition} \label{prop:paraGood}

Set
\begin{align*} 
\hat{\Phi} &= \hat{\varphi} - T_{\partial_w^\varrho \hat{\varphi}} \eta, 
\end{align*}
and let $T_E:= ( \delta^{-2} \mathrm{Id} + T_{ \check{\mathtt{a}} -1/\delta^2 } ) \partial_w^2 +  \partial_x^2 + T_{\check{\mathtt{b}}}  \partial_x \, \partial_w - T_{ \check{\mathtt{c}} } \,\partial_w$, then $\hat{\Phi}$ satisfies 
\begin{align} \label{eq:ParaIntEq}
T_{E} \hat{\Phi} &= \hat{F}_0 + \hat{F}_1 ,
\end{align}
where $\hat{F}_0 \in C^1( [-h,0];H^{s+1/2}(\mathbb{T}) ) $ and $\hat{F}_1 \in C^0( [-h,0];H^{2s-3}(\mathbb{T}) ) $ . 

\end{proposition}

\begin{proof}

We follow the strategy of the proof of Proposition 4.12 in \cite{alazard2009paralinearization}. First, due to the regularity assumptions on $\eta$ and $\psi$ we have that
\begin{align*}
\partial_w^k \hat{\varphi} &\in  C^0( [-h,0] ; H^{s-k}(\mathbb{T}) )  , \;\; k=1,2,3,
\end{align*}
and that by the regularity assumption of $\eta$ we have
\begin{align*}
\check{\mathtt{a}}-\frac{1}{\delta^2} &\in H^{s-1/2}(\mathbb{T})  , \quad \check{\mathtt{b}} \in   H^{s-1/2}(\mathbb{T})  , \quad \check{\mathtt{c}} \in H^{s-3/2}(\mathbb{T})  .
\end{align*}

By Lemma \ref{lem:Paraprod2} we have that if we set
\begin{align*}
\mathcal{T}(\hat{\varphi}) &:= T_{ \partial_w^2 \hat{\varphi} } \, (\check{\mathtt{a}}-1/\delta^2)  + T_{ \partial_x \, \partial_w \hat{\varphi} } \cdot \check{\mathtt{b}} - T_{ \partial_w \hat{\varphi} } \, \check{\mathtt{c}} ,
\end{align*}
then
\begin{align*}
E \hat{\varphi} - T_{E} \hat{\varphi} - \mathcal{T}(\hat{\varphi}) &\in C^0([-h,0];H^{2s-3}(\mathbb{T})) .
\end{align*}

Next, we substitute $\hat{\varphi} = \hat{\Phi} + T_{\partial_w^\varrho \hat{\varphi}} \eta$ in the interior equation \eqref{eq:StrIntEq}; we get (up to a smoother remainder)

\begin{align*}
\partial_w^2 \hat{\Phi} &=  \partial_w^2 \hat{\varphi} - T_{\partial_w^2 \partial_w^\varrho \hat{\varphi}} \eta  , \\
\partial_x^2 \hat{\Phi} &= \partial_x^2 \hat{\varphi} - T_{\partial_x^2 \partial_w^\varrho \hat{\varphi}} \eta -2 T_{ \partial_x \partial_w^\varrho \hat{\varphi}} \, \eta_x - T_{\partial_w^\varrho \hat{\varphi}} \eta_{xx} , \\
\partial_x \partial_w \hat{\Phi} &= \partial_x \partial_w \hat{\varphi} - T_{\partial_x \partial_w \partial_w^\varrho \hat{\varphi}} \eta - T_{ \partial_w \partial_w^\varrho \hat{\varphi}} \eta_x , \\
\partial_w \hat{\Phi} &= \partial_w \hat{\varphi} - T_{\partial_w \partial_w^\varrho \hat{\varphi}} \eta ,
\end{align*}
hence by \eqref{eq:StrIntEq} and by Lemma \ref{lem:Paraprod2} we get 
\begin{align*}
T_{E} \hat{\Phi} + T_E T_{\partial_w^{\varrho} \hat{\varphi}}\eta + \mathcal{T}(\hat{\varphi}) &= g \in C^0( [-h,0]; H^{s+1/2}(\mathbb{T}) )  ,
\end{align*}
up to a remainder in $C^0( [-h,0]; H^{2s-3}(\mathbb{T}) )$. Notice that
\begin{align*}
T_E T_{\partial_w^{\varrho} \hat{\varphi}}\eta &= T_{ E \partial_w^{\varrho} \hat{\varphi}}\eta + 2 T_{ \partial_x \partial_w^{\varrho} \hat{\varphi} } \eta_x  + T_{ \partial_w^{\varrho} \hat{\varphi}  } \eta_{xx} + T_{ \check{\mathtt{b}} } T_{ \partial_w \partial_w^{\varrho} \hat{\varphi}  } \eta_{x},
\end{align*}
up to a remainder in $C^0( [-h,0]; H^{2s-3}(\mathbb{T}) )$. Next, observe that $s > \frac{5}{2}$ implies
\begin{align*}
T_{E  \partial_w^{\varrho} \hat{\varphi} } \eta &\in C^0([-h,0] ; H^{s+1/2}(\mathbb{T}) ) .
\end{align*}
Since
\begin{align*}
T_E T_{\partial_w^{\varrho} \hat{\varphi}}\eta + \mathcal{T}(\hat{\varphi}) &\in C^0([-h,0] ; H^{s+1/2}(\mathbb{T}) ),
\end{align*}
up to a remainder in $C^0( [-h,0]; H^{2s-3}(\mathbb{T}) )$, we can deduce the thesis.

\end{proof}

Next we consider the reduction to the boundary. We want to find two symbols $\mathtt{M}$, $\mathtt{N}$ such that (up to a remainder)
\begin{align*}
T_{E}  &= \check{\mathtt{a}} \left( \partial_w  - T_{ \mathtt{N} } \right) \left( \partial_w  - T_{ \mathtt{M} } \right),
\end{align*}
see Sec. 3.4.2 of \cite{alazard2011water} and Sec. 4.3 of \cite{alazard2009paralinearization} for a corresponding argument for the irrotational case, and Sec. 3.1 of \cite{groves2024analytical} for the case of Beltrami fields.

\begin{lemma} \label{lem:factorisation}
Let $\eta \in H^{s+1/2}(\mathbb{T})$, let $\sigma = \min \left( \frac{1}{2} , s-\frac{5}{2} \right) >0$. Then there exist operators $M^{(1)}$, $N^{(1)}$, $M^{(0)}$, $N^{(0)}$ with associated symbols 
\begin{equation*}
\mathtt{M}^{(1)},\mathtt{N}^{(1)} \in  \Gamma^1_{3/2+\sigma}(\mathbb{T}) , \quad \mathtt{M}^{(0)},\mathtt{N}^{(0)} \in  \Gamma^0_{1/2+\sigma}(\mathbb{T}) ,
\end{equation*}
such that if we write 
\begin{equation*}
\mathtt{M} = \mathtt{M}^{(1)} + \mathtt{M}^{(0)}, \;\; \mathtt{N} = \mathtt{N}^{(1)} + \mathtt{N}^{(0)}, \;\; M = M^{(1)} + M^{(0)}, \;\; N = N^{(1)} + N^{(0)},
\end{equation*}
then
\begin{itemize}
\item[(i)] there exist operators $R_0$ and $R_1$ of order $\frac{1}{2}-\sigma$ and $-\frac{1}{2}-\sigma$ respectively, such that
\begin{equation*}
T_{E} = \check{\mathtt{a}} (\partial_w  - T_{ \mathtt{N} } )(\partial_w  - T_{ \mathtt{M} } ) + R_0 + R_1 \partial_w ,
\end{equation*}
\item[(ii)]
the symbols $\mathtt{M}^{(1)}$, $-\mathtt{N}^{(1)}$ are strongly elliptic, namely there exist $c,C>0$ such that
\begin{align*}
\mathrm{Re} \, \mathtt{M}^{(1)} &\geq c \langle \xi \rangle , \quad \forall | \xi | \geq C, \\
\mathrm{Re} \, ( -\mathtt{N}^{(1)} ) &\geq c \langle \xi \rangle , \quad \forall | \xi | \geq C .
\end{align*}

\end{itemize}
\end{lemma}

\begin{remark}

Actually, arguing as in Sec. 4 of \cite{alazard2009paralinearization} one can prove that if $s \in \mathbb{N}$, then $\mathtt{M}^{(1)}, \mathtt{N}^{(1)} \in \Sigma^1_{s-1}(\mathbb{T}) $.

\end{remark}

\begin{proof}

\begin{itemize}
\item[$\mathrm{(i)}$] We first set $ L_1 :=  \check{\mathtt{b}} \partial_x -  \check{\mathtt{c}}   $, $L_0 :=  \partial_x^2$, so that $E = \check{\mathtt{a}} \partial_w^2 + L_1 \, \partial_w + L_0$.  Recall that $\eta \in H^{s+1/2}(\mathbb{T})$. Because
\begin{align*}
E - \check{\mathtt{a}} (\partial_w - N)(\partial_w  - M ) &= (L_1+ \check{\mathtt{a}} (M +N)) \partial_w + (L_0 - \check{\mathtt{a}} N M)
\end{align*}
we set
\begin{align*}
N &= - \check{\mathtt{a}}^{-1} L_1 - M
\end{align*}
and, by following Proposition \ref{prop:AsympParadiffOp}, we seek $M$ such that
\begin{align*}
L_0 + L_1 M + \check{\mathtt{a}} M^2 &= 0 
\end{align*}
by constructing a symbol $\mathtt{M}=\mathtt{M}^{(1)} + \mathtt{M}^{(0)}$ such that 
\begin{align*}
\mathtt{M}^{(1)} \in  \Gamma^1_{3/2+\sigma}(\mathbb{T}) , &\quad \mathtt{M}^{(0)} \in  \Gamma^0_{1/2+\sigma}(\mathbb{T}) ,
\end{align*}
and such that
\begin{align*}
&  -  \xi^2  + \left[ \mathrm{i} \, \check{\mathtt{b}}  \xi - \check{\mathtt{c}} \right] \, \mathtt{M}   +  \check{\mathtt{b}} \, \partial_x \mathtt{M}  + \check{\mathtt{a}} \sum_{\alpha \in \mathbb{N} } \partial_{\xi}^{\alpha}  \mathtt{M} \, D_{x}^{\alpha}  \mathtt{M} = 0 ,
\end{align*}
up to smoother remainders.

We proceed by computing the terms in the expansion of $\mathtt{M}$.

\begin{itemize}
\item
At first order we have 
\begin{align*}
-  \xi^2  + \mathrm{i} \check{\mathtt{b}} \xi \mathtt{M}^{(1)} + \check{\mathtt{a}} (\mathtt{M}^{(1)})^2 &= 0 ,
\end{align*}
so that
\begin{align*}
\mathtt{M}^{(1)}:= \frac{1}{ 2 \check{\mathtt{a}} } \left[  - \mathrm{i} \check{\mathtt{b}} \, \xi +  \sqrt{ - ( \check{\mathtt{b}} \, \xi )^2 + 4 \check{\mathtt{a}} \, \xi^2  } \right] .
\end{align*}
Note that 
\begin{align*}
\mathtt{M}^{(1)} &= \delta \, \frac{ \mathrm{i} \, \xi \, \eta_x + \lambda^{(1)} }{ 1+\eta_x^2 }, \;\; \lambda^{(1)} \, := \, |\xi| ,
\end{align*}
where $\lambda^{(1)}$ is the leading order symbol of the classical Dirichlet--Neumann operator.
\item
The subprincipal symbol of $\mathtt{M}$ is found from the equation
\begin{align*}
& - \check{\mathtt{c}} \, \mathtt{M}^{(1)}  +\mathrm{i} \check{\mathtt{b}} \, \xi \, \mathtt{M}^{(0)} +  \check{\mathtt{a}} \mathtt{M}^{(0)}\mathtt{M}^{(1)} + \check{\mathtt{a}} \mathtt{M}^{(1)}\mathtt{M}^{(0)} \\
&\;\;\;\; + \check{\mathtt{b}} \, \partial_x \mathtt{M}^{(1)} - \check{\mathtt{a}} \mathrm{i}  \partial_{\xi} \mathtt{M}^{(1)} \partial_x \mathtt{M}^{(1)}  = 0,
\end{align*}
which yields
\begin{align*}
\mathtt{M}^{(0)} &= \frac{1}{ \mathrm{i} \check{\mathtt{b}} \, \xi + 2 \check{\mathtt{a}} \mathtt{M}^{(1)} }  \left[ \mathrm{i} \, \check{\mathtt{a}} \partial_{\xi} \mathtt{M}^{(1)} \, \partial_x \mathtt{M}^{(1)} -  \check{\mathtt{b}} \, \partial_x \mathtt{M}^{(1)} + \check{\mathtt{c}} \, \mathtt{M}^{(1)}  \right] .
\end{align*}
Observe that
\begin{align*}
\mathrm{i} \check{\mathtt{b}} \, \xi + 2 \check{\mathtt{a}} \mathtt{M}^{(1)} = 2\delta^{-1} \lambda^{(1)}  , &\quad \mathrm{i} \, \check{\mathtt{a}} \partial_{\xi} \mathtt{M}^{(1)}  -  \check{\mathtt{b}} = \delta^{-1} ( \eta_x + \mathrm{i} \partial_{\xi} \lambda^{(1)} ) ,
\end{align*}
hence
\begin{align*}
\mathtt{M}^{(0)} &= \delta \, \frac{ 1 }{ 2\lambda^{(1)} }  \left( ( \eta_x + \mathrm{i} \partial_{\xi}\lambda^{(1)} ) \, \partial_x + \eta_{xx} \right) \left( \frac{ \mathrm{i} \, \xi \, \eta_x + \lambda^{(1)} }{ 1+\eta_x^2 } \right) .
\end{align*}

\end{itemize}

Similarly, we find that
\begin{align*}
\mathtt{N}^{(1)} &:= \frac{1}{2} \left[  - \mathrm{i} \check{\mathtt{b}} \, \xi -  \sqrt{ - ( \check{\mathtt{b}} \, \xi )^2 + 4 \check{\mathtt{a}} \, \xi^2  } \right] ,
\end{align*}
so that
\begin{align*}
\mathtt{N}^{(1)}  &= \delta \, \frac{ \mathrm{i} \xi \, \eta_x - \lambda^{(1)} }{1+\eta_x^2} ,
\end{align*}
while $\mathtt{N}^{(0)}$ can be written explicitly arguing as for $\mathtt{M}^{(0)}$.

Moreover, we have that 
\begin{align*}
R_0 &\coloneqq \check{\mathtt{a}}  T_{\mathtt{N}} T_{\mathtt{M}} -  \partial_x^2
\end{align*}
is of order $2 - \frac{3}{2} - \sigma = \frac{1}{2}-\sigma$, while
\begin{align*}
R_1 &\coloneqq - \check{\mathtt{a}} (T_{\mathtt{N}} + T_{\mathtt{M}}) + T_{ \check{\mathtt{b}} } \partial_x - T_{ \check{\mathtt{c}} }
\end{align*}
is of order $1 - \frac{3}{2} - \sigma = -\frac{1}{2}-\sigma$.

\item[$\mathrm{(ii)}$] Finally, observe that
$$\mathrm{Re}\, \mathtt{M}^{(1)}  = - \mathrm{Re}\, \mathtt{N}^{(1)}  = \delta \frac{ \lambda^{(1)}}{1+\eta_x^2} \gtrsim \langle \xi \rangle$$
for sufficiently large $|\xi|$, so that $\mathtt{M}^{(1)} $ and $-\mathtt{N}^{(1)}$ are strongly elliptic.
\end{itemize}

\end{proof}

Next, we state an elliptic regularity result (for the proof see Proposition 3.19 in \cite{alazard2011water}; see also Proposition 4.13 in \cite{alazard2009paralinearization}).

\begin{lemma} \label{lemma:EqGoodUnknown}

Let $\eta \in H^{s+1/2}(\mathbb{T})$, let $\sigma = \min \left( \frac{1}{2} , s-\frac{5}{2} \right) >0$. 

Introduce $\hat{W} := (\partial_w-T_{\mathtt{M}})\hat{\Phi}$, then
\begin{align*}
\partial_w \hat{W} - T_{ N^{(1)} } \hat{W} &= T_{ N^{(0)} } \hat{W} + \hat{G}_0 ,
\end{align*}
where $\hat{G}_0 \in  C^0([-h,0];H^{s+\sigma-1/2}(\mathbb{T})) $. Moreover,
\begin{align} \label{eq:ParaGoodEq} 
\partial_w \hat{\Phi} - T_{\mathtt{M}} \hat{\Phi}|_{w=0} &\in H^{s+1/2}(\mathbb{T}). 
\end{align}

\end{lemma}

Now, in order to simplify the notation we omit until the end of the proof the fact that all functions depending on $w$ are evaluated on $w=0$. Recall that 
\begin{align*}
G(\eta,\beta,\gamma)(\psi) &= - \eta_x \, \partial_x^{\varrho}\hat{\varphi} + \partial_w^{\varrho}\hat{\varphi} |_{w=0} \, = \, - \eta_x \, \partial_x\hat{\varphi} +  \frac{1+\eta_x^2}{\delta}  \partial_w\hat{\varphi} |_{w=0} .
\end{align*}

By Lemma \ref{lem:Paraprod2} we obtain that
\begin{align*}
& - T_{ \eta_x } \, \partial_x\hat{\varphi}  - T_{ \partial_x\hat{\varphi} } \, \eta_x +  T_{   \frac{1+\eta_x^2}{\delta} } \partial_w\hat{\varphi}  + T_{  \partial_w\hat{\varphi} } \frac{1+\eta_x^2}{\delta} - G(\eta,\beta,\gamma)(\psi) \in H^{2s-2}(\mathbb{T}) .
\end{align*}

Now, introducing the good unknown $\hat{\Phi} = \hat{\varphi} - T_{ \partial_w^{\varrho} \hat{\varphi} } \eta$ we get 
\begin{align*}
& - T_{ \eta_x } \, \partial_x\hat{\Phi} + T_{ \eta_x } \,  T_{ \partial_x \partial_w^{\varrho} \hat{\varphi} } \eta  + T_{ \eta_x } \,  T_{ \partial_w^{\varrho} \hat{\varphi} } \eta_x  - T_{ \partial_x\hat{\varphi} } \, \eta_x \\
&\qquad +  T_{ \frac{1+\eta_x^2}{\delta} } \partial_w\hat{\Phi} - T_{ \frac{1+\eta_x^2}{\delta} }  T_{ \partial_w \partial_w^{\varrho} \hat{\varphi} } \eta + T_{ \partial_w\hat{\varphi} } \frac{ 1+\eta_x^2}{\delta}   - G(\eta,\beta,\gamma)(\psi) \in H^{2s-2}(\mathbb{T}) , 
\end{align*}
equivalently
\begin{align*}
& T_{ \frac{1+ \eta_x^2}{\delta} } \partial_w \hat{\Phi} - T_{\eta_x} \partial_x \hat{\Phi} - T_{ \partial_x \hat{\varphi} - \eta_x \partial_w^{\varrho} \hat{\varphi} } \eta_x - T_{ \partial_x ( \partial_x \hat{\varphi} - \eta_x \partial_w^{\varrho} \hat{\varphi} ) } \eta \\
&\qquad - G(\eta,\beta,\gamma)(\psi) \in H^{2s-2}(\mathbb{T}) , 
\end{align*}

Notice that by \eqref{eq:ParaGoodEq} we have
\begin{align*}
T_{ \frac{1+\eta_x^2}{\delta} } \partial_w\hat{\Phi} - T_{ \eta_x } \, \partial_x\hat{\Phi} &= T_{\lambda} \omega + \mathcal{R}_1,
\end{align*}
where $\lambda \in \Sigma^1_{s-1}(\mathbb{T})$ is given by
\begin{align} \label{eq:lambdaOp}
\lambda &= \frac{1+\eta_x^2}{\delta}  \, \mathtt{M} - \mathrm{i} \, \eta_x \, \xi ,
\end{align} 
and $\mathcal{R}_1 \in H^{s+1/2}(\mathbb{T})$. Therefore, we obtain
\begin{align*}
G(\eta,\beta,\gamma)(\psi) &= T_{\lambda} \omega - T_{ \partial_x \hat{\varphi} - \eta_x \partial_w^{\varrho} \hat{\varphi} } \eta_x - T_{ \partial_x ( \partial_x \hat{\varphi} - \eta_x \partial_w^{\varrho} \hat{\varphi} ) } \eta + \mathcal{R}_2 ,
\end{align*}
where $\mathcal{R}_2 \in H^{s+1/2}(\mathbb{T})$, and since
\begin{align*}
V(\eta,\beta,\gamma)(\psi) &= \partial_x \hat{\varphi} - \eta_x \partial_w^{\varrho} \hat{\varphi} |_{w=0} , \quad T_{ \partial_x V(\eta,\beta,\gamma)(\psi) } \eta \in H^{s+1/2}(\mathbb{T}) , 
\end{align*}
we can deduce the thesis.

\section{Local well-posedness} \label{sec:LWP}

In this section we prove the local well-posedness result Theorem \ref{thm:LWPMain} for the system \eqref{eq:WWsys}. Following the notation used in Theorem \ref{thm:LWPMain}, we recall that $s_0 = \frac{5}{2}$.

\subsection{Paralinearization of the water waves system} \label{subsec:paralin}

Let us consider the system \eqref{eq:WWsys}, namely
\begin{equation*} 
\begin{cases}
\eta_t &= G(\eta,\beta,\gamma)(\psi) + \gamma \, \eta \, \eta_x , \\
\psi_t &= - \frac{1}{2} \, \psi_x^2 + \frac{ 1 }{ 2(1+\eta_x^2)}   ( G(\eta,\beta,\gamma)(\psi) + \eta_x \psi_x )^2 + \gamma \, \left( \eta \, \psi_x +  \partial_x^{-1}G(\eta,\beta,\gamma)(\psi) \right) \\
&\quad - g \, \eta + \kappa \, \left(  \frac{\eta_x}{ \sqrt{1+\eta_x^2}} \right)_x   .
\end{cases}
\end{equation*}

In this subsection we consider a given solution $(\eta,\psi)$ of the system \eqref{eq:WWsys} on the interval $[0,T]$, $T < +\infty$, such that
\begin{equation} \label{eq:SolSys}
(\eta,\psi) \in C^0( [0,T]; H^{s+1/2}_0(\mathbb{T}) \times H^{s}(\mathbb{T}) ) , \;\; \text{for some} \;\; s > s_0.
\end{equation}

We want to paralinearize system \eqref{eq:WWsys}. In Proposition \ref{prop:OpGpara} we performed the paralinearization of the operator $G(\eta,\beta,\gamma)$. Now we deal with the other terms in \eqref{eq:WWsys}.

\begin{lemma} \label{lem:ParaTech1}

Let $s > s_0$, and assume that $(\eta,\psi)$ satisfies \eqref{eq:SolSys} and solves the system \eqref{eq:WWsys}. Then 
\begin{align*}
\left(  \frac{\eta_x}{ \sqrt{1+\eta_x^2}} \right)_x &= - T_{\mathtt{h}} \eta + g_1 ,
\end{align*}
where
\begin{align}
\mathtt{h} &= \mathtt{h}^{(2)} + \mathtt{h}^{(1)}, \;\; \mathtt{h}^{(2)} = \frac{ \xi^2 }{ (1+\eta_x^2)^{3/2} }  , \;\; \mathtt{h}^{(1)} = - \frac{\mathrm{i}}{2} \, (\partial_x \, \partial_{\xi}) \mathtt{h}^{(2)} , \label{eq:hSymbol} 
\end{align}
and where $g_1 \in L^{\infty}([0,T]; H^{2s-5/2}(\mathbb{T}) )$ satisfies
\begin{align*}
\| g_1 \|_{ L^{\infty}([0,T]; H^{2s-5/2}(\mathbb{T}) ) } &\leq C \left( \| \eta \|_{ L^{\infty}([0,T]; H^{s+1/2}_0(\mathbb{T}) ) } \right),
\end{align*}
for some non-decreasing function $C$.

\end{lemma}

For the proof of Lemma \ref{lem:ParaTech1} see the proof of Lemma 3.25 in \cite{alazard2011water}.

Now, let us recall the operators $B(\eta,\beta,\gamma)$ and $V(\eta,\beta,\gamma)$ given by
\begin{align*}
B(\eta,\beta,\gamma)(\psi) &= \frac{ G(\eta,\beta,\gamma)(\psi) + \eta_x \psi_x }{1+\eta_x^2} , \quad V(\eta,\beta,\gamma)(\psi) = \psi_x - B(\eta,\beta,\gamma)(\psi) \, \eta_x .
\end{align*}

\begin{lemma} \label{lem:ParaTech2}

Let $s > s_0$. Assume that $(\eta,\psi)$ satisfies \eqref{eq:SolSys} and solves the system \eqref{eq:WWsys}, and that $\eta$ and $\beta$ satisfy \eqref{eq:StrConnected} uniformly on $[0,T]$. Then 
\begin{align*}
& \frac{1}{2} \, \psi_x^2 - \frac{ 1 }{ 2(1+\eta_x^2)}   ( G(\eta,\beta,\gamma)(\psi) + \eta_x \psi_x )^2 - \gamma \, \left( \eta \, \psi_x +  \partial_x^{-1}G(\eta,\beta,\gamma)(\psi) \right) \\
&= T_{V(\eta,\beta,\gamma)(\psi) - \gamma \, \eta} \, \psi_x - T_{ V(\eta,\beta,\gamma)(\psi) } T_{ B(\eta,\beta,\gamma)(\psi) } \eta_x  - T_{ B(\eta,\beta,\gamma)(\psi) } G(\eta,\beta,\gamma)(\psi) + g_2 ,
\end{align*}
where $g_2 \in L^{\infty}([0,T]; H^{s}(\mathbb{T}) )$ satisfies
\begin{align*}
\| g_2 \|_{ L^{\infty}([0,T]; H^{s}(\mathbb{T}) ) } &\leq  C \left( \| (\eta,\psi) \|_{ L^{\infty}( [0,T]; H^{s+1/2}_0(\mathbb{T}) \times H^{s}(\mathbb{T}) )  }  \right) ,
\end{align*}
for some non-decreasing function $C$ depending only on $\eta_0$, $h$, $h_0$, $|\gamma|$ and $\|\beta\|_{ H^{s+1/2}(\mathbb{T}) }$. 

\end{lemma}

\begin{proof}

First, we notice that 
\begin{align*}
F(a,b,c) &= \frac{1}{2} \frac{ (a b +c)^2 }{1+a^2} , \;\; a,b,c \in \mathbb{R}, \\
\partial_a F &= \frac{a b +c}{1+a^2} \left[ b - \frac{a b +c}{1+a^2} a \right] , \;\; \partial_b F = \frac{ a b +c }{1+a^2} a, \;\; \partial_c F = \frac{ a b +c }{1+a^2} .
\end{align*}
Using the above relations with $a=\eta_x$, $b=\psi_x$ and $c = G(\eta,\beta,\gamma)(\psi)$, by Lemma \ref{lem:Paraprod2} we obtain
\begin{align*}
& \frac{1}{2} \, \psi_x^2 - \frac{ 1 }{ 2(1+\eta_x^2)}   ( G(\eta,\beta,\gamma)(\psi) + \eta_x \psi_x )^2 - \gamma \, \left( \eta \, \psi_x +  \partial_x^{-1}G(\eta,\beta,\gamma)(\psi) \right) \\
&= T_{V(\eta,\beta,\gamma)(\psi) - \gamma \, \eta} \, \psi_x - T_{ V(\eta,\beta,\gamma)(\psi) \, B(\eta,\beta,\gamma)(\psi) } \eta_x  - T_{ B(\eta,\beta,\gamma)(\psi) } G(\eta,\beta,\gamma)(\psi) + g_{2,0} ,
\end{align*}
where $g_{2,0} \in L^{\infty}([0,T]; H^{s}(\mathbb{T}) )$ satisfies
\begin{align*}
\| g_{2,0} \|_{ L^{\infty}([0,T]; H^{s}(\mathbb{T}) ) } &\leq  C \left( \| (\eta,\psi) \|_{ L^{\infty}( [0,T]; H^{s+1/2}_0(\mathbb{T}) \times H^{s}(\mathbb{T}) )  }  \right) ,
\end{align*}
for some non-decreasing function $C$ depending only on $\eta_0$, $h$, $h_0$ and $|\gamma|$. 

Since $T_{ B(\eta,\beta,\gamma)(\psi) V(\eta,\beta,\gamma)(\psi) } - T_{ B(\eta,\beta,\gamma)(\psi)  } T_{  V(\eta,\beta,\gamma)(\psi) }$ is of order $- (s-3/2)$, we can deduce the thesis.
\end{proof}

\begin{lemma} \label{lem:ParaTech3}

Let $s > s_0$. Assume that $(\eta,\psi)$ satisfies \eqref{eq:SolSys} and solves the system \eqref{eq:WWsys}, and that $\eta$ and $\beta$ satisfy \eqref{eq:StrConnected} uniformly on $[0,T]$. Then 
\begin{align*}
\| T_{ \partial_t B(\eta,\beta,\gamma)(\psi) } \eta \|_{ H^{s}(\mathbb{T}) ) } &\leq  C \left( \| (\eta,\psi) \|_{ L^{\infty}( [0,T]; H^{s+1/2}_0(\mathbb{T}) \times H^{s}(\mathbb{T}) )  }  \right) ,
\end{align*}
for some non-decreasing function $C$ depending only on $\eta_0$, $h$, $h_0$, $|\gamma|$ and $\|\beta\|_{ H^{s+1/2}(\mathbb{T}) }$. 
\end{lemma}

\begin{proof}

By \eqref{eq:EstGUSob} and by exploiting the structure of \eqref{eq:WWsys}, we obtain
\begin{align*}
\| \eta_t \|_{ H^{s-1}(\mathbb{T}) } + \| \psi_t \|_{ H^{s-3/2}(\mathbb{T}) } &\leq  C \left( \| (\eta,\psi) \|_{ L^{\infty}( [0,T]; H^{s+1/2}_0(\mathbb{T}) \times H^{s}(\mathbb{T}) )  }  \right) .
\end{align*}
Using Proposition \ref{prop:SysOpDer} we have
\begin{align*}
\partial_t \left[ G(\eta,\beta,\gamma)(\psi) \right] &= G^{DN}\left( \psi_t - B(\eta,\beta,\gamma)(\psi) \eta_t \right) - \partial_x ( V(\eta,\beta,\gamma)(\psi) \eta_t ) \\
&\;\; + \gamma \, G^{NN}(\eta,\beta)( (-h+\beta) \beta_x ),
\end{align*}
therefore
\begin{align*}
\| \partial_t \left[ G(\eta,\beta,\gamma)(\psi) \right] \|_{ H^{s-5/2}(\mathbb{T}) ) } &\leq  C \left( \| (\eta,\psi) \|_{ L^{\infty}( [0,T]; H^{s+1/2}_0(\mathbb{T}) \times H^{s}(\mathbb{T}) )  }  \right) ,
\end{align*}
and we can deduce the thesis.
\end{proof}

Combining Proposition \ref{prop:OpGpara}, Lemma \ref{lem:ParaTech1}, Lemma \ref{lem:ParaTech2} and Lemma \ref{lem:ParaTech3} we obtain

\begin{proposition} \label{prop:ParalinFull}

Let $s > s_0$. Assume that $(\eta,\psi)$ satisfies \eqref{eq:SolSys} and solves the system \eqref{eq:WWsys}, and that $\eta$ and $\beta$ satisfy \eqref{eq:StrConnected} uniformly on $[0,T]$. Then, recalling the good unknown $\omega$ defined in \eqref{eq:defGU}, we have that $(\eta,\omega)$ solves the system
\begin{equation} \label{eq:ParalinFull}
\begin{cases}
\eta_t + T_{V(\eta,\beta,\gamma)(\psi) - \gamma \eta} \, \eta_x - T_{\lambda}\omega &= f , \\
\omega_t + T_{V(\eta,\beta,\gamma)(\psi) - \gamma \eta} \; \omega_x + \kappa \, T_{\mathtt{h}} \eta &= g ,
\end{cases}
\end{equation}
where
\begin{align*}
f \in L^{\infty}([0,T]; H^{s+1/2}(\mathbb{T}) ) , &\;\; g \in L^{\infty}([0,T]; H^{s}(\mathbb{T}) ) , \\
\| (f,g) \|_{ L^{\infty}([0,T]; H^{s+1/2}(\mathbb{T}) \times H^{s}(\mathbb{T}) ) } &\leq  C \left( \| (\eta,\psi) \|_{ L^{\infty}( [0,T]; H^{s+1/2}_0(\mathbb{T}) \times H^{s}(\mathbb{T}) )  }  \right) ,
\end{align*}
for some non-decreasing function $C$ depending on $\eta_0$, $h$, $h_0$, $|\gamma|$, $\kappa$ and $\|\beta\|_{ H^{s+1/2}(\mathbb{T}) }$.

\end{proposition}

\subsection{Symmetrization} \label{subsec:symm}

Let $s > s_0$. Assume that $(\eta,\psi)$ satisfies \eqref{eq:SolSys} and solves the system \eqref{eq:WWsys}, and that $\eta$ and $\beta$ satisfy \eqref{eq:StrConnected} uniformly on $[0,T]$. In this subsection we consider the system \eqref{eq:ParalinFull}, namely
\begin{equation*}
(\partial_t + T_{V(\eta,\beta,\gamma)(\psi) - \gamma \eta} \, \partial_x) \,
\begin{pmatrix}
\eta \\ \omega
\end{pmatrix}
+ 
\begin{pmatrix}
0 & - T_{\lambda} \\
\kappa T_{\mathtt{h}} & 0
\end{pmatrix}
\begin{pmatrix}
\eta \\ \omega
\end{pmatrix}
=
\begin{pmatrix}
f \\ g
\end{pmatrix}
,
\end{equation*}
where $f \in L^{\infty}([0,T]; H^{s+1/2}(\mathbb{T}) )$ and $g \in L^{\infty}([0,T]; H^{s}(\mathbb{T}) )$, and we use paradifferential calculus in order to construct a symmetrizer 
\begin{align*}
S &=
\begin{pmatrix}
T_p & 0 \\
0 & T_q
\end{pmatrix}
,
\end{align*}
such that (up to a smoother remainder)
\begin{align*}
S 
\begin{pmatrix}
0 & - T_{\lambda} \\
\kappa T_{\mathtt{h}} & 0
\end{pmatrix}
&= 
\begin{pmatrix}
0 & -T_{\vartheta} \\
(T_{\vartheta})^{\ast} & 0 
\end{pmatrix}
S ,
\end{align*}
and such that the unknown
\begin{align} \label{eq:SymmUn}
\Phi &\coloneqq S
\begin{pmatrix}
\eta \\ \omega
\end{pmatrix}
\end{align}
satisfies
\begin{align} \label{eq:SymmSys}
\Phi_t + T_{V(\eta,\beta,\gamma)(\psi) - \gamma \, \eta} \Phi_x + 
\begin{pmatrix}
0 & -T_{\vartheta} \\
(T_{\vartheta})^{\ast} & 0 
\end{pmatrix}
\Phi &= F ,
\end{align}
with $F \in L^{\infty}([0,T] ; H^s(\mathbb{T}) \times H^s(\mathbb{T}) )$. Moreover, we have
\begin{align*}
\| F \|_{ L^{\infty}([0,T] ; H^s(\mathbb{T}) \times H^s(\mathbb{T}) ) } &\leq C \left( \| (\eta,\psi) \|_{ L^{\infty}( [0,T]; H^{s+1/2}_0(\mathbb{T}) \times H^{s}(\mathbb{T}) )  }  \right) ,
\end{align*}
for some non-decreasing function $C$ depending on $\eta_0$, $h$, $h_0$, $|\gamma|$, $\kappa$ and $\|\beta\|_{ H^{s+1/2}(\mathbb{T}) }$.

First, we introduce the main definitions and results for symbolic calculus with low regularity. In the following we write $\mathbb{R}_{\neq} := \mathbb{R}\setminus\{0\}$.

All the symbols which we consider below are of the form
\begin{align*}
a &= a^{(m)}+a^{(m-1)} ,
\end{align*} 
where
\begin{itemize}
\item[i.] $a^{(m)}$ is a real-valued  elliptic symbol, homogenous of degree $m$ in $\xi$ and depends only on the first order-derivatives of $\eta$;
\item[ii.] $a^{(m-1)}$ is homogenous of degree $m-1$ in $\xi$ and depends also, but only linearly, on the second order-derivatives of $\eta$.
\end{itemize}
We recall that $\eta\in C^{0}([0,T];H^{s+\frac{1}{2}}(\mathbb{T}))$ is a fixed given function.

\begin{definition} \label{def:defiSigma}
Given $m\in \mathbb{R}$, $\Sigma^m$ denotes the class of symbols $a$ of the form
$$
a=a^{(m)}+a^{(m-1)}
$$
with
$$
a^{(m)}(t,x,\xi)= F(\eta_x(t,x),\xi), \quad 
a^{(m-1)}(t,x,\xi)= G_2(\eta_x(t,x),\xi)\partial_x^2 \eta (t,x),
$$
such that
\begin{itemize}
\item[i] $T_a$ maps real-valued functions to real valued functions;
\item[ii] $F$ is a $C^\infty$ real-valued function of $(\zeta,\xi)\in \mathbb{R}\times \mathbb{R}_{\neq}$, 
homogeneous of order $m$  in $\xi$, and such that there exists a continuous function $K=K(\zeta)>0$ such that 
$$
F(\zeta,\xi)\geq K (\zeta) |\xi|^m,
$$
for all $(\zeta,\xi)\in \mathbb{R} \times  \mathbb{R}_{\neq}$;
\item[iii] $G_2$ is a $C^\infty$ complex-valued function of $(\zeta,\xi)\in \mathbb{R}\times \mathbb{R}_{\neq}$, 
homogeneous of order $m-1$ in $\xi$. 
\end{itemize}
\end{definition}

Notice that, since we only assume $s > s_0$, we have some additional technical difficulties compared for example to the case $s > \frac{7}{2}$. In order to overcome this problem, we exploit that for all our symbols the sub-principal terms have only a linear dependence on the second order derivative of $\eta$.

Next, we state a result ensuring that the previous class of symbols is stable by the standard rules of symbolic calculus, and a symbolic calculus result modulo admissible remainders; before stating the results, we introduce the following notation.

\begin{definition} \label{def:defisym}
Let $m\in\mathbb{R}$ and consider two families of operators order $m$, 
\begin{equation*}
\{ A_1(t) \,:\, t\in [0,T]\},\quad 
\{ A_2(t) \,:\, t\in [0,T]\}.
\end{equation*}
We say that $A_1 \sim A_2$ if $A_1-A_2$ is of order $m-3/2$ (see Definition \ref{def:OrdOp}) and satisfies the following estimate: for all $\mu\in\mathbb{R}$, 
there exists a continuous function $C$ such that
\begin{equation*}
\| A_1(t)-A_2(t)\|_{H^{\mu}(\mathbb{T}) \rightarrow H^{\mu-m+\frac{3}{2}}(\mathbb{T}) }\le C\left( \| \eta(t) \|_{H^{s+\frac{1}{2}}(\mathbb{T}) }\right), \;\; \forall t \in [0,T].
\end{equation*}
\end{definition}

First, we have the following result (see Proposition 4.3 in \cite{alazard2011water}).

\begin{proposition} \label{prop:pcs}

Let $m,m'\in\mathbb{R}$. Then
\begin{itemize}
\item[i] If $a\in \Sigma^m$ and $b\in \Sigma^{m'}$ then 
$T_a T_b\sim T_{a\sharp b}$ where $a\sharp b\in \Sigma^{m+m'}$ is given by 
\begin{equation*}
a\sharp b =a^{(m)}b^{(m')}+a^{(m-1)}b^{(m')}+a^{(m)}b^{(m'-1)}+\frac{1}{i}
\partial_\xi a^{(m)}\cdot \partial_x b^{(m')}.
\end{equation*}
\item[ii] If $a\in \Sigma^m$ then $(T_a)^{\ast} \sim T_b$ where $b\in \Sigma^m$ is given by
\begin{equation*}
b=a^{(m)}+\overline{a^{(m-1)}} + \frac{1}{i}(\partial_x\cdot \partial_\xi) a^{(m)}.
\end{equation*}
\end{itemize}

\end{proposition}

Given that $a\in \Sigma^{m}$, since $a^{(m-1)}$ involves two derivatives of 
$\eta$, Proposition \ref{prop:AsympParadiffOp} and the embedding 
$H^{s}(\mathbb{T})\subset W^{2,\infty}(\mathbb{T})$ implies that
\begin{equation} \label{eq:ong}
\| T_{a(t)} \|_{H^{\mu}(\mathbb{T}) \rightarrow H^{\mu-m}(\mathbb{T})} \leq C \, \sup_{| \xi |=1}\| a(t,\cdot,\xi)\|_{L^\infty} \leq C\left( \| \eta(t)\|_{H^s(\mathbb{T})}\right). 
\end{equation}
Our second observation concerning the class $\Sigma^m$ is that one can prove a continuity result which requires 
only an estimate of $\| \eta \|_{ H^{s-1}(\mathbb{T}) }$.

\begin{proposition} \label{prop:2d21}
Let $m\in\mathbb{R}$ and $\mu\in \mathbb{R}$. Then 
there exists a function $C$ such that for all symbol $a\in \Sigma^m$ and all $t\in [0,T]$,
\begin{align*}
\| T_{a(t)} u \|_{H^{\mu-m}(\mathbb{T}) } \le C( \| \eta(t) \|_{ H^{s-1}(\mathbb{T} }) \| u \|_{ H^{\mu}(\mathbb{T}) }.
\end{align*}
\end{proposition}

For the proof of the above proposition, see Proposition 4.4 in \cite{alazard2011water}. We point out that the above result follows by Proposition \ref{prop:OrdParadiffOp} for $s>\frac{7}{2}$, since the $L^\infty$-norm of $a(t,\cdot,\xi)$ is controlled by $\| \eta(t) \|_{ H^{s-1}(\mathbb{T}) }$ in this case; the above proposition is necessary to overcome this issue for $s > s_0$.

Similarly, we have the following elliptic regularity result, where one controls the constants in terms of the $H^{s-1}$-norm of $\eta$ only (see Proposition 4.6 in \cite{alazard2011water}).

\begin{proposition}\label{prop:2d22}

Let $m\in\mathbb{R}$ and $\mu\in \mathbb{R}$. Then 
there exists a function $C$ such that for all 
$a\in \Sigma^m$ 
and all $t\in [0,T]$, we have
\begin{align*}
\| u \|_{ H^{\mu+m}(\mathbb{T}) } &\leq  C(\| \eta(t) \|_{ H^{s-1}(\mathbb{T}) }) \left[ \| T_{a(t)} u \|_{H^\mu(\mathbb{T})} +  \| u \|_{L^2(\mathbb{T})} \right] .
\end{align*}

\end{proposition}

The classical result is that, for all elliptic symbol $a\in \Gamma^{m}_\rho(\mathbb{T})$ with $\rho>0$, the following estimate holds true,
\begin{align*}
\| f \|_{H^m(\mathbb{T})} &\leq C(\| \eta(t) \|_{ W^{2+\rho,\infty}(\mathbb{T}) } ) \left[ \| T_a f \|_{L^2(\mathbb{T}} + \| f \|_{L^2(\mathbb{T}} \right] ,
\end{align*}
for $\rho>0$ small enough, which would lead to a worse estimate for $\frac{5}{2} < s <\frac{7}{2}$.

Arguing as in Sec. 4.2 of \cite{alazard2011water}, one can deduce the following result.

\begin{proposition} \label{prop:key}
There exist $q\in\Sigma^0$, $p\in \Sigma^{1/2}$, $\vartheta \in \Sigma^{3/2}$ such that
\begin{equation} \label{eq:sy}
\begin{cases}
T_p T_\lambda &\sim T_\vartheta T_q ,\\
T_{q} \, \kappa T_{\mathtt{h}} &\sim T_{\vartheta} T_{p} ,\\
T_\vartheta &\sim (T_{\vartheta})^*,
\end{cases}
\end{equation}
where the notation $A_1 \sim A_2$ is the one introduced in Definition \ref{def:defisym}.

\end{proposition}

\begin{remark} \label{rem:key}

By following the argument of Sec. 4.2 of \cite{alazard2011water}, we obtain that
\begin{align*}
q &= \left( 1+\eta_x^2 \right)^{-1/2},\\
p &= \left( 1+\eta_x^2 \right)^{-\frac{5}{4}}\sqrt{\lambda^{(1)}}+p^{(-1/2)},\\
\vartheta &= \sqrt{ \kappa \mathtt{h}^{(2)}\lambda^{(1)}}+\sqrt{ \frac{\kappa \mathtt{h}^{(2)}}{\lambda^{(1)}} } \frac{ \mathrm{Re} \, \lambda^{(0)}}{2}
-\frac{\mathrm{i}}{2} \partial_\xi \partial_x ( \sqrt{\kappa \mathtt{h}^{(2)}\lambda^{(1)}} ), 
\end{align*}
where $\lambda^{(1)}$ is the principal symbol of the classical Dirichlet--Neumann operator, $\lambda^{(0)}$ is the sub-principal of the symbol $\lambda$ given by \eqref{eq:lambdaOp}, $\mathtt{h}^{(2)}$ is given by \eqref{eq:hSymbol}, and $p^{(-1/2)}$ is a symbol homogeneous in $\xi$ of order $-\frac{1}{2}$.

\end{remark}

We now prove a technical lemma.

\begin{lemma} \label{lem:comdtp}
For all $\mu \in \mathbb{R}$ there exists a non-decreasing function $C$ such that, for all $t\in [0,T]$,
\begin{align*}
\| T_{ \partial_{t}p(t)} \|_{H^{\mu}(\mathbb{T}) \rightarrow H^{\mu-1/2}(\mathbb{T}) } +
\| T_{\partial_{t}q(t)} \|_{H^{\mu}(\mathbb{T}) \rightarrow H^{\mu}(\mathbb{T}) } &\leq C\left( \| (\eta(t),\psi(t)) \|_{H^{s+1/2}_0(\mathbb{T}) \times H^s(\mathbb{T})} \right).
\end{align*}
\end{lemma}

\begin{proof}

First, observe that
\begin{align*}
\| \eta_t \|_{W^{1,\infty}(\mathbb{T})} \leq C_1 \, \| \eta_t \|_{H^{s-1}(\mathbb{T})} &\leq
C_2 \left( \|  (\eta,\psi) \|_{H^{s+1/2}_0(\mathbb{T}) \times H^s(\mathbb{T})} \right).
\end{align*}
This implies that
\begin{align*}
\| q_t(t,\cdot) \|_{ L^{\infty}(\mathbb{T}) } + \sup_{|\xi|=1} \| p^{(1/2)}_t(t,\cdot,\xi) \|_{L^{\infty}(\mathbb{T}) } &\leq C\left( \| ( \eta ,\psi ) \|_{H^{s+1/2}_0(\mathbb{T}) \times H^s(\mathbb{T}) }\right).
\end{align*}

By Proposition \ref{prop:OrdParadiffOp}, this implies that 
\begin{align*}
\| T_{ p^{(1/2)}_t(t)} \|_{H^{\mu}(\mathbb{T}) \rightarrow H^{\mu-1/2}(\mathbb{T}) } + 
\| T_{q_t(t)} \|_{H^{\mu}(\mathbb{T}) \rightarrow H^{\mu}(\mathbb{T})} &\leq C\left( \| (\eta(t),\psi(t))\|_{H^{s+1/2}_0(\mathbb{T}) \times H^s(\mathbb{T})}\right).
\end{align*}

We now need to estimate $\| T_{ p^{(-1/2)}_t(t)} \|_{H^{\mu}(\mathbb{T}) \rightarrow H^{\mu-1/2}(\mathbb{T}) }$: we rewrite $p^{(-1/2)}$ as
\begin{align*}
p^{(-1/2)} &= P_{2}(\eta_x,\xi)\eta_{xx},
\end{align*}
where $P_{2}$ is a smooth functions for $\xi\neq 0$, homogeneous of degree $- \frac{1}{2}$ in $\xi$. We write
\begin{equation} \label{eq:ComSdt}
T_{ p^{(-1/2)}_t} 
= T_{(P_{2,t}(\eta_x,\xi) ) \eta_{xx}} + T_{P_2( \eta_x,\xi)\eta_{xxt} } ,
\end{equation}
hence
\begin{align*}
\sup_{| \xi |=1} \| P_{2,t}(\eta_x(t),\xi) \|_{L^{\infty}(\mathbb{T}) } &\leq C\left( \| ( \eta,\psi) \|_{H^{s+1/2}_0(\mathbb{T}) \times H^s(\mathbb{T}) }\right) ,
\end{align*}
and since $s>s_0$, then the first term in the right hand side of \eqref{eq:ComSdt} is uniformly of order $-\frac{1}{2}$. To estimate the second term in the right-hand side of \eqref{eq:ComSdt}, notice that
\begin{align*}
\| T_{P_{2}( \eta_x,\xi)\eta_{xxt} } \|_{H^\mu(\mathbb{T}) \rightarrow H^{\mu-1/2}(\mathbb{T}) }
&\leq \| P_{2}( \eta_x,\xi)\eta_{xxt} \|_{H^{s-3}(\mathbb{T}) } ,
\end{align*}
and by standard product rule in Sobolev spaces we have
\begin{align*}
&\| P_{2}( \eta_x,\xi)\eta_{xxt} \|_{H^{s-3}(\mathbb{T}) }\\
&\leq C \, \left[ | P_2(0,\xi) | + | P_{2}( \eta_x,\xi)-P_2(0,\xi) |_{H^{s-1}(\mathbb{T}) } \right] \, \| \eta_{xxt} \|_{H^{s-3}(\mathbb{T}) },
\end{align*}
and thus
\begin{align*}
\| T_{P_{2}( \eta_x,\xi) \eta_{xxt} } \|_{H^\mu(\mathbb{T}) \rightarrow H^{\mu-m}(\mathbb{T}) }
&\leq C\left( \| ( \eta,\psi) \|_{H^{s+1/2}_0(\mathbb{T}) \times H^s(\mathbb{T}) } \right) ,
\end{align*}
which leads to the thesis.
\end{proof}

By combining Proposition \ref{prop:key} with the paralinearization obtained in Proposition \ref{prop:ParalinFull}, we obtain the following symmetrization of the system \eqref{eq:WWsys}.

\begin{corollary}\label{cor:psym}
Introduce the new unknowns
\begin{align*}
\Phi_{1} \coloneqq T_{p} \eta , &\quad \Phi_{2} \coloneqq T_{q} \omega.
\end{align*}
Then $\Phi_1,\Phi_2 \in C^{0}([0,T];H^s(\mathbb{T}))$ and
\begin{equation} \label{eq:systred}
\begin{cases}
\Phi_{1,t}+T_{V(\eta,\beta,\gamma)(\psi) - \gamma \, \eta} \Phi_{1,x} - T_{\vartheta} \Phi_{2} &= F_1, \\
\Phi_{2,t}+T_{V(\eta,\beta,\gamma)(\psi) - \gamma \, \eta} \Phi_{2,x} + T_{\vartheta} \Phi_{1} &= F_2,
\end{cases}
\end{equation}
where $F_{1},F_{2}\in L^{\infty}([0,T];H^{s}(\mathbb{T}))$. 
Moreover 
\begin{align*}
\| (F_1,F_2) \|_{L^{\infty}([0,T];H^s(\mathbb{T}) \times H^s(\mathbb{T}) )} &\leq C \left( \| (\eta,\psi) \|_{L^\infty([0,T];H^{s+1/2}_0(\mathbb{T}) \times H^s(\mathbb{T}) )} \right),
\end{align*}
for some function $C$ depending only on $\eta_0$, $h$, $h_0$, $|\gamma|$ and $\|\beta\|_{ H^{s+1/2}(\mathbb{T}) }$.
\end{corollary}

\begin{proof}

We first observe that by Proposition \ref{prop:key} and Proposition \ref{prop:ParalinFull} that
\begin{equation*}
\begin{cases}
\Phi_{1,t} + T_{V(\eta,\beta,\gamma)(\psi) - \gamma \, \eta }\Phi_{1,x} - T_{\vartheta} \Phi_2 &= B_{1}\eta + f_1,
\\
\Phi_{2,t} + T_{V(\eta,\beta,\gamma)(\psi) - \gamma \, \eta }\Phi_{2,x} + T_{\vartheta} \Phi_1 &= B_{2}\omega + f_2,
\end{cases}
\end{equation*}
with $f_{1},f_{2}\in L^{\infty}([0,T];H^{s}(\mathbb{T}))$ satisfying
\begin{align*}
\| (f_1,f_2) \|_{L^{\infty}([0,T];H^{s}(\mathbb{T}) )} &\leq C \left( \| (\eta,\psi) \|_{L^\infty([0,T];H^{s+1/2}_0(\mathbb{T})\times H^{s}(\mathbb{T}))}\right),
\end{align*}
and 
\begin{align*}
B_{1} &\coloneqq [ \partial_{t}, T_p]+\left[ T_{V(\eta,\beta,\gamma)(\psi) - \gamma \, \eta} \,  \partial_x,T_{p}\right] , \;\; B_{2} \coloneqq  [ \partial_{t}, T_q] +\left[ T_{V(\eta,\beta,\gamma)(\psi) - \gamma \, \eta} \, \partial_x,T_{q}\right] .
\end{align*}
Since
\begin{align*}
\| B_{1}\eta \|_{H^s(\mathbb{T})} &\leq \| B_1 \|_{H^{s+1/2}(\mathbb{T}) \rightarrow H^s(\mathbb{T}) } \| \eta \|_{H^{s+1/2}_0(\mathbb{T}) },\\
\| B_{2} \omega \|_{H^s(\mathbb{T})} &\leq \| B_2 \|_{H^{s}(\mathbb{T}) \rightarrow H^{s}(\mathbb{T})} \| \omega \|_{H^{s}(\mathbb{T})},
\end{align*}
and since we can estimate $\| B_1 \|_{H^{s+1/2}(\mathbb{T}) \rightarrow H^s(\mathbb{T}) }$ and $\| B_2 \|_{H^{s}(\mathbb{T}) \rightarrow H^s(\mathbb{T}) }$ by Lemma \ref{lem:comdtp}, we can deduce the thesis. 
\end{proof}

In the irrotational case the symmetrization \eqref{eq:systred} of the water waves system has been useful in order to prove a controllability result, see Sec. 2 of \cite{alazard2018control}.

\subsection{A priori estimates} \label{subsec:apriori}

Let $s > s_0$, and assume that $\eta \in C([0,T]; H^{s+1/2}_0(\mathbb{T}) )$ and $\beta \in H^{s+1/2}(\mathbb{T})$ satisfy \eqref{eq:StrConnected} uniformly on $[0,T]$. Now we prove a priori estimates for the system \eqref{eq:WWsys} with initial data \eqref{eq:InData}. These estimates are based on the following reformulation, which follows by direct computation.

\begin{lemma} \label{lem:formulation}

$(\eta,\psi)$ solves \eqref{eq:WWsys} if and only if 
\begin{align*}
&
\begin{pmatrix}
\mathrm{Id} & 0 \\
-T_{B(\eta,\beta,\gamma)(\psi)} & \mathrm{Id} 
\end{pmatrix}
\left( \partial_t+ T_{V(\eta,\beta,\gamma)(\psi) - \gamma \, \eta} \partial_x \right) 
\begin{pmatrix} 
\eta \\ \psi 
\end{pmatrix}
\\
&\quad +
\begin{pmatrix}
0 & -T_{\lambda} \\
\kappa T_{\mathtt{h}} & 0 \end{pmatrix}
\begin{pmatrix}
I & 0 \\
-T_{B(\eta,\beta,\gamma)(\psi)} & \mathrm{Id} 
\end{pmatrix}
\begin{pmatrix} 
\eta \\ \psi 
\end{pmatrix}
=
\begin{pmatrix} 
f_1 \\ f_2 
\end{pmatrix}
,
\end{align*}
where
\begin{align} 
f_1 &= G(\eta,\beta,\gamma)(\psi) - \left[ T_{\lambda} ( \psi -T_{B(\eta,\beta,\gamma)(\psi)} \eta ) - T_{V(\eta,\beta,\gamma)(\psi)} \eta_x \right] +\gamma \eta \eta_x -T_{\gamma \, \eta} \eta_x , \nonumber \\
f_2 &= -\frac{1}{2} \psi_x^2 + \frac{\left( G(\eta,\beta,\gamma)(\psi) + \eta_x \psi_x \right)^2}{2(1+\eta_x^2)} \nonumber \\
&\quad +\gamma \left[\eta \psi_x+\partial_x^{-1}G(\eta,\beta,\gamma)(\psi) \right] -g \, \eta + \kappa \left( \frac{\eta_x}{\sqrt{1+\eta_x^2}} \right)_x \nonumber \\
&\quad + T_{V(\eta,\beta,\gamma)(\psi) - \gamma \, \eta} \psi_x - T_{B(\eta,\beta,\gamma)(\psi)} T_{V(\eta,\beta,\gamma)(\psi) - \gamma \, \eta} \eta_x  \nonumber \\
&\quad - T_{B(\eta,\beta,\gamma)(\psi)} \left[ G(\eta,\beta,\gamma)(\psi) +\gamma \eta \eta_x \right] + \kappa T_{\mathtt{h}}\eta . \label{eq:f1f2}
\end{align}

\end{lemma}

Using Lemma \ref{lem:formulation} we can deduce that $(\eta,\psi)$ solves \eqref{eq:WWsys}  if and only if
\begin{align} \label{eq:ref}
\left( \partial_t+ T_{V(\eta,\beta,\gamma)(\psi)-\gamma \, \eta } \partial_x + \mathcal{L} \right) 
\begin{pmatrix}
\eta \\ \psi
\end{pmatrix}
&= f(\eta,\psi) ,
\end{align}
where
\begin{align*}
\mathcal{L} &\coloneqq
\begin{pmatrix}
\mathrm{Id} & 0 \\
T_{B(\eta,\beta,\gamma)(\psi)} & \mathrm{Id} 
\end{pmatrix}
\begin{pmatrix}
0 & -T_{\lambda}\\
\kappa T_{\mathtt{h}} & 0 
\end{pmatrix}
\begin{pmatrix}
\mathrm{Id} & 0 \\
-T_{B(\eta,\beta,\gamma)(\psi)} & \mathrm{Id} 
\end{pmatrix}
, \\
f(\eta,\psi) &\coloneqq
\begin{pmatrix}
\mathrm{Id} & 0 \\
T_{B(\eta,\beta,\gamma)(\psi)} & \mathrm{Id}
\end{pmatrix}
\begin{pmatrix} 
f_1 \\ f_2 
\end{pmatrix}
.
\end{align*}

We seek solutions of the Cauchy problem \eqref{eq:ref}, \eqref{eq:InData} as limits of solutions of approximating systems. The definition depends on two operators: the first one is a suitable mollifier, whereas the second one is a right-parametrix for the symmetrizer 
\begin{equation*}
S =
\begin{pmatrix}
T_p & 0 \\ 0 & T_q 
\end{pmatrix}
.
\end{equation*}

In order to regularize the equations, we define the following mollifiers: given $\epsilon \in [0,1]$, we define $J_{\epsilon}$ as the paradifferential operator with symbol $j_{\epsilon}=j_{\epsilon}(t,x,\xi)$ given by
\begin{align*}
j_{\epsilon} = j_{\epsilon}^{(0)} + j_{\epsilon}^{(-1)} = \exp \big( -\epsilon \gamma^{(3/2)}\big) -\frac{\mathrm{i}}{2} (\partial_x \partial_\xi)\exp \big( -\epsilon \gamma^{(3/2)}\big).
\end{align*}
Observe that
\begin{align*}
j_{\epsilon} \in C^{0}([0,T];\Gamma^0_{3/2}(\mathbb{T})), \quad \{ j_{\epsilon}^{(0)}, \gamma^{(3/2)} \}=0, \quad
\mathrm{Im} \, j_{\epsilon}^{(-1)} = -\frac{1}{2} \, \partial_x\partial_\xi j_{\epsilon}^{(0)}.
\end{align*}
For any $\epsilon >0$, $j_{\epsilon} \in C^{0}([0,T];\Gamma^m_{3/2}(\mathbb{T}))$ for all $m \leq 0$; moreover, $j_{\epsilon}$ is uniformly bounded in $C^{0}([0,T];\Gamma^0_{3/2}(\mathbb{T}))$ for all $\epsilon \in [0,1]$.  Therefore, we have
\begin{align*}
\| J_{\epsilon} T_{\vartheta} - T_{\vartheta} J_{\epsilon} \|_{H^\mu(\mathbb{T}) \rightarrow H^\mu(\mathbb{T}) } &\le\ C(\| \eta_x \|_{W^{3/2,\infty}(\mathbb{T})}),\\
\| (J_{\epsilon})^{\ast} -J_{\epsilon} \|_{H^\mu(\mathbb{T}) \rightarrow H^{\mu+3/2}(\mathbb{T}) } &\leq C(\| \eta_x \|_{W^{3/2,\infty}(\mathbb{T})}),
\end{align*}
for some non-decreasing function $C$ independent of $\epsilon \in [0,1]$. Equivalently, we have
\begin{equation*}
J_{\epsilon} T_{\vartheta} \sim T_{\vartheta} J_{\epsilon}, \quad (J_{\epsilon})^{\ast} \sim J_{\epsilon},
\end{equation*}
uniformly in $\epsilon$. Next, recalling Definition \ref{def:defiSigma}, we look for
\begin{align*}
\mathtt{P} &= \mathtt{P}^{(-1/2)} + \mathtt{P}^{(-3/2)}\in \Sigma^{-1/2}
\end{align*}
such that 
\begin{align*}
p \sharp \mathtt{P} &= p^{(1/2)} \mathtt{P}^{(-1/2)} + p^{(1/2)} \mathtt{P}^{(-3/2)} + p^{(-1/2)} \mathtt{P}^{(-1/2)} -\mathrm{i} \partial_\xi p^{(1/2)} \cdot \partial_x \mathtt{P}^{(-1/2)} =1.
\end{align*}
In order to solve this equation we set
\begin{equation*}
\mathtt{P}^{(-1/2)} = \frac{1}{p^{(1/2)}}, \quad \mathtt{P}^{(-3/2)} = -\frac{1}{p^{(1/2)}} \left( \mathtt{P}^{(-1/2)}p^{(-1/2)} -\mathrm{i} \partial_\xi \mathtt{P}^{(-1/2)}  \partial_x p^{(1/2)} \right) ,
\end{equation*}
from which we deduce that $T_p T_{\mathtt{P}} \sim \mathrm{Id}$, where recall that the notation $A \sim B$ introduced in Definition \ref{def:defisym}. On the other hand, since by Remark \ref{rem:key} we have that $q$ does not depend on $\xi$, it follows by Lemma \ref{lem:Paraprod3} that $T_q T_{1/q} \sim \mathrm{Id}$. Hence,
\begin{align*}
\begin{pmatrix} 
T_p & 0 \\ 
0 & T_q
\end{pmatrix} 
\begin{pmatrix} 
T_{\mathtt{P}} & 0 \\ 
0 & T_{1/q} 
\end{pmatrix}
&\sim
\begin{pmatrix} 
\mathrm{Id} & 0 \\ 
0 & \mathrm{Id}
\end{pmatrix}
,
\end{align*}
namely, we constructed a parametrix for the symmetrizer. Next, we define
\begin{align*}
\mathcal{L}_{\epsilon} &\coloneqq 
\begin{pmatrix}
\mathrm{Id} & 0 \\
T_{B(\eta,\beta,\gamma)(\psi)} & \mathrm{Id} 
\end{pmatrix}
\begin{pmatrix}
0 & -  T_{\lambda}   \\
\kappa T_{\mathtt{h}}  & 0 
\end{pmatrix}
\begin{pmatrix}
T_{\mathtt{P}}  J_{\epsilon} T_p & 0\\
0 & T_{1/q}J_{\epsilon} T_q   
\end{pmatrix}
\begin{pmatrix}
\mathrm{Id} & 0 \\
-T_{B(\eta,\beta,\gamma)(\psi)} & \mathrm{Id} 
\end{pmatrix}.
\end{align*}

We construct $(\eta,\psi)$ of \eqref{eq:ref} as limits of solutions of the Cauchy problems
\begin{equation} \label{eq:refeps} 
\begin{cases}
\left( \partial_t+   T_{V(\eta,\beta,\gamma)(\psi) - \gamma \, \eta} \partial_x J_{\epsilon}
+ \mathcal{L}_{\epsilon} \right)
\begin{pmatrix}
\eta \\ \psi
\end{pmatrix}
&= f( J_{\epsilon} \eta,J_{\epsilon} \psi),\\
(\eta,\psi)|_{t=0} &= (\eta_{0},\psi_0).
\end{cases}
\end{equation} 

\begin{proposition} \label{prop:apriori}

Let $s > s_0$. There exists a non-decreasing function $C$ such that, for all $\epsilon \in [0, 1]$, all $T\in ]0,1]$ and all solution $(\eta,\psi)$ of \eqref{eq:refeps} such that $(\eta,\psi)\in C^{1}( [0,T];H^{s+1/2}_0(\mathbb{T}) \times H^{s}(\mathbb{T}) )$, the quantity
\begin{align*}
M(T) &= \| (\eta,\psi) \|_{L^{\infty}( [0,T];H^{s+1/2}_0(\mathbb{T}) \times H^{s}(\mathbb{T}) )}
\end{align*}
satisfies the estimate
\begin{align*}
M(T) &\leq C(M_0) + T C(M(T)), \quad M_0 \coloneqq  \| (\eta_0,\psi_0) \|_{H^{s+1/2}_0(\mathbb{T}) \times H^{s}(\mathbb{T}) } .
\end{align*}
\end{proposition}

The proof of Proposition \ref{prop:apriori} follows the argument of Sec. 5.5 in \cite{alazard2011water}. Observe that the estimates of Proposition \ref{prop:apriori} holds for $\epsilon = 0$; hence, the above result proves also a priori estimates for the water waves system \eqref{eq:WWsys}. 

Now consider a solution $(\eta,\psi)\in C^{0}([0,T];H^{s+1/2}(\mathbb{T}) \times H^{s}(\mathbb{T}))$ to the system \eqref{eq:refeps}. We want to state uniform estimates for solutions $(\tilde{\eta},\tilde{\psi})$ to the linear system  
\begin{equation} \label{eq:refLin}
\begin{cases}
\left( \partial_t +   T_{V(\eta,\beta,\gamma)(\psi) - \gamma \, \eta} \partial_x J_{\epsilon}
+ \mathcal{L}_{\epsilon} \right)
\begin{pmatrix}
\tilde{\eta} \\ \tilde{\psi}
\end{pmatrix}
&= F,\\
(\tilde{\eta},\tilde{\psi})|_{t=0} &= (\tilde{\eta}_{0},\tilde{\psi}_0).
\end{cases}
\end{equation}

We rewrite the systems \eqref{eq:refeps} and \eqref{eq:refLin}, respectively, in the form
\begin{align*}
E(\epsilon,\eta,\psi)
\begin{pmatrix}
\eta \\ \psi
\end{pmatrix}
&= 
f (J_{\epsilon} \eta,J_{\epsilon} \psi),  \qquad \text{and} \qquad
E(\epsilon,\eta,\psi)
\begin{pmatrix}
\tilde{\eta} \\ \tilde{\psi}
\end{pmatrix}
= F.
\end{align*}

\begin{proposition} \label{prop:apriori2}

Let $s > s_0$ and $0\leq \sigma \leq s$. 
There exists a non-decreasing function $C$ such that, for all $\epsilon \in [0,1]$, all $T\in ]0,1]$ and all $\tilde{\eta},\tilde{\psi},\eta,\psi,F$ such that
\begin{align*}
E(\epsilon,\eta,\psi)
\begin{pmatrix}
\eta \\ \psi
\end{pmatrix}
&= f (J_{\epsilon} \eta,J_{\epsilon} \psi), \qquad \text{and} \qquad
E(\epsilon,\eta,\psi)
\begin{pmatrix}
\tilde{\eta} \\ \tilde{\psi}
\end{pmatrix}
= F,
\end{align*}
and such that
\begin{align*}
(\eta,\psi) &\in C^{0}([0,T];H^{s+1/2}_0(\mathbb{T}) \times H^{s}(\mathbb{T}) ),\\
(\tilde{\eta},\tilde{\psi}) &\in C^{1}([0,T];H^{\sigma+1/2}_0(\mathbb{T}) \times H^{\sigma}(\mathbb{T}) ), \\
F &=(F_1,F_2)\in L^{\infty}([0,T];H^{\sigma+1/2}(\mathbb{T}) \times H^{\sigma}(\mathbb{T}) ),
\end{align*}
we have
\begin{align} 
&\| (\tilde{\eta},\tilde{\psi}) \|_{ L^{\infty}([0,T];H^{\sigma+1/2}_0(\mathbb{T}) \times H^{\sigma}(\mathbb{T}) ) } \nonumber \\
&\leq
\widetilde{C}
\| (\tilde{\eta}_0,\tilde{\psi}_0) \|_{H^{\sigma+1/2}_0(\mathbb{T}) \times H^\sigma(\mathbb{T}) } \nonumber \\
&\quad +T C\left( \| (\eta,\psi) \|_{L^{\infty}([0,T];H^{s+1/2}_0(\mathbb{T}) \times H^{s}(\mathbb{T}) )}\right)
\| (\tilde{\eta},\tilde{\psi}) \|_{L^{\infty}([0,T];H^{\sigma+1/2}_0(\mathbb{T}) \times H^{\sigma}(\mathbb{T}) )} \nonumber \\
&\quad +T \| F \|_{L^{\infty}([0,T];H^{\sigma+1/2}(\mathbb{T}) \times H^{\sigma}(\mathbb{T}) )}, \label{eq:estuni2}
\end{align}
where 
\begin{align*}
\widetilde{C} &\coloneqq C \left(\| (\eta_0,\psi_0) \|_{H^{s+1/2}_0(\mathbb{T}) \times H^s(\mathbb{T}) } \right)+ T C \left( \| (\eta,\psi) \|_{ H^{s+1/2}_0(\mathbb{T}) \times H^s(\mathbb{T}) } \right) .
\end{align*}
\end{proposition}

The proof of Proposition \ref{prop:apriori2} follows the one of Proposition 5.4 in \cite{alazard2011water}. By applying Proposition \ref{prop:apriori2} with $(\eta,\psi)=(\tilde{\eta},\tilde{\psi})$ we obtain Proposition \ref{prop:apriori}.

\subsection{Proof of Theorem \ref{thm:LWPMain}} \label{subsec:LWPproof}

In this section we prove the local well-posedness result Theorem \ref{thm:LWPMain}.

We first prove the existence of solutions. First, we state that for any $\epsilon >0$, the approximate systems \eqref{eq:refLin} are well-posed locally in time.

\begin{lemma} \label{lem:TechLem1}

For all $(\eta_0,\psi_0)\in H^{s+1/2}(\mathbb{T}) \times H^s(\mathbb{T})$, and any $\epsilon \in ]0,1]$, there exists $T_{\epsilon}>0$ such that the Cauchy problem \eqref{eq:refeps} has a unique maximal solution $(\eta_{\epsilon},\psi_{\epsilon} )\in C^{0}([0,T_{\epsilon}[;H^{s+1/2}(\mathbb{T}) \times H^{s}(\mathbb{T}) )$.
\end{lemma}

The proof of Lemma \ref{lem:TechLem1} follows by observing that the operator $J_{\epsilon}$ is smoothing and by a standard well-posedness argument for ODEs (see Lemma 6.1 in \cite{alazard2011water}).

Next, we notice that by using the a priori estimates of Proposition \ref{prop:apriori} that the solutions $(\eta_{\epsilon},\psi_{\epsilon})$ of the approximate system \eqref{eq:refLin} are uniformly bounded with respect to $\epsilon$. 

\begin{lemma} \label{lem:TechLem2}

There exists $T_0>0$ such that $T_{\epsilon} \geq T_0$ for all $\epsilon \in ]0,1]$ and such that 
$\{ (\eta_{\epsilon},\psi_{\epsilon}) \}_{\epsilon \in ]0,1]}$ is bounded in 
$C^{0}([0,T_0]; H^{s+1/2}(\mathbb{T}) \times H^{s}(\mathbb{T}) )$. 
\end{lemma}

For the proof of the above lemma see the proof of Lemma 6.2 in \cite{alazard2011water}. Due to the above result, we can deduce that the sequence $\{ (\eta_{\epsilon},\psi_{\epsilon}) \}_{\epsilon \in ]0,1]}$ converges.

Next, we show that the solutions $(\eta_{\epsilon},\psi_{\epsilon})_{0 < \epsilon \leq 1}$ form a Cauchy sequence in a suitable space. In order to prove this fact, we need the following technical lemma.

\begin{lemma} \label{lem:TechLem30}
Let $0 < \epsilon_1 < \epsilon_2$, consider $s^{\prime}$ such that
\begin{equation*}
1 < s^{\prime} < s-\frac{3}{2} ,
\end{equation*}
and set $a=s-s^{\prime}-\frac{3}{2}$.

Then $\delta\eta \coloneqq \eta_{\epsilon_1}-\eta_{\epsilon_2}$ and $\delta\psi \coloneqq \psi_{\epsilon_1}-\psi_{\epsilon_2}$ satisfy a system of the form
\begin{equation}\label{eq:deltaEq}
\left[ \partial_t+ T_{V_{\epsilon_1}(\eta,\beta,\gamma)(\psi) - \gamma \, \eta} \partial_x J_{\epsilon_1} + \mathcal{L}_{\epsilon_1} \right] 
\begin{pmatrix}
\delta\eta \\ \delta\psi
\end{pmatrix}
= 
\begin{pmatrix}
r_1 \\ r_2
\end{pmatrix}
,
\end{equation}
where
\begin{align*}
&\| (r_1,r_2) \|_{L^{\infty}([0,T]; H^{s^{\prime}+1/2}(\mathbb{T}) \times H^{s^{\prime}}(\mathbb{T}) )} \\
&\leq C \left[ \| (\delta\eta,\delta\psi) \|_{L^{\infty}([0,T]; H^{s^{\prime}+1/2}(\mathbb{T}) \times H^{s^{\prime}}(\mathbb{T}) ) } + (\epsilon_2-\epsilon_1)^a \right],
\end{align*}
for some constant $C$ depending only on $\sup_{0 < \epsilon \leq 1} \| (\eta_{\epsilon},\psi_{\epsilon}) \|_{L^{\infty}([0,T]; H^{s+1/2}(\mathbb{T}) \times H^{s}(\mathbb{T}) ) }$.
\end{lemma}
 
The proof of Lemma \ref{lem:TechLem30} follows along the lines of the proof of Lemma 6.3 in \cite{alazard2011water}.

\begin{lemma} \label{lem:TechLem3}

Let $s^{\prime} < s- \frac{3}{2}$. Then there exists $0 < T_1 \leq T_0$ such that $\{ (\eta_{\epsilon},\psi_{\epsilon}) \}_{0 < \epsilon \leq 1}$ is a Cauchy sequence in 
$C^{0}([0,T_1]; H^{s^{\prime}+1/2}(\mathbb{T}) \times H^{s^{\prime}}(\mathbb{T}) )$. 
\end{lemma}

\begin{proof}

Let $0< \epsilon_1 <\epsilon_2$ and consider two solutions 
$(\eta_{\epsilon_j},\psi_{\epsilon_j})\in C^{0}([0,T];H^{s+1/2}(\mathbb{T}) \times H^{s}(\mathbb{T}) )$ of \eqref{eq:refeps}.  Let
\begin{equation*}
B_{\epsilon_j} = \frac{\partial_x\eta_{\epsilon_j} \, \partial_x\psi_{\epsilon_j} + G(\eta_{\epsilon_j},\beta,\gamma)(\psi_{\epsilon_j}) }{ 1+ (\partial_x\eta_{\epsilon_j})^2 }, \quad V_{\epsilon_j} =\partial_x\psi_{\epsilon_j} - B_{\epsilon_j}\partial_x\eta_{\epsilon_j},
\end{equation*}
and denote by $\lambda_j, \mathtt{h}_j$ the symbols obtained by replacing $\eta$ by $\eta_{\epsilon_j}$ in \eqref{eq:lambdaOp} and \eqref{eq:hSymbol}, respectively. Since for $t=0$ we have $\delta\eta=0=\delta\psi$, it follows from  Lemma \ref{lem:TechLem3} and from Proposition \ref{prop:apriori2} applied with $\sigma=s^{\prime}$, $\epsilon=\epsilon_1$, $\eta=\delta\eta$ and $\psi=\delta\psi$, that 
\begin{align*}
N &\leq T C \left[ N + (\epsilon_2-\epsilon_1)^a \right] ,
\end{align*}
so that by choosing $T$  and $\epsilon_2$ small enough, this implies $N = \mathcal{O}((\epsilon_2-\epsilon_1)^a)$, which proves the thesis.
\end{proof}

In order to complete the proof of Theorem \ref{thm:LWPMain}, we have to prove the uniqueness. 

\begin{proposition} \label{prop:Uniq}
Let $(\eta_j,\psi_j)\in C^0([0,T]; H^{s+1/2}(\mathbb{T}) \times H^s(\mathbb{T}) )$, $j=1,2$, be two solutions of system \eqref{eq:WWsys} with the same initial data, and such that the assumption $\mathrm{(A1)}$ is satisfied for all $t\in [0,T]$. Then $(\eta_1,\psi_1)=(\eta_2,\psi_2)$. 
\end{proposition}

We split the proof Proposition \ref{prop:Uniq} into few intermediate results. Recall that $(\eta,\psi)$ solves \eqref{eq:WWsys} if and only if they solve the system \eqref{eq:ref}, where the right-hand side $f=(f_1, f_2)^T$ is given by \eqref{eq:f1f2}. Let us introduce the notation
\begin{equation} \label{eq:BjVj}
B_j = \frac{\partial_x\eta_j \partial_x\psi_j + G(\eta_j,\beta,\gamma)(\psi_j)}{1+ (\partial_x\eta_j)^2}, \quad V_j = \partial_x\psi_j - B_{j}\partial_x\eta_j,
\end{equation}
and denote by $\lambda_j, \mathtt{h}_j$ the symbols obtained by replacing $\eta$ by $\eta_j$ in \eqref{eq:lambdaOp}, \eqref{eq:hSymbol} respectively. Similarly, we denote by $\mathcal{L}_j$ the operator obtained by replacing $B(\eta,\beta,\gamma)(\psi)$, $\lambda$ and $\mathtt{h}$ in \eqref{eq:ref} by $B_1$, $\lambda_1$ and $\mathtt{h}_1$.

In the following we denote $\delta\eta \coloneqq \eta_1 - \eta_2$, $\delta\psi \coloneqq \psi_1 - \psi_2$; moreover, we also write $M_j \coloneqq \| (\eta_j,\psi_j) \|_{ H^{s+1/2}(\mathbb{T}) \times H^s(\mathbb{T}) }$.

\begin{lemma} \label{lem:UnTechLem1}
We have
\begin{align*}
\| V_1 - V_2 \|_{ H^{s-5/2}(\mathbb{T}) } &\leq C(M_1,M_2,|\gamma|,\|\beta\|_{H^{s+1/2}(\mathbb{T})}) \, \|  (\delta\eta,\delta\psi) \|_{H^{s-1}(\mathbb{T}) \times H^{s-3/2}(\mathbb{T}) },\\
\| B_1 - B_2 \|_{ H^{s-5/2}(\mathbb{T}) } &\leq C(M_1,M_2,|\gamma|,\|\beta\|_{H^{s+1/2}(\mathbb{T})}) \, \|  (\delta\eta,\delta\psi) \|_{H^{s-1}(\mathbb{T}) \times H^{s-3/2}(\mathbb{T}) },
\end{align*}
\begin{align*}
& \sup_{|\xi|=1} \left[ |\partial_\xi^\alpha \big(\lambda_1^{(1)}(\cdot,\xi) -\lambda_{2}^{(1)}(\cdot,\xi) \big) | 
+ \| \partial_\xi^\alpha \big(\lambda_1^{(0)}(\cdot,\xi) -\lambda_{2}^{(0)}(\cdot,\xi)\big) \|_{ H^{s-3}(\mathbb{T}) } \right] \\
&\leq C(M_1,M_2,|\gamma|,\|\beta\|_{H^{s+1/2}(\mathbb{T})},\alpha) \, \| \delta\eta\|_{ H^{s-1}(\mathbb{T}) }, \; \; \forall \alpha \in \mathbb{N}, \\
&\sup _{|\xi|=1} \left[ |\partial_\xi^\alpha \big( \mathtt{h}_1^{(2)}(\cdot,\xi) - \mathtt{h}_2^{(2)}(\cdot,\xi) \big) | 
+ \| \partial_\xi^\alpha \big( \mathtt{h}_1^{(1)}(\cdot,\xi) - \mathtt{h}_2^{(1)}(\cdot,\xi)\big) \|_{ H^{s-3}(\mathbb{T}) } \right] \\
&\leq C(M_1,M_2,|\gamma|,\alpha) \, \| \delta\eta\|_{ H^{s-1}(\mathbb{T}) }, \; \; \forall \alpha \in \mathbb{N} .
\end{align*}

\end{lemma}

\begin{proof}

In order to prove the first two estimates, by \eqref{eq:BjVj} the non trivial point is to prove that
\begin{align*}
& \| G(\eta_1,\beta,\gamma)(\psi_1) - G(\eta_2,\beta,\gamma)(\psi_2) \|_{ H^{s-5/2}(\mathbb{T}) } \\
&\leq C(M_1,M_2,|\gamma|,\|\beta\|_{H^{s+1/2}(\mathbb{T})}) \, \|  (\delta\eta,\delta\psi) \|_{H^{s-1}(\mathbb{T}) \times H^{s-3/2}(\mathbb{T}) } .
\end{align*}
From \eqref{eq:OpG} we have 
\begin{align*}
G(\eta_1,\beta,\gamma)(\psi_1) - G(\eta_2,\beta,\gamma)(\psi_2) &= G^{DN}(\eta_1,\beta)\psi_1 - G^{DN}(\eta_2,\beta)\psi_2 \\
&\quad +\gamma \, \left[ G^{NN}(\eta_1,\beta) - G^{NN}(\eta_2,\beta) \right] ( (-h+\beta)\beta_x ) .
\end{align*}
Setting $\check{\eta}(t) \coloneqq t \eta_1 + (1 - t) \eta_2$, we have
\begin{align}
& G(\eta_1,\beta,\gamma)(\psi_1) - G(\eta_2,\beta,\gamma)(\psi_2)   \nonumber \\
&= G^{DN}(\eta_1,\beta)\delta\psi + \int_0^1 \mathrm{d}_{\eta} G^{DN}(\check{\eta},\beta)(\delta\eta)\psi_2  + \partial_{\eta} G^{NN}(\check{\eta},\beta)(\delta\eta)( -\gamma (-h+\beta)\beta_x ) \, \mathrm{d}t ,\label{eq:EstDiffOpG}
\end{align}
where
\begin{align*}
&\| G^{DN}(\eta_1,\beta)\delta\psi \|_{ H^{s-5/2}(\mathbb{T}) } \\
&\leq C(M_1,\|\beta\|_{H^{s+1/2}(\mathbb{T})}) \, \| \delta\psi \|_{ H^{s-3/2}(\mathbb{T}) } ,
\end{align*}
and where by Proposition \ref{prop:SysOpDer} i. we have that the $H^{s-5/2}$-norm of the integral in \eqref{eq:EstDiffOpG} is bounded by
\begin{align*}
& \int_0^1 \left\|  G^{DN}(\check{\eta},\beta) \left[ \delta\eta \, \underline{W}(\check{\eta},\beta)(\psi_2, \gamma (-h+\beta)\beta_x )  \right] \right\|_{ H^{s-5/2}(\mathbb{T}) } \mathrm{d}t \\
&+  \int_0^1 \left\| \partial_x  \left[ \delta\eta \, \underline{V}(\check{\eta},\beta)(\psi_2, \gamma (-h+\beta)\beta_x) \right] \right\|_{ H^{s-5/2}(\mathbb{T}) } \mathrm{d}t \\
&\leq C(M_1,M_2,|\gamma|,\|\beta\|_{H^{s+1/2}(\mathbb{T})}) \, \|  \delta\eta \|_{ H^{s-1}(\mathbb{T})  } ,
\end{align*}
which leads to the thesis.

The last two estimates are obtained by using the product rule in Sobolev spaces (see also Lemma 6.6 in \cite{alazard2011water}). 

\end{proof}

\begin{corollary} \label{cor:UnTechCor}
The following estimates hold true
\begin{align*}
\| T_{V_1-V_2} \eta_{2,x} \|_{ H^{s-1}(\mathbb{T}) } &\leq C(M_1,M_2,|\gamma|,\|\beta\|_{H^{s+1/2}(\mathbb{T})}) \, \|  (\delta\eta,\delta\psi) \|_{ H^{s-1}(\mathbb{T}) \times H^{s-3/2}(\mathbb{T}) },\\
\| T_{V_1-V_2} \psi_{2,x} \|_{ H^{s-3/2}(\mathbb{T}) } &\leq C(M_1,M_2,|\gamma|,\|\beta\|_{H^{s+1/2}(\mathbb{T})}) \,  \|  (\delta\eta,\delta\psi) \|_{ H^{s-1}(\mathbb{T}) \times H^{s-3/2}(\mathbb{T}) }, \\
\| T_{\lambda_1-\lambda_2}\psi_2 \|_{ H^{s-1}(\mathbb{T}) } &\leq C(M_1,M_2,|\gamma|,\|\beta\|_{H^{s+1/2}(\mathbb{T})}) \,  \|  (\delta\eta,\delta\psi) \|_{ H^{s-1}(\mathbb{T}) \times H^{s-3/2}(\mathbb{T}) },\\
\| T_{\mathtt{h}_1-\mathtt{h}_2}\eta_2 \|_{ H^{s-3/2}(\mathbb{T}) } &\leq C(M_1,M_2,|\gamma|,\|\beta\|_{H^{s+1/2}(\mathbb{T})}) \, \|  (\delta\eta,\delta\psi) \|_{ H^{s-1}(\mathbb{T}) \times H^{s-3/2}(\mathbb{T}) } .
\end{align*}
\end{corollary}

\begin{proof}
According to Proposition \ref{prop:OrdParadiffOp}, we have
\begin{align*}
\| T_a u \|_{H^\mu(\mathbb{T})} &\leq C \| a \|_{ L^{2}(\mathbb{T}) } \| u \|_{ H^{\mu+1/2}(\mathbb{T}) },
\end{align*}
which leads to the first two estimates. Using the estimates for $\lambda_1-\lambda_2$ and $\mathtt{h}_1-\mathtt{h}_2$ given by Lemma \ref{lem:UnTechLem1} we can deduce also the other two estimates.
\end{proof}

\begin{lemma} \label{lem:UnTechLem2}
The functions $\delta\eta$ and $\delta\psi$ satisfy a system of the form
\begin{align*}
\left( \partial_t+ T_{V_1} \, \partial_x + \mathcal{L}_1 \right) 
\begin{pmatrix}
\delta\eta \\ \delta\psi
\end{pmatrix}
&= f,
\end{align*}
where
\begin{align*}
\| f \|_{ L^\infty([0,T]; H^{s-1}(\mathbb{T}) \times H^{s-3/2}(\mathbb{T}) )} &\leq C (M_1,M_2,|\gamma|,\|\beta\|_{H^{s+1/2}(\mathbb{T})} ) N, \\
N &\coloneqq \| (\delta\eta,\delta\psi) \|_{L^\infty([0,T]; H^{s-1}(\mathbb{T}) \times H^{s-3/2}(\mathbb{T}) )}.
\end{align*}
\end{lemma}

\begin{proof}
Recall that for any $u \in H^{s+1/2}(\mathbb{T})$, 
\begin{align*}
\| T_{B_1-B_2} u \|_{ H^{s}(\mathbb{T}) } &\leq C \, \|  (\delta\eta,\delta\psi) \|_{ H^{s-1}(\mathbb{T}) \times H^{s-3/2}(\mathbb{T}) } \, \| u \|_{ H^{s+3/2}(\mathbb{T}) }.
\end{align*}
Now, we want to estimate the difference $f(\eta_1,\psi_1)-f(\eta_2,\psi_2)$, where $f(\eta,\psi)$ is defined in \eqref{eq:f1f2}. We prove an estimate for $f_1(\eta_1,\psi_1)-f_1(\eta_2,\psi_2)$, the other term being easier: arguing as in the proof of Lemma \ref{lem:UnTechLem1}, we have that by \eqref{eq:EstOpG}, by Lemma \ref{lem:EstGU}, by Proposition \ref{prop:OpGpara} and by Proposition \ref{prop:SysOpDer} i. 
\begin{align*}
&\| f_1(\eta_1,\psi_1) - f_1(\eta_2,\psi_2) \|_{ H^{s-1}(\mathbb{T}) } \\
&\leq C(M_1,M_2,|\gamma|,\|\beta\|_{H^{s+1/2}(\mathbb{T})} ) \, \| (\delta\eta,\delta\psi) \|_{ H^{s-1}(\mathbb{T}) \times H^{s-3/2}(\mathbb{T})} .
\end{align*}
\end{proof}

\begin{proof}[ Proof of Proposition \ref{prop:Uniq}: ]
Set
\begin{align*}
\delta \omega &\coloneqq \delta\psi - T_{B_1}\delta\eta, \qquad
\delta\Phi \coloneqq 
\begin{pmatrix} 
T_{p_1}\delta\eta \\ 
T_{q_1}\delta\omega 
\end{pmatrix}
,
\end{align*}
we obtain that $\delta\Phi$ solves a system of the form
\begin{align*}
\partial_{t}\delta\Phi + T_{V_1} \, \partial_x \delta\Phi  +
\begin{pmatrix} 
0 &  - T_{\vartheta_1} \\ 
T_{\vartheta_1} & 0
\end{pmatrix}
\delta\Phi  &= F ,
\end{align*}
with
\begin{align*}
\| F \|_{L^\infty([0,T]; H^{s-3/2}(\mathbb{T}) \times H^{s-3/2}(\mathbb{T}) ) } &\leq C(M_1,M_2,|\gamma|,\|\beta\|_{ H^{s+1/2}(\mathbb{T}) } ) N.
\end{align*}
Then by Proposition \ref{prop:apriori2} 
we have that $N$ satisfies and estimate of the form
\begin{align*}
N &\leq T C(M_1,M_2,|\gamma|,\|\beta\|_{ H^{s+1/2}(\mathbb{T}) }) N. 
\end{align*}
By choosing $T$ small enough this implies $N=0$, which leads to the thesis.
\end{proof}

\section*{Acknowledgements}

The author would like to thank Thomas Alazard for suggesting the setting of the problem, and  Emanuele Haus, Alberto Maspero, Federico Murgante and Erik Wahl\'en for useful comments and suggestions. The author would also like to thank the anonymous referee for suggesting improvements to the paper.

S. Pasquali is supported by PRIN 2022 ``Turbulent effects vs Stability in Equations from Oceanography" (TESEO), project number: 2022HSSYPN. S. Pasquali would like to thank INdAM-GNAMPA.

\begin{appendix}

\section{Paradifferential calculus} \label{sec:Paradiff}

Here we remind some rules of paradifferential calculus (see Appendix A in \cite{alazard2018control} and Sec. 4 of \cite{alazard2009paralinearization}; see also Sec. 5 of \cite{metivier2008paradifferential} and Sec. 3 of \cite{alazard2011water} for the non-periodic case). 

Observe that, due to \eqref{eq:chiSymm}, if $a$ and $u$ are real-valued functions, so is $T_a u$.

\begin{definition} \label{def:OrdOp}
Let $m \in \mathbb{R}$, then an operator $T$ is of order $m$ if, for all $s \in \mathbb{R}$, it is bounded from $H^{s}(\mathbb{T})$ to $H^{s-m}(\mathbb{T})$.
\end{definition}

\begin{proposition} \label{prop:OrdParadiffOp}
Let $m \in \mathbb{R}$ and let $a \in \Gamma^m_0(\mathbb{T})$, then $T_a$ is of order $m$. Moreover, for any $s \in \mathbb{R}$ there exists $K > 0$ such that
\begin{align*}
\| T_a \|_{L( H^s(\mathbb{T}) ,  H^{s-m}(\mathbb{T}) )} &\leq K \, M^m_0(a) .
\end{align*}
\end{proposition}

Recall the following properties
\begin{proposition} \label{prop:CompParadiffOp}
Let $0 <\varrho \leq 1$, $m_1,m_2 \in \mathbb{R}$ and let $a \in \Gamma^{m_1}_{\varrho}(\mathbb{T})$, $b \in \Gamma^{m_2}_{\varrho}(\mathbb{T})$, then $T_a T_b - T_{ab}$ is of order $m_1+m_2-\varrho$. Moreover, for any $s \in \mathbb{R}$ there exists $K > 0$ such that
\begin{align*}
\| T_a T_b - T_{ab} \|_{L( H^s(\mathbb{T}) ,  H^{s-m_1-m_2+\varrho}(\mathbb{T}) )} &\leq K \, M^{m_1}_{\varrho}(a) \, M^{m_2}_{\varrho}(b).
\end{align*}
\end{proposition}

\begin{proposition} \label{prop:AsympParadiffOp}
Let $\varrho >0$, $m_1,m_2 \in \mathbb{R}$ and let $a \in \Gamma^{m_1}_\varrho(\mathbb{T})$, $b \in \Gamma^{m_2}_\varrho(\mathbb{T})$. Set
\begin{align*}
a \sharp b (x,\xi) &:= \sum_{|\alpha| < \varrho} \frac{ (-\mathrm{i})^\alpha }{\alpha!} \; \; \pd_{\xi}^\alpha a(x,\xi) \; \; \pd_{x}^\alpha b(x,\xi) \in \sum_{j < \varrho} \Gamma^{m_1+m_2-j}_{\varrho-j}(\mathbb{T}),
\end{align*}
then $T_a T_b - T_{a \sharp b}$ is of order $\leq m_1+m_2-\varrho$. Moreover, for any $s \in \mathbb{R}$ there exists $K > 0$ such that
\begin{align*}
\| T_a T_b - T_{a \sharp b} \|_{L( H^s(\mathbb{T}) ,  H^{s-m_1-m_2+\varrho}(\mathbb{T}) )} &\leq K \, M^{m_1}_{\varrho}(a) \, M^{m_2}_{\varrho}(b).
\end{align*}
\end{proposition}

If $a=a(x)$ does not depend on the variable $\xi$, then the paradifferential operator $T_a$ is called \emph{paraproduct}. From Proposition \ref{prop:OrdParadiffOp} we have that if $b \in H^\beta(\mathbb{T})$ with $\beta > \frac{1}{2}$, then $T_b$ is of order $0$. Moreover, we recall the following properties

\begin{lemma} \label{lem:Paraprod1}
Let $m_1 \in \mathbb{R}$, $m_2 < 1/2$, $a \in H^{m_1}(\mathbb{T})$, $b \in H^{m_2}(\mathbb{T})$. Then $T_b a \in H^{m_1+m_2-1/2}(\mathbb{T})$.

Moreover, if $a \in L^\infty(\mathbb{T})$, then $T_a$ is an operator of order $\leq 0$, and there exists $C>0$ such that
\begin{align*}
\| T_a u \|_{ H^s(\mathbb{T}) } &\leq C \, \|a\|_{L^\infty(\mathbb{T})} \, \|  u \|_{ H^s(\mathbb{T}) } , \; \; \forall u \in H^s(\mathbb{T}), \; \; \forall s \in \mathbb{R}.
\end{align*}
\end{lemma}

\begin{lemma} \label{lem:Paraprod2}
Let $m > 1/2$, $a \in H^{m}(\mathbb{T})$, $F \in C^{\infty}(\mathbb{R})$. Then 
\begin{equation*}
F(a)-F(0)-T_{F'(a)} a \in H^{2m-1/2}(\mathbb{T}) . 
\end{equation*}

Moreover, let $m_1,m_2 \in \mathbb{R}$ be such that $m_1 + m_2 >0$, and let $a \in H^{m_1}(\mathbb{T})$, $b \in H^{m_2}(\mathbb{T})$. Then 
\begin{align*}
R(a,b) &:= ab - T_a b - T_b a \in H^{m_1 + m_2 - 1/2}(\mathbb{T}) ,
\end{align*}
and there exists $K>0$ such that
\begin{align*}
\| R(a,b) \|_{ H^{m_1 + m_2 - 1/2}(\mathbb{T}) } &\leq K \, \| a \|_{H^{m_1}(\mathbb{T})} \, \| b \|_{H^{m_2}(\mathbb{T})} .
\end{align*}

\end{lemma}

\begin{lemma} \label{lem:Paraprod3}
Let $m_1,m_2 > \frac{1}{2}$. If $a \in H^{m_1}(\mathbb{R})$, $b \in H^{m_2}(\mathbb{R})$, then $T_a T_b - T_{ab}$ is of order $ - \left( \min(m_1,m_2) - \frac{1}{2} \right)$.  Moreover, for any $s \in \mathbb{R}$ there exists $K > 0$ such that
\begin{align*}
\| T_a T_b - T_{a \sharp b} \|_{L( H^s(\mathbb{T}) ,  H^{s+\min(m_1,m_2) - \frac{1}{2}}(\mathbb{T}) )} &\leq K \, \|a\|_{H^{m_1}(\mathbb{T})} \,  \|b\|_{H^{m_2}(\mathbb{T})} .
\end{align*}
\end{lemma}

\begin{lemma} \label{lem:Paraprod4}
Let $m,s \in \mathbb{R}$ be such that $m+s>0$. If $\ell \in \mathbb{R}$ satisfies
\begin{align*}
\ell \leq m , &\;\; \ell < m+s-\frac{1}{2},
\end{align*}
then there exists $K > 0$ such that for all $a \in H^m(\mathbb{T})$ and all $u \in H^s(\mathbb{T})$ we have
\begin{align*}
\| a u - T_a u \|_{ H^{\ell}(\mathbb{T}) } &\leq K \, \|a\|_{H^{m}(\mathbb{T})} \,  \|u\|_{H^{s}(\mathbb{T})} .
\end{align*}
\end{lemma}

\section{On the Laplace equation with nonhomogeneous Neumann condition at the bottom} \label{sec:BVP}

In this section we consider the following boundary value problem with nonhomogeneous Neumann equation at the bottom 
\begin{equation} \label{eq:NonhomBVP}
\begin{cases}
\Delta \varphi = 0 \; \; &\text{in} \; \; D_{\eta,\beta}, \\
\varphi = \psi \; \; &\text{at} \; \; y=\eta, \\
\varphi_x \, \beta_x - \varphi_y = \theta  \; \; &\text{at} \; \; y=-h+\beta ,
\end{cases}
\end{equation}
where both $\psi:\mathbb{T} \to \mathbb{R}$ and $\theta:\mathbb{T} \to \mathbb{R}$ are given functions. Clearly, \eqref{eq:NonhomBVP} is a generalization of the boundary value problem \eqref{eq:EllBVP}, where $\theta = \gamma (-h+\beta)\beta_x$. 

The system \eqref{eq:NonhomBVP} has been studied by Iguchi \cite{iguchi2011mathematical} in the context of the description of tsunamis, generalizing some properties of the classical Dirichlet--Neumann operator for the homogeneous case $\theta=0$ (see also Appendix A.5 in \cite{lannes2013water}). We now state some key results, and we defer to Secc. 3--5 of \cite{iguchi2011mathematical} for their proofs.

Recalling the straightening diffeomorphism \eqref{eq:straight}, we rewrite \eqref{eq:NonhomBVP} into the following elliptic boundary value problem on $D_0$,
\begin{equation} \label{eq:FlatNonhomBVP}
\begin{cases}
\Delta^{\Sigma} \tilde{\varphi} = 0 \; \; &\text{in} \; \; D_0, \\
\tilde{\varphi} = \psi \; \; &\text{at} \; \; w=0, \\
\partial_x^{\Sigma}\tilde{\varphi} \, \beta_x - \partial_w^{\Sigma}\tilde{\varphi} = \theta  \; \; &\text{at} \; \; w=-h .
\end{cases}
\end{equation}





Regarding the existence of a variational solution of \eqref{eq:FlatNonhomBVP} see Proposition A.20 of \cite{lannes2013water}. Following Sec. 3 of \cite{iguchi2011mathematical}, we define the Dirichlet--Neumann operator $G^{DN}(\eta,\beta)$, the Neumann--Neumann operator $G^{NN}(\eta,\beta)$, the Dirichlet--Dirichlet operator $G^{DD}(\eta,\beta)$ and the Neumann--Dirichlet operator $G^{ND}(\eta,\beta)$ such that
\begin{equation} \label{eq:SysOpDef}
\begin{cases}
G^{DN}(\eta,\beta)\psi + G^{NN}(\eta,\beta)\theta &= \nabla\varphi \cdot \bmN |_{y=\eta(x)} \\
G^{DD}(\eta,\beta)\psi + G^{ND}(\eta,\beta)\theta &= \varphi |_{y=-h+\beta(x)} \\
\end{cases}
.
\end{equation}

Now we state some results about the operators defined in \eqref{eq:SysOpDef}.

\begin{proposition} \label{prop:Op00}

Let us consider the operators defined in \eqref{eq:SysOpDef}, then
\begin{align*}
G^{DN}(0,0) &= |\bmD| \, \tanh(h \, |\bmD|) , \quad G^{ND}(0,0) = |\bmD|^{-1} \, \tanh(h \, |\bmD|) , \\
G^{NN}(0,0) &= - G^{DD}(0,0) = - \frac{1}{ \cosh(h \, |\bmD|) } .
\end{align*}

\end{proposition}

\begin{proof}

Following the strategy of the proof of Lemma \ref{lem:SolFlatBVP}, we can prove that for $\eta=\beta=0$ the boundary value problem \eqref{eq:NonhomBVP} has the solution
\begin{align*}
\varphi(\cdot,y) &= - \frac{ \sinh(|\bmD| \, y) }{ |\bmD| \, \cosh(h \, |\bmD|) } \theta + \frac{ \cosh(|\bmD| (y+h)) }{ \cosh(h \, |\bmD|) } \psi ,
\end{align*}
and by \eqref{eq:SysOpDef} we can deduce the thesis.
\end{proof}

\begin{proposition} \label{prop:SysOpAdj}

Let $s > \frac{3}{2}$, and let us assume that $\eta, \beta \in H^s(\mathbb{T})$ satisfy \eqref{eq:StrConnected}. Then the operators $G^{DN}(\eta,\beta)$ and $G^{ND}(\eta,\beta)$ are symmetric in $L^2(\mathbb{T})$, while the adjoint of $G^{NN}(\eta,\beta)$ on $L^2(\mathbb{T})$ is given by $- G^{DD}(\eta,\beta)$. Namely, we have
\begin{align*}
\left( G^{DN}(\eta,\beta)\psi_1 , \psi_2 \right) &= \left( \psi_1 , G^{DN}(\eta,\beta)\psi_2 \right) , \qquad \forall \psi_1,\psi_2 \in H^1(\mathbb{T}) , \\
\left( G^{ND}(\eta,\beta)\theta_1 , \theta_2 \right) &= \left( \theta_1 , G^{ND}(\eta,\beta)\theta_2 \right) , \qquad \forall \theta_1,\theta_2 \in L^2(\mathbb{T}) , \\
\left( G^{NN}(\eta,\beta)\theta_1 , \psi_1 \right) &= - \left( \theta_1 , G^{DD}(\eta,\beta)\psi_1 \right) , \qquad \forall \psi_1 \in H^1(\mathbb{T}), \;\; \theta_1 \in L^2(\mathbb{T}) .
\end{align*}

\end{proposition}

For the proof of Proposition \ref{prop:SysOpAdj} see the proof of Proposition 3.2 in \cite{iguchi2011mathematical}.


\begin{proposition} \label{prop:SysOpDer}

Let $s > 5/2$. Assume that $\eta,\beta \in H^{s+1/2}(\mathbb{T})$ satisfy \eqref{eq:StrConnected}. Then 
\begin{enumerate}
\item[i.] the shape derivative of $G^{DN}(\eta,\beta)+ G^{NN}(\eta,\beta)$ with respect to free surface variations is given by
\begin{align*}
& \partial_{\eta} G^{DN}(\eta,\beta)(\delta\eta)\psi + \partial_{\eta} G^{NN}(\eta,\beta)(\delta\eta)\theta \\
&= - G^{DN}(\eta,\beta) \left[ \delta\eta \, \underline{W}(\eta,\beta)(\psi,\theta)  \right] - \partial_x  \left[ \delta\eta \, \underline{V}(\eta,\beta)(\psi,\theta) \right] ,
\end{align*}
where
\begin{align*}
\underline{V}(\eta,\beta)(\psi,\theta) &\coloneqq \psi_x - \underline{W}(\eta,\beta)(\psi,\theta) \, \eta_x , \\
\underline{W}(\eta,\beta)(\psi,\theta) &\coloneqq \frac{1}{1+\eta_x^2} \left[ G^{DN}(\eta,\beta)\psi + G^{NN}(\eta,\beta)\theta + \eta_x \psi_x  \right] ;
\end{align*}

\item[ii.] the shape derivative of $G^{DN}(\eta,\beta)+ G^{NN}(\eta,\beta)$ with respect to bottom variations is given by
\begin{align*}
& \partial_{\beta} G^{DN}(\eta,\beta)(\delta\beta)\psi + \partial_{\beta} G^{NN}(\eta,\beta)(\delta\beta)\theta \\
&= - G^{NN}(\eta,\beta)\left[ \partial_x \left( \delta\beta \, w(\eta,\beta)(\psi,\theta)  \right) \right] ,
\end{align*}
where
\begin{align*}
w(\eta,\beta)(\psi,\theta) &\coloneqq \partial_x \left( G^{DD}(\eta,\beta)\psi + G^{ND}(\eta,\beta)\theta  \right) - W_b(\eta,\beta)(\psi,\theta) \, \beta_x, \\
W_b(\eta,\beta)(\psi,\theta) &\coloneqq \frac{1}{1+\beta_x^2} \left[ - \theta + \beta_x \, \partial_x \left( G^{DD}(\eta,\beta)\psi + G^{ND}(\eta,\beta)\theta  \right) \right] ;
\end{align*}

\item[iii.] the shape derivative of $G^{DD}(\eta,\beta)+ G^{ND}(\eta,\beta)$ with respect to free surface variations is given by
\begin{align*}
\partial_{\eta} G^{DD}(\eta,\beta)(\delta\eta)\psi + \partial_{\eta} G^{ND}(\eta,\beta)(\delta\eta)\theta &= - G^{DD}(\eta,\beta) \left[ \delta\eta \, \underline{W}(\eta,\beta)(\psi,\theta)  \right] ;
\end{align*}

\item[iv.] the shape derivative of $G^{DD}(\eta,\beta)+ G^{ND}(\eta,\beta)$ with respect to bottom variations is given by
\begin{align*}
& \partial_{\beta} G^{DD}(\eta,\beta)(\delta\beta)\psi + \partial_{\beta} G^{ND}(\eta,\beta)(\delta\beta)\theta \\
&= \delta\beta \, W_b(\eta,\beta)(\psi,\theta) - G^{ND}(\eta,\beta) \left[ \partial_x \left( \delta\beta \, w(\eta,\beta)(\psi,\theta) \right)  \right] .
\end{align*}
\end{enumerate}

\end{proposition}

For the proof of Proposition \ref{prop:SysOpDer} see Theorem 3.4-Theorem 3.7 in \cite{iguchi2011mathematical}.

\begin{proposition} \label{prop:EstOp1}

Let $s > 5/2$, assume that 
\begin{align*}
\| (\eta,\beta) \|_{ H^{s+1}(\mathbb{T}) \times H^{s+1}(\mathbb{T}) } &\leq M ,
\end{align*}
and assume that the fluid domain $D_{\eta,\beta}$ is strictly connected (see \eqref{eq:StrConnected}). Then there exists $C=C(M,h,h_0,s) > 0$ such that
\begin{align*}
\| G^{DN}(\eta,\beta)\psi \|_{ H^s(\mathbb{T}) } &\leq C \, \| \psi \|_{ H^{s+1}(\mathbb{T}) } ,  \qquad \forall \psi \in H^{s+1}(\mathbb{T}).
\end{align*}

\end{proposition}

For the proof of Proposition \ref{prop:EstOp1} see the proof of Proposition 5.2 in \cite{iguchi2011mathematical}.

\begin{proposition} \label{prop:EstOp2}

Let $s > 5/2$, assume that 
\begin{align*}
\| (\eta,\beta) \|_{ H^{s}(\mathbb{T}) \times H^{s}(\mathbb{T}) } &\leq M ,
\end{align*}
and assume that the fluid domain $D_{\eta,\beta}$ is strictly connected. Then there exists $C=C(M,h,h_0,s) > 0$ such that
\begin{align*}
\| G^{NN}(\eta,\beta)\theta \|_{ H^s(\mathbb{T}) } &\leq C \, \| \theta \|_{ H^{s}(\mathbb{T}) } , \qquad \forall \theta \in H^s(\mathbb{T}) .
\end{align*}

\end{proposition}

For the proof of Proposition \ref{prop:EstOp2} see the proof of Proposition 5.11 in \cite{iguchi2011mathematical}.

\begin{proposition} \label{prop:EstOp3}

Let $s > 3$, assume that 
\begin{align*}
\| (\eta,\beta) \|_{ H^{s+1/2}(\mathbb{T}) \times H^{s+1/2}(\mathbb{T}) } &\leq M ,
\end{align*}
and assume that the fluid domain $D_{\eta,\beta}$ is strictly connected. Then there exists $C=C(M,h,h_0,s) > 0$ such that
\begin{align*}
\| G^{DD}(\eta,\beta)\psi \|_{ H^s(\mathbb{T}) } &\leq C \, \| \psi \|_{ H^{s}(\mathbb{T}) } , \qquad \forall \psi \in H^s(\mathbb{T}) .
\end{align*}

\end{proposition}

For the proof of Proposition \ref{prop:EstOp3} see the proof of Proposition 5.13 in \cite{iguchi2011mathematical}.

\end{appendix}

\bibliography{Water_Waves_LWP}
\bibliographystyle{alpha}

\end{document}